\documentclass[a5paper,11pt]{book}

\usepackage{setspace} 

\usepackage{amsmath} 
\usepackage{amssymb}
\usepackage{theorem}

\usepackage{xspace} 
\usepackage{cite} 

\usepackage{endnotes} 

\renewcommand{\notesname}{Примечания}   
\makeatletter
\renewcommand{\enoteheading}%
{\chapter*{\notesname
  \@mkboth{\notesname}{\notesname}
  \addcontentsline{toc}{chapter}{Примечания} }%
     \leavevmode\par\vskip-\baselineskip}

\renewcommand{\@makeenmark}{\hbox{$^{\@theenmark )}$}} 

\renewcommand{\enoteformat}{\rightskip\z@ \leftskip\z@ \parindent=1.8em
     \leavevmode\llap{\hbox{$^{\@theenmark )}$}}} 

\makeatother

\renewcommand{\thechapter}{\Roman{chapter}}

\renewcommand{\chaptermark}[1]{\markboth{\chaptername\ \thechapter. \ #1}{}}

\makeatletter
\renewcommand{\@makechapterhead}[1]{%
  \vspace*{20\p@}
  {\parindent \z@ \raggedright \normalfont
    \ifnum \c@secnumdepth >\m@ne
      \if@mainmatter
        \large \bfseries \@chapapp\space \thechapter
        \par\nobreak
        \vskip 10\p@
      \fi
    \fi
    \interlinepenalty\@M
    \LARGE \bfseries #1\par\nobreak
    \vskip 40\p@
  }}
\renewcommand{\@makeschapterhead}[1]{%
  \vspace*{20\p@}%
  {\parindent \z@ \raggedright
    \normalfont
    \interlinepenalty\@M
    \LARGE \bfseries  #1\par\nobreak
    \vskip 40\p@
  }}
\makeatother

\makeatletter
\def\@chapter[#1]#2{\ifnum \c@secnumdepth >\m@ne
                       \if@mainmatter
                         \refstepcounter{chapter}%
                         \typeout{\@chapapp\space\thechapter.}%
                         \addcontentsline{toc}{chapter}%
                                   {\protect \chaptername\ \thechapter.\ \ #1}
                       \else
                         \addcontentsline{toc}{chapter}{#1}%
                       \fi
                    \else
                      \addcontentsline{toc}{chapter}{#1}%
                    \fi
                    \chaptermark{#1}%
                    \addtocontents{lof}{\protect\addvspace{10\p@}}%
                    \addtocontents{lot}{\protect\addvspace{10\p@}}%
                    \if@twocolumn
                      \@topnewpage[\@makechapterhead{#2}]%
                    \else
                      \@makechapterhead{#2}%
                      \@afterheading
                    \fi}
\makeatother

\makeatletter
\def\@sect#1#2#3#4#5#6[#7]#8{%
  \ifnum #2>\c@secnumdepth
    \let\@svsec\@empty
  \else
    \refstepcounter{#1}%
    \protected@edef\@svsec{\@seccntformat{#1}\relax}%
  \fi
  \@tempskipa #5\relax
  \ifdim \@tempskipa>\z@
    \begingroup
      #6{%
        \@hangfrom{\hskip #3\relax\@svsec}%
          \interlinepenalty \@M #8\@@par}%
    \endgroup
    \csname #1mark\endcsname{#7}%
    \addcontentsline{toc}{#1}{%
      \ifnum #2>\c@secnumdepth \else
        \protect\numberline{\S \csname the#1\endcsname}
      \fi
      #7}%
  \else
    \def\@svsechd{%
      #6{\hskip #3\relax
      \@svsec #8}%
      \csname #1mark\endcsname{#7}%
      \addcontentsline{toc}{#1}{%
        \ifnum #2>\c@secnumdepth \else
          \protect\numberline{\csname the#1\endcsname}%
        \fi
        #7}}%
  \fi
  \@xsect{#5}}
\makeatother

\usepackage[cp1251]{inputenc}
\usepackage[russian,english]{babel}

\newtheorem{Pa}{Paper}[section]

\newtheorem{corollary}[Pa]{Следствие}   

\newtheorem{definition}[Pa]{Определение}  

\newtheorem{lemma}[Pa]{Лемма}  

\begin{document}

\selectlanguage{russian}

\renewcommand{\contentsname}{Содержание}
\renewcommand{\bibname}{Литература}

\thispagestyle{empty}

\begin{center}
{\large БЕЛОРУССКИЙ ГОСУДАРСТВЕННЫЙ УНИВЕРСИТЕТ}
\quad \\
\quad \\
\quad \\
\bigskip {\bf \LARGE Александр Киселев }\bigskip\\
\quad \\
\quad \\
\quad \\
{\bf \Huge Недостижимость}
\quad \\
\quad \\
{\bf \Huge и}
\quad \\
\quad \\
{\bf \Huge субнедостижимость}\\
\quad \\
\quad \\
{\Large В двух частях} \\
\smallskip
{\Large Часть I} \\
\quad \\
\quad \\

\quad \\
\quad \\
\quad \\
\quad \\
\quad \\
\quad \\
\quad \\
\quad \\
{Минск} \\
\smallskip
{``Издательский центр БГУ''} \\
\smallskip
{2011}

\end{center}

\newpage

\thispagestyle{empty}

\noindent УДК 510.227
\\

{\small \textbf{Киселев, А. А.} Недостижимость и
субнедостижимость. В 2 ч. Ч. 1 / Александр Киселев. -- Минск :
Изд. центр БГУ, 2011. -- \pageref{end} с. -- ISBN
978-985-476-931-8.}

{\footnotesize Данное издание представляет собой перевод с
английского языка монографии А. А. Киселева под тем же названием,
содержащей доказательство (в $ZF$) несуществования недостижимых
кардиналов, 1-е издание которой вышло в свет в 2000 г. Часть I
содержит аппарат субнедостижимых кардиналов и его основные
средства -- теории редуцированных формульных спектров и матриц,
теорию диссеминаторов и другие -- которые используются в этом
доказательстве и представлены теперь в их более прозрачной и
подробной форме. Большое внимание уделяется более глубокой
разработке и культивированию базовых идей, служащих основаниями
для основных конструкций и рассуждений.

Предназначено для специалистов по теории множеств и математической
логике, а также для преподавателей и студентов факультетов
математического профиля.

Библиогр.: 26 назв.
\\
\\
Перевод осуществлён по изданию:\textbf{ Kiselev, Alexander}.
Iaccessibility and Subinaccessibility. In 2 pt. Pt 1 / Alexander
Kiselev. -- 2nd ed., enrich. and improv. -- Minsk: Publ. center of
BSU, 2008.}
\begin{center}

 {Р\;е\;ц\;е\;н\;з\;е\;н\;т\;ы}
\\
профессор {\em П. П. Забрейко;}
\\
профессор {\em А. В. Лебедев}
\quad \\
\quad \\
\quad \\
\noindent Классификация математических тем (2000): \\
03E05, 03E15, 03E35, 03E55, 03E60

\end{center}

\quad \\
\quad \\
\quad \\

\noindent {\scriptsize {\bf ISBN 978-985-476-931-8 (ч. 1)} \hfill
$\copyright$  Киселев А. А., 2011}

\noindent {\scriptsize {\bf ISBN 978-985-476-596-9 (pt. I)} \hfill
$\copyright$  Alexander Kiselev, 2008}

\newpage
\chapter*{Благодарности}

Автор высказывает свои первые слова глубокой благодарности Хан\-не
 Калиендо за понимание значительности темы и за сердечную
воодушевляющую помощь в продвижении работы.

Особая признательность выражается проф. С. Р. Когаловскому,
который научил автора Теории Иерархий, и  проф. Акихиро Канамори
за их бесценную воодушевляющую поддержку, придавшую необходимый
импульс завершению работы.

Автор хотел бы также выразить особую благодарность проф. А. В.
Лебедеву и проф. П. П. Забрейко за многолетнюю поддержку его
работ; большой интеллектуальный и моральный  долг высказывается им
обоим за их практическую и духовную помощь.

Глубокая благодарность выражается  проф. А. В. Тузикову и Ю.
Прокопчуку, которые оказывали автору интенсивную экспертную помощь
в наборе текстов.

Проф. В. М. Романчак  оказывал автору материальную и моральную
помощь в продолжении самого сложного периода исследований темы и
автор адресует ему много слов глубокой благодарности. Большая
благодарность высказывается также Людмиле Лаптёнок, которая
осуществила много первоначальных трудных наборов предыдущих
текстов автора, подготовивших  эту работу.

Глубокая сердечная благодарность высказывается Надежде Забродиной
за многолетние воодушевление и поддержку, без которых эта работа
была бы значительно затруднена.

Автор благодарит всех других специалистов за их поддержку и
воодушевление.

\newpage
{} \thispagestyle{empty}

\newpage

\begin{spacing}{1.5}

\tableofcontents

\end{spacing}

\newpage
{} \thispagestyle{empty}

\newpage
\chapter*{Введение}
\addcontentsline{toc}{chapter}{Введение}
\markboth{Введение}{Введение}

\setcounter{equation}{0}

\quad \ Понятие недостижимости в его различных формах принадлежит
к наиболее важным концепциям  культурных традиций с древнейших
времён и состоит в представлении о {\it экстраординарно} огромном
феномене, который не может быть достигнут средствами
 ``менее мощными'', чем сам этот феномен.
\\
В этом понимании он представляется одним из базовых архетипов
человечества.
\\

В начале 20 века эта концепция получила своё достаточно адекватное
выражение в основаниях математики следующим образом.
Конфинальностью кардинала \ $k$ \ называется минимальный кардинал
\ $cf(k)$, \ который может быть отображён некоторой функцией \ $f$
\ в \ $k$ \ так, что $\sup rng(f)=k$; \ кардинал \ $k$ \
называется регулярным, если  \ $cf(k)=k$, \ иначе он называется
сингулярным; он называется слабо недостижимым, если он несчётный,
регулярный и для каждого \ $\alpha < k$ \ следующий после \
$\alpha$ \ кардинал снова \ $ < k $. Понятие сильно недостижимого
кардинала получается заменой последнего условия на более сильное:
\ $\forall \alpha < k \ \ 2^{\alpha}< k $.
\\

Впервые понятие слабо недостижимого кардинала было введёно
Феликсом Хаусдорфом ~\cite{Hausdorf} в 1908 году; Пауль
Мало~\cite{Mahlo1, Mahlo2, Mahlo3} исследовал более сильные
предельные образования, кардиналы Мало в 1911-1913 годах. Позже
появились сильно недостижимые кардиналы, введённые Вацлавом
Серпинским и Альфредом Тарским~\cite{Sierp} и Эрнстом
Цермело~\cite{Zermelo} в 1930годах. После этого начального периода
появилось много видов недостижимых кардиналов -- слабо компактные,
неописываемые, кардиналы Эрдёша, Йонсона, Рауботтома, кардиналы
Рамсея, измеримые, сильные, кардиналы Вудин, суперсильные, сильно
компактные, суперкомпактные, главные кардиналы Вопенки, почти
огромные, огромные, суперогромные кардиналы и т.д., образуя
иерархию, соответствующую существенному усилению ``степени
недостижимости'' -- все они получили общее название  {\it
``большие кардиналы''.}
\\

Таким образом в основаниях математики прочно установилась традиция
недостижимости, служащая обоснованием веры в то, что главное
направление развития современной теории множеств состоит в
формировании систем гипотез больших кардиналов, естественных
моделей таких гипотез и генерических расширений таких моделей. Во
всяком случае, утвердилась вера в то, что:

``Для всех гипотез больших кардиналов ... можно предложить более
или менее убедительные обоснования ... -- убедительные, во всяком
случае настолько, что, пожалуй, никто не ожидает опровержений этих
гипотез.

Существуют, однако, другие, более проблематичные гипотезы больших
кардиналов. Для них до сих пор нет убедительного обоснования.
Однако и эти гипотезы далеки от того, чтобы быть отвергнутыми, но
привели к построению чрезвычайно красивых математических теорий.''
(Кюнен~\cite{Kunen}).

В дополнение к очень красивому обоснованию непротиворечивости
теории множеств в целом, эта вера базируется на большом количестве
хорошо известных результатов, указывающих на близкие взаимные
связи (по относительной непротиворичивости) между Гипотезами
больших кардиналов, Аксиомой детерминированности, регулярными
свойствами континуальных множеств, инфинитарными комбинаториками,
инфинитарными языками и другими. Большой набор таких результатов
представлен в монографиях Фрэнка Дрейка~\cite{Drake} и Акихиро
Канамори~\cite{Kanamori}; последняя, замечательно исчерпывающая и
современная, содержит также выдающуюся демонстрацию развития идей
различных гипотез больших кардиналов.
\\
``Исследования гипотез больших кардиналов действительно составляет
основное направление развития современной теории множеств и они
признаются играющими ключевую роль в изучении определимых
континуальных множеств, в частности их лебеговой измеримости. Хотя
и сформированные на различных стадиях развития теории множеств и с
различными мотивациями, эти гипотезы оказались формирующими
\emph{линейную} иерархию, простирающуюся до противоречивых
усилений мотивирующих концепций. Все известные
теоретико-мно\-жест\-вен\-ные утверждения были вовлечены в эту
иерархию в терминах силы совместности, и получающаяся структура
импликаций обеспечивает замечательно богатую, детальную и
согласованную картину сильнейших утверждений математики, вложенных
в теорию множеств.'' (Канамори~\cite{Kanamori}).
\\

\noindent Данная работа составляет первую часть второго издания
работы~\cite{Kiselev8}, которая представляет доказательство
теоремы:
\\

\noindent {\bf Основная теорема } ($ZF$) \\
\quad \\
\hspace*{1cm} {\it Не существует слабо недостижимого кардинала}
\\
\\

Доказательство этой теоремы получается в результате использования
аппарата субнедостижимых кардиналов, который автор разрабатывал
начиная с 1976 года, однако некоторые его технические детали
разрабатывались начиная с 1970 года. Идея этого доказательства
возникла в 1984 году (хотя начальные приближения к ней были
предприняты автором начинная с 1973 года) и в 1996 году она
приняла свою настоящую форму. Основные технические средства этого
доказательства -- формульные спектры и субнедостижимые кардиналы
-- были разработаны в~\cite{Kiselev1}; редуцированные спектры и
матрицы -- в~\cite{Kiselev2, Kiselev3, Kiselev4}; диссеминаторы и
\ $\delta$-матрицы -- в~\cite{Kiselev4}; матричная информативность
была разработана автором в~\cite{Kiselev5}; матричные  функции --
в~\cite{Kiselev6}. Доказательство основной теоремы было
представлено впервые в~\cite{Kiselev6}, а его более прозрачный и
полный вариант -- в~\cite{Kiselev7}; более систематическое
изложение всего этого исследования было представлено
в~\cite{Kiselev8}.

Однако высказывалась критика, упрекающая эту
работу~\cite{Kiselev8} в том, что она педставляет собой теорию
слишком сложную и перегруженную технической стороной исследования
в ущерб явному развитию базовых идей, которое должно
предшествовать их технической реализации.

Эти обстоятельства делают желательной определённую предварительную
экспозицию всего материала, свободную от подобных упрёков.

Итак, настоящая работа представляет собой первую часть второго
издания~\cite{Kiselev8}, призванную их преодолеть.

Доказательство основной теоремы получает в этой работе свою более
прозрачную и упрощенную форму, хотя и вызванную некоторым
естественным улучшением в заключительном определении
8.2~\cite{Kiselev8} понятия  \ $\alpha$-функции, и понятие нулевой
характеристики подвергнуто здесь более детальному рассмотрению.

\noindent Эти обстоятельства делают возможным упрощение обсуждения
\ $\alpha$-функций и их свойств. В результате финальная часть
доказательства основной теоремы представляется в её более
естественной и прозрачной форме. Но в целом  техническая сторона
материала повторяет здесь~\cite{Kiselev8}, но более
систематическим образом.
\\

Особое внимание уделяется существенному предварительному развитию
и описанию идей всех основных конструкций и рассуждений.
\\
Следует подчеркнуть, что доказательство основной теоремы занимает
здесь только  \S~11 второй части этой работы. Предшествующие же
разделы посвящены базовой и специальной теориям субнедостижимости,
которые хотя и используются для разработки идей и техники этого
доказательства, однако имеют самостоятельное значение.
\\
 Предварительно делается обзор материала~\cite{Kiselev1, Kiselev2, Kiselev3,
Kiselev4, Kiselev5, Kiselev6, Kiselev7}, необходимого для этого
доказательства. Технически простые детали будут опускаться.
Сведения из теории множеств, понятия и символика, необходимые для
дальнейших рассуждений, общеприняты и находятся, например, в
замечательной работе Йеха~\cite{Jech}, которая обеспечивает
начальное изложение материала и гораздо более того, поэтому они
будут считаться известными и часто будут использоваться без
комментариев.
\\

План этой работы можно описать следующим образом.
\\
Для того, чтобы излагать идеи этой работы, требуются специальные
понятия и термины.
\\
С этой целью первая часть этой работы, глава~I, посвящена развитию
базовой теории, необходимой для введения требуемого для этого
материала.
\\
Специальная цель всех этих исследований делает возможным
ограничить их  самой узкой версией.
\\
Итак, предполагается существование некоторого недостижимого
кардинала \ $k$ \, который остаётся недостижимым в конструктивном
классе $L$ и все рассмотрения и рассуждения проводятся в модели
$L_k$ или в её генерических расширениях, вызывая заключительное
противоречие.
\\

В \S~1 представляется краткое изложение идеи доказательства
основной теоремы в её первом приближении, без вхождения в детали.
Однако строгое введение этой идеи требует последовательного
введения новых понятий и средств, необходимых на каждой стадии,
так как вся теория развивается как серия последовательных
приближений, вызванных некоторой недостаточностью этой идеи в её
предшествующих формах. Поэтому она получает свою более уточнённую
форму последовательно, вовлекая на
каждой стадии небходимую технику по мере необходимости.\\
Всякий другой способ представления этой идеи, то есть с самого
начала в её \textit{окончательной форме}, даже без детализации и
строгостей, был бы неестественным и немотивированным.
\\
Также \S~1 содержит некоторые хорошо известные классические
сведения.
\\

В \S~2 вводятся \textit{формульные спектры} (определения 2.3,
2.4). Говоря здесь в общих чертах, спектр данной формулы \
$\varphi$ \ некоторого уровня \ $n$ \ -- это функция, содержащая
{\it всю информацию\/} о её существенных истинностных свойствах во
{\it всех\/} генерических расширениях модели \ $L_k$ \
определённого вида посредством генерических ультрафильтров на \
$(\omega_0, k)$-алгебре $B$ Леви. Такой спектр имеет область
значений, состоящую из \textit{всех существенных} булевых значений
в $B$ этой формулы, то есть её \textit{булев спектр}, и область
определения, состоящую из соответствующих ординалов, называемых
\textit{ординалами скачка} этой формулы, обеспечивающих такие
значения, то есть её ординальный спектр.
\\

Среди всех спектров наиболее интересны так называемые
\textit{универсальные спектры} (определение 2.6), объединяющие
\textit{все} формульные ординальные спетры данного уровня \ $n$ \
в единый спектр (лемма~\ref{2.7.}).
\\
Эти понятия рассматриваются также в их релятивизированной форме; в
этом случае ограничивающий ординал \ $\alpha < k$ \ вызывает все
истинностные свойства формулы \ $\varphi$ \ во всех
рассматриваемых расширениях модели \ $L_k$, \ ограниченных
ординалом \ $\alpha$, \, потому \ $\alpha$ \ называется
\textit{носителем} соответствующего спектра этой формулы.
\\

После этого, имея в распоряжении формульные спектры, мы вводим в
\S~3 центральное понятие базовой теории -- понятие
\textit{субнедостижимости} (определение 3.1). Главное свойство
субнедостижимого ординала \ $\alpha$ \ уровня \ $n$ \ состоит в
\textit{невозможности достижения} этого ординала ``средствами
языка'' уровня \ $n$. Так как истинностные свойства формул
заключены в их спектрах, то это означает, что ординал \ $\alpha$ \
включает в себя все формульные спектры уровня \ $n$ \ с
константами \ $ < \alpha$ \ и, таким образом, не может быть
достигнут такими спектрами.
\\
Здесь можно увидеть определённые параллели с самим понятием
\textit{недостижимости}  (см. замечание после определения 3.1).

Некоторые простые свойства понятия субнедостижимости делают
возможным введение так называемых \textit{субнедостижимо
универсальных} формул и их спектров, обладающих только
субнедостижимыми ординалами скачка меньших уровней (определение
3.9). В последующей теории такие спектры доставляют очень удобные
технические средства. \\

В \S~4 возникает проблема ``сложности'' спектров, которую
естественно представлять как порядковый тип спектра. Почти
очевидно, что такая характеристика спектра должна неограниченно
возрастать вместе с увеличением его носителя до \ $k$. \ Поэтому
наиболее важные аспекты этой ситуации исследуются для
\textit{редуцированных спектров} (определение 4.1), полученных
редукцией их булевых значений к некоторому заданному кардиналу \
$\chi$, после чего такие спектры получают порядковые типы \ $<
\chi^{+}$.\
\\
Основной результат этого раздела, лемма 4.6 о спектральном типе,
показывает, что при некоторых естественных и креативных условиях
порядковый тип субнедостижимо универсального спектра,
редуцированного к \ $\chi$ \ на носителе \ $\alpha$, \ превосходит
всякий ординал \ $< \chi^{+}$, \ определимый ниже \ $\alpha$ \
(более точно, определимый ниже его некоторых ординалов скачка).
\\

Однако тем не менее такие спектры на различных носителях трудно
сравнивать друг с другом, потому что их области определения --
ординальные спектры -- состоят из ординалов, возрастающих
неограниченно, когда их носители возрастают до недостижимого
кардинала \ $k$.
\\
Поэтому в \S~5 мы обращаемся к \textit{редуцированным матрицам}
(определение 5.1), получаемых из редуцированных субнедостижимо
универсальных спектров изоморфной заменой их областей определение
на соответствующие ординалы. Термин ``матрица'' здесь вполне
уместен, так как возможно использование двухмерных, трёхмерных
спектров и матриц и т. д. для более тонкого анализа истинностных
свойств формул (именно на этом пути автор в течении долгого
времени проводил свои начальные исследования этой темы).
\\
В~1977 году автор установил основной результат этого раздела,
лемму 5.11 о матричной информативности, поволяющую видеть, что
такая матрица может содержать \textit{всю информацию} о
\textit{каждой} начальной части конструктивного универсума и,
более того, о \textit{всех} их генерических расширениях,
ограниченных её носителем (и снова, более точно, ограниченных
ординалами скачка соответствующего спектра на этом носителе).
Таким образом, эта информация сохраняется при переходе от одного
носителя этой матрицы к другому её носителю.
\\

Затем, после формирования необходимых понятий и терминов, главные
инструменты доказательства основной теоремы -- \textit{матричные
функции} -- вступают в действие. Такие матричные функции
представляют собой последовательности матриц, редуцированных к так
называемому \textit{полному кардиналу} \ $\chi^{\ast}$ \
(определение 5.4), точной верхней грани универсального
ординального спектра данного уровня. Сначала простейшая версия
такой функции (определения 5.7, 5.14) представляется в конце \S~5.
Затем эта версия уточняется посредством придания этим матрицам так
называемых \textit{диссеминаторов} (определение 6.1) --
специальных кардиналов, распространяющих информацию о нижних
уровнях универсума до их носителей (более точно, до ординалов
скачка соответствующего спектра на таких носителях).
\\
Поэтому \S~6 посвящён исследованию этого понятия диссеминаторов в
его самой бедной форме в соответствии с основной целью; в этом
разделе мы предсталяем некоторые методы получения матриц и
матричных функций, снабжённых такими диссеминаторами.
\\

На этом этапе завершается построение базовой теории. Здесь
естественно кратко представить последующую вторую часть этой
работы, её главы~2 и~3. Используя аппарат, разработанный в первой
части, начинает действовать специальная теория матричных \
\mbox{$\delta$-}функ\-ций (определения 7.1, 7.2), вырабатывающая
некоторые инструменты коррекции необходимых технических средств.

В \S~8 эти функции развиваются до \ $\alpha$-функций,
представляющих собой инструмент доказательства основной теоремы в
его окончательной форме.

В \S\S~9, 10 исследуются свойства \ $\alpha$-функций.

Эта информация обеспечивает доказательство основной теоремы,
представленное в \S~11. После этого в \S~12 излагаются некоторые
следствия этой теоремы и других хорошо известных результатов.

Следует также обратить особое внимание на примечания в конце всей
работы.

\newpage
\chapter{Базовая теория: субнедостижимость, формульные спектры и матрицы, диссеминаторы}

\markboth{\chaptername\ \thechapter. \ Базовая теория}%
{\chaptername\ \thechapter. \ Базовая теория} 


\section{Предварительные сведения}
\setcounter{equation}{0}

Мы докажем, что аксиоматическая система
\[
    ZF+\exists k \hspace{2mm} (k \; \mbox{\it \ это слабо недостижимый кардинал})
\]
противоречива.
\\
В дальнейшем все рассуждения проводятся в этой теории.
\\

Как уже было  отмечено выше, идея доказательства основной теоремы
довольно сложна и непрозрачна, поэтому желателен некоторый её
предварительный анализ.

Наилучший образ действий состоит в её введении в трёх
последовательных стадиях, уточняя её на каждой последующей стадии
посредством разработки соответствующих новых понятий и
инструментов, призванных преодолеть некоторые препятствующие
особенности её предшествующих версий.
\\
Таким образом она становится существенно более усложнённой на
каждой стадии и получает свою окончательную форму в  \S~8.

Краткое описание этой идеи состоит в следующем.

Впервые эта идея возникает в конце \S~5 и опирается на матричные
функции, которые представляют собой последовательности матриц,
редуцированных к определённому фиксированному кардиналу.
\\
Можно получить достаточно адекватное представление о такой
функции, если предварительно ознакомиться со следующими понятиями
в общих чертах: понятиями формульного спектра (определения 2.3,
2.4), универсального спектра (определение 2.6, лемма 2.7),
субнедостижимости (определение 3.1), субнедостижимо универсального
спектра (определение 3.9), редуцированных спектров и матриц
(определения 4.1, 5.1) и сингулярных матриц (определение 5.7).
\\
На этих основаниях простейшие матричные функции
\[
    S_{\chi f}=( S_{\chi \tau } )_{\tau }
\]
вводятся как последовательности таких сингулярных матриц
специального вида (определения 5.13, 5.14).
\\
Такая функция имеет область значений, состоящую из сингулярных
матриц, редуцированных к некоторому фиксированному кардиналу \
$\chi$ \ и минимальных в смысле гёделевой функции \ $Od$ \ на
соответствующих носителях; это задание матричной функции очевидно
вызывает её монотонность в том же смысле (лемма 5.17~~1)~).
Область определения такой функции состоит из ординалов и
конфинальна недостижимому кардиналу \ $k$ \ (лемма \ref{5.18.})
\\
Роль редуцирующего кардинала \ $\chi$ \ в дальнейшем играет
главным образом полный кардинал \ $\chi^{\ast}$ \ (определение
5.4).

Теперь возникает идея доказательства основной теоремы в её
\textit{начальной форме}:
\\
\\
\textit{Требуемое противоречие следует получить посредством
создания матричной функции, которая должна обладать несовместными
свойствами: она должна быть монотонна и в то же время должна быть
лишена этой монотонности.}
\\
\\
Однако на этой стадии непосредственному доказательству
немонотонности этой функции препятствует следующее явление:
свойства рассматриваемого универсума изменяются после его
ограничения носителями значений этой функции -- редуцированных
матриц \ $S_{\chi \tau}$ \ (это можно видеть из обсуждения в конце
\S~5 после леммы 5.18).

Именно для того, чтобы преодолеть это препятствие, начинается
процесс трансформации матричных функций.
\\
Вводятся специальные кардиналы -- \textit{диссеминаторы}
(определение 6.1, леммы 6.2, 6.3). Это понятие используется здесь
лишь для трансформации матричной функции в наших целях, хотя имеет
более широкое самостоятельное значение.
\\
С этой целью её значения, матрицы \ $S_{\chi \tau}$, \ снабжаются
определёнными диссеминаторами, продолжающие необходимые свойства
универсума с его нижних уровней до их носителей (точнее, до их
особых кардиналов предскачка, определения 6.1, 5.9, 6.9).
\\
После этого в последующей второй части настоящей работы
совершается второе приближение к идее доказательства основной
теоремы и матричная функция трансформируется в \ $\delta$\textit{-
функцию,} которая также определена на множестве ординалов,
конфинальном \ $k$ \ (определения 7.1, 7.2, лемма 7.6).
\\

Но в результате эта новая функция, наоборот, теряет свою
монотонность (это можно видеть из обсуждения этой новой ситуации в
конце \S~7 после леммы 7.7).
\\
Выход из этого нового затруднительного положения состоит в третьем
приближении к идее доказательства, то есть в трансформации этой
последней функции в её более сложную рекурсивную форму -- так
называемую \ $\alpha$\textit{-функцию} (определения 8.1-8.3),
которая также определена на множестве ординалов, конфинальном $k$
(лемма 8.9). Это рекурсивное определение сформировано таким
образом, что на каждом шагу случаи подтверждения монотонности
требуются в \textit{первую} очередь (так что они снабжаются
\textit{``единичной характеристикой''}), в то время как случаи
нарушения монотонности разрешаются во вторую очередь, за неимением
ничего лучшего (и они относятся к \textit{``нулевой
характеристике''}).
\\
Таким образом, приоритет в этой рекурсии принадлежит значениям
матричной функции единичной характеристики и именно поэтому
избегаются случаи нарушения её монотонности.

В результате \ $\alpha$-функция, наконец, доставляет требуемое
противоречие: она не может быть монотонной (теорема 1) и в то же
время она обладает свойством монотонности (теорема 2).
\\

Такова идея доказательства основной теоремы в общих чертах; более
детальное её описание вовлекло бы слишком много технических
деталей. Поэтому вместо этого более удобно последовать
ознакомлению с упомянутыми выше понятиями (до конца \S~6) без
вхождения в детали и обходясь, может быть, общими представлениями.
Используя другой подход (представление этой идеи с самого начала в
её окончательной форме) можно получить картину ещё более сложную,
чем изложенная выше, затемнённую или даже бессодержательную.
\\

После этого перейдём к реализации намеченной выше программы.
\\

Слабо недостижимые кардиналы становятся сильно недостижимыми в
конструктивном классе Гёделя \ $L$; \ напомним, что это класс
значений конструктивной функции Гёделя \ $F$, \ определённой на
классе всех ординалов. Каждое множество \ $a\in L$ \ получает свой
ординальный номер
\[
    Od ( a ) =\min \{ \alpha :F (\alpha  ) =a\}.
\]
Если \ $\alpha $ \ это ординал, то \ $L_{\alpha }$ \ обозначает
начальный сегмент
\[
    \{ a\in L:Od ( a )<\alpha \}
\]
этого класса. Начальная структура в дальнейших рассуждениях -- это
счётный начальный сегмент
\[
    \mathfrak{M}=(L_{\chi ^{0}},\; \in , \; =)
\]
класса \ $L $, \ служащий стандартной моделью теории
\[
    ZF+V=L+\exists k \hspace{2mm} (k \; \mbox{\it \ слабо недостижимый
    кардинал}).
\]
В действительности будет использована только конечная часть этой
системы, так как будут использоваться формулы некоторой
ограниченной длины, как это будет ясно из дальнейшего. Более того,
счётность этой модели используется только для некоторого
технического удобства (см. ниже) и можно вполне обойтись без неё,
действуя в булевом универсуме.

Везде далее \ $k$ \ -- это  \textit{минимальный недостижимый
кардинал в} \ $\mathfrak{M}$. Мы будем исследовать его
``изнутри'', рассматривая иерархию субнедостижимых кардиналов\ $ <
k $; последние ``недостижимы'' посредством формул определённого
элементарного языка. Для получения достаточно богатой иерархии
таких кардиналов следует соответственно использовать некоторую
достаточно богатую алгебру истинности \ $B$. Хорошо известно
(Крипке~\cite{Kripke}), что каждая булева алгебра вкладывается в
подходящую
 \ $ ( \omega
_{0},\mu  ) $-алгебру смещения и поэтому для наших целей в
качестве \ $B$ \ естественно использовать сумму семейства таких
алгебр мощности \ $k$, \ то есть \ $ ( \omega _{0},k ) $-алгебру
Леви \ $B$.

Для наших целей удобно использовать эту алгебру в следующем виде.
Пусть \ $P\in \mathfrak{M}$ \ это множество вынуждающих условий --
конечных функций \ $p\subset k\times k$ \ таких, что для каждого
предельного \ $\alpha <k$ \ и \ $n\in \omega_{0}$ \vspace{-6pt}
\[
    \alpha +n\in dom ( p ) \longrightarrow p (\alpha +n ) <\alpha ;
\]
также пусть \ $p(n)\leq n $ \ для \ $\alpha=0 $. На \ $P$ вводится
отношение \ $\leq $ \ частичного порядка по обратному включению:
\vspace{-6pt}
\[
    p_{1} \leq p_{2} \longleftrightarrow p_{2} \subseteq p_{1} .
\]
После этого \ $P$ \ плотно вкладывается в булеву алгебру \ $B\in
\mathfrak{M} $, состоящую из регулярных сечений \ $\subseteq P$ и
полную в \ $\mathfrak{M}$. \ На \ $B$ определяется отношение
частичного порядка \ $\leq $:
\[
    A_{1} \leq A_{2} \longleftrightarrow A_{1} \subseteq A_{2} ,
\]
а также булевы операции \ $\cdot$, $+$, $\prod$, $\sum $ \
(см.~\cite{Jech}) определяются на \ $B$. \ Каждое условие \
\mbox{$p\in P$} \ отождествляется с регулярным сечением
\[
    \left[ p\right] =\{ p_{1}\in P:p_{1}\leq p\}
\]
и поэтому \ $P$ \ изоморфно вкладывается в \ $B$. \ Далее
приводятся хорошо известные результаты Коэна~\cite{Cohen},
Леви~\cite{Levy} (см. также Йех~\cite{Jech}).

\begin{lemma} \label{1.1.} \quad \\
\hspace*{1em} Алгебра \ $B$ \ выполняет \ $k$-цепное условие, то
есть каждое множество \  $X\subseteq B$, \ $X\in \mathfrak{M}$,\
состоящее из попарно несовместных булевых значений, имеет мощность
\ $<k$ \ в \ $\mathfrak{M\;}$:
\[
    \left ( \forall A_{1},A_{2}\in X\quad A_{1}\cdot A_{2}=0 \right)
    \longrightarrow \left| X\right| <k.
\]
\end{lemma}

Следуя этой лемме можно рассматривать вместо значений \  $A\in B$
\ только множества
\[
    P_{A}=\{ p\in A:dom (p) \subseteq \chi \}
\]
где
    \[\chi =\min \{ \chi^{\prime }:\forall p\in A\quad p\left|
    \chi ^{\prime }\right. \leq A\}
\]
(здесь \ $p\left| \chi ^{\prime }\right.$ \ это ограничение \ $p$
\  на \  $\chi ^{\prime }$). Так как \  $A=\sum P_{A}$, \ мы будем
отождествлять \  $A$ \ и \ $P_{A}$, \ то есть мы будем всегда
рассматривать \ $P_{A}$ \ вместо самого \ $A$.
\\
Благодаря этому соглашению все булевы значения \ $A \in B$ \ это
множества в \ $L_k$, \ не классы, и это делает возможным все
последующие рассуждения в целом.
\\

Мы будем изучать иерархию субнедостижимых кардиналов посредством
булевых значений в \ $B$ \  некоторых утверждений об их свойствах,
то есть действуя внутри булевозначного универсума \
$\mathfrak{M}^{B}$. Счётность структуры \ $\mathfrak{M}$ \
желательна здесь только для сокращения и упрощения рассуждений,
используя её генерические расширения посредством \
$\mathfrak{M}$-генерических ультрафильтров на \ $B$. \ Можно
обойтись без этого, используя соответствующие рассуждения в
булевозначном универсуме \ $L^{B}$ \ (см., например,~\cite{Jech}).

Будет более удобно производить генерические расширения \
$\mathfrak{M}$ \ не посредством ультрафильтров, а посредством
функций. Именно, \  \mbox{$\mathfrak{M}$-генерической функцией}
или функцией Леви на \ $k $ \ будем называть функцию \ $l\in
{}^{k}k$ \ такую, что каждое множество \ $X \in \mathfrak{M}  ,
X\subseteq P$, плотное в \ $B$, содержит некоторое \ $p\subset l$.
Все функции этого вида будут обозначаться общим символом \ $l$. \
Очевидно, \ $\mathfrak{M}$-генерические ультрафильтры \ $G$ \ на \
$B$ \ и такие функции взаимнооднозначно определяют друг друга:
\[
    \quad l=\cup ( P\cap G) \quad , \quad G= \{ A\in B:\exists
    p \in P( p\leq A\wedge p\subset l) \}.
\]
В этом случае интерпретация \ $i_{G}$ \ универсума \
$\mathfrak{M}^{B}$ \  обозначается через \ $i_{l}$. \ Как обычно,
если \ $\underline{a}\in \mathfrak{M}^{B}$, $a\in
\mathfrak{M}\left[ l\right] $ \  и \ $i_{l}(\underline{a}) =a$, то
\ $\underline{a}$ \ называется этикеткой или именем \ $a$. \
Символом \ $\left\| \varphi \right\| $, \ как обычно, обозначается
булево значение утверждения \ $\varphi $ \ с константами из \
$\mathfrak{M}^{B}$  \  в алгебре \ $B$.
\\
Для некоторого удобства вводится отношение \ $\stackrel{\ast}{\in
}~$: \ для всяких \ $l\in {}^{k}k$, \ $A\in B$
\[
    l\stackrel{\ast }{\in }A\longleftrightarrow \exists p\in P
    ( p\subset l\wedge p\leq A ).
\]

\begin{lemma} \label{1.2.} \hfill {} \\
\hspace*{1em} Пусть \ $l$ \ -- это  \ $\mathfrak{M}$-генерическая
функция на \ $k$ \ и \ $\varphi  ( a_{1}...,a_{n} ) $ \ --
утверждение, содержащее константы \ $a_{1}...,a_{n}\in
\mathfrak{M}\left[ l \right] $ \ с именами \
$\underline{a}_{1}...,\underline{a}_{n}$, \ тогда
\[
    \mathfrak{M}[l] \vDash \varphi ( a_{1}...,a_{n} )
    \longleftrightarrow l\ \stackrel{\ast }{\in }\left\| \varphi
    ( \underline{a}_{1}...,\underline{a}_{n} ) \right\| .
\]
\end{lemma}
\begin{lemma} \label{1.3.} \hfill {} \\
\hspace*{1em} Пусть \ $l$ \ --  \
$\mathfrak{M}$-генерическая функция на $k$, \ тогда:\quad \\
\quad
\\ {\em 1)} \quad \quad \quad \quad \quad \quad $
\mathfrak{M}\left[ l\right]
\vDash ZF+V=L\left[ l\right] +GCH+k=\omega _{1} ; $ \\ \quad \\
{\em 2)} \quad Для каждого \ $\chi<k$ \ пусть \ $\chi_{1}=\chi$ \
если \ $\chi$ \ регулярный и \ $\chi_{1}=\left( \chi^{+} \right)^
\mathfrak{M} $ \ если \ $\chi$ \ сингулярный кардинал в \
$\mathfrak{M}$, тогда
\[
    \mathfrak{M}\left[ l|\chi \right]\models \forall \alpha < \chi_{1}
    ~ |\alpha |\leq \omega_{0} \wedge \forall \alpha \geq
    \chi_{1} ~ |\alpha |=|\alpha |^\mathfrak{M}.
\]
\end{lemma}

\vspace{-0pt}

\begin{lemma} \label{1.4.} \hfill {} \\
\hspace*{1em} Пусть
\[
    t\in \mathfrak{M}^{B}, \: dom(t)\subseteq \left\{ \check{a} :
    a\in \mathfrak{M} \right\} \:\: \mbox{\it \ и \ } \:\: |rng(t)|<k.
\]
и пусть \ $B_{t}$ \ это подалгебра \ $B$, \ порождённая \
$rng(t)$. \\
Тогда для всякого утверждения \ $\varphi (t) $
\[
    \parallel \varphi (t) \parallel \in B_{t}.
\]
\end{lemma}

Таковы предварительные сведения, необходимые чтобы приступить к
разработке теории субнедостижимости. \label{c1}
\endnote{
\ стр. \pageref{c1}. \ Здесь желательно прояснить вопрос,
используется ли гипотеза недостижимого кардинала  в этих
исследованиях существенным образом.
\\
Эта гипотеза использется здесь на протяжении всего исследования.
Аксиома множества подмножеств используется в \ $L_k$ \ почти на
каждом шагу, в частности, при использовании сколемовских функций
для преобразования формул и их спектров, а также в каждом случае
разработки некоторого нового понятия.
\\
Кроме того, недостижимость
кардинала \ $k$ \ весьма существенно используется в $k$-цепном
условии алгебры \ $B$ \ (лемма~1.1), без которого булевы значения,
спектры и матрицы -- это классы в \ $L_k$, не множества, и они не
могут сравниваться друг с другом требуемым образом. В результате
базовая теория и специальная теория становятся невозможными в этой
ситуации.
\\
\quad \\
} 

\newpage

\section{Спектры формул}
\setcounter{equation}{0}

Основной инструмент дальнейших рассуждений -- это понятие
формульного спектра. В этом разделе представлена базовая теория
таких спектров, содержащая обсуждение простейших спектральных
свойств. Можно ввести это понятие в более общей форме (для
произвольных булевой алгебры \  $B$ \ и частично упорядоченной
структуры \ $\mathfrak{M\ }$); однако здесь достаточно
использовать его самый бедный вариант.
\\

Среди имён из \  $\mathfrak{M}^{B}$ \  выделяются канонические
имена, которые дают аналогичные результаты при любой
интерпретации. Например, таковы имена \ $\check{a}$ \ множеств \ $
a\in \mathfrak{M}$, \ которые будут отождествлятся с этими
множествами:
\[
    dom ( \check{a} ) =\{ \check{b}:b\in a\}, \quad
    rng ( \check{a} ) =\{ 1\}.
\]
Каноническое имя всякой \ $\mathfrak{M}$-генерической функции \
$l$ \ на \  $k$ \ это функция
\[
    \underline{l}=\{ ( ( \alpha, \beta ) ,\{ ( \alpha ,\beta ) \}):
    \{( \alpha ,\beta ) \} \in P\}.
\]
Легко видеть, что всегда \quad $i_{l}( \check{a}) = a$, \ $i_{l}(
\underline{l}) =l$.

Вводится следующий элементарный язык \ $\mathcal{L}$ \ над
стандартной структурой
\[
    ( L_{k}\left[ l\right] ,\; \in, \; =,\; l ).
\]
Его алфавит состоит из обычных логических символов: кванторов \
$\forall $, $\exists $, \ связок \  $ \wedge $, $\vee $, $\neg $,
$\longrightarrow $, $\longleftrightarrow $, \ скобок $( \; , \;
)$, индивидных переменных \ $x,y,z$... \ (с индексами или без
них), всех имён из булевозначного универсума \ $L_{k}^{B} $ ,\
служащих индивидными константами, и символов \ $\in $,
$=$, $\underline{l}$.\\
При интерпретации этого языка в генерическом расширении \
$L_{k}\left[ l\right] $ \ переменные пробегают \ $L_{k}\left[
l\right] $, \ индивидные константы \  $a\in L_{k}^{B}$ \
обозначают \ $i_{l} ( a ) $ \  и константы \  $\in $, $=$,
$\underline{l}$ \ обозначают соответственно стандартное отношение
принадлежности, равенство и функцию \ $l$. \ Если это расширение
фиксировано, тогда константы и их интерпретации как обычно
отождествляются.

Формулы языка \  $\mathcal{L}$ \  определяются обычным рекурсивным
образом начиная с атомарных формул вида \ $ t_{1}=t_{2}$, $\quad
t_{1}\in t_{2}$, \ где \ $t_{1}$, $t_{2}$ \ любые термы, которые
рекурсивно образуются из индивидных переменных, констант \ $\in
L_{k}^{B}$ \ и \ $\underline{l}$ \ последовательными
суперпозициями. Однако следуя традиции некоторые удобные
обозначения, отношения и термы будут использоваться в записях
формул, если это не повлечёт трудностей. Например, кортеж \ $m$ \
индивидных переменных или констант \ $x_{1},...,x_{m}$  \
обозначается через \ $( x_{1},...,x_{m} ) $, или, кратко, через \
$\overrightarrow{a}$; отношение порядка \ $ x_{1}\in x_{2} \  $ на
множестве ординалов обозначается через \ $x_{1}<x_{2}$ \ и так
далее.

Далее формулы будут рассматриваться как формулы языка \
$\mathcal{L}$, \ которые будут обозачаться малыми буквами конца
греческого алфавита (если не указано иное).

Формула \ $\varphi $ \ с индивидными  свободными переменными и
константами, образующими кортеж \  $ \overrightarrow{a}$, \ будет
обозначаться через \ $\varphi ( \overrightarrow{a} ) $. \  Если в
дополнение к этому \ $\varphi $ \ содержит символ, на который
нужно обратить внимание, он тоже может быть указан особо;
например, запись \ $\varphi ( \overrightarrow{a},\underline{l} ) $
\ указывает, что \  $ \varphi $ \ содержит вхождение имени \
$\underline{l}$.

Как обычно, вхождения кванторов \ $\exists x,\forall x$ \ в
формулу \ $\varphi$ \ называются ограниченными термом \ $t$, \
если они имеют вид
\[
    \exists x~ ( x\in t\wedge ... )
    ,\forall x~ ( x\in t\longrightarrow ... );
\]
формула называется ограниченной, если все её кванторы ограничены
некоторыми термами; она называется пренексной формулой, если все
вхождения её \textit{неограниченных} кванторов расположены слева
от вхождения её остальных кванторов и связок; эта
последовательность её неограниченных кванторов называется её
\textit{кванторным префиксом}.
\\
Формулы языка \  $\mathcal{L}$  \  будут интерпретироваться в
расширениях \ $L_{k}[l] $ \ и поэтому мы будем называть формулы \
$\varphi,\psi $ \ \textit{эквивалентными} и писать \ $\varphi
\longleftrightarrow \psi $ \ \textit, если они эквивалентны в
теории   $ZFC^{-}+V=L[\underline{l}]$, где $ZFC^{-}$ это $ZFC$ без
аксиомы множества всех подмножеств.
\\
Кроме того, интерпретируя формулы \ $ \varphi, \psi $ \ в  \
$L_{k}$, \  мы будем называть их \textit{конструктивно
эквивалентными} и использовать прежнее обозначение \ $\varphi
\longleftrightarrow \psi $ \ \textit, если они эквивалентны в
теории  $ZFC+V=L$.
\\
Мы будем называть их \textit{генерическими эквивалентными}\ и
писать \ $\varphi \dashv \vdash \psi $, \ \textit если
\[
    \| \varphi \longleftrightarrow \psi \| =1
\]
для любых значений в \ $L_{k}^{B}$  \  их свободных переменных.
\\
Далее рассматриваеся эквивалентность  в
  $ZFC^{-}+V=L[\underline{l}]$, если контекст не указывает
на другой вариант; и в любом случае значение понятия
эквивалентности и символа \ ``$\longleftrightarrow$'' \ будет
всегда очевидно определяться контекстом.

\begin{definition}
\label{2.1.} \ \\
{\em 1)}\quad Класс все пренексных формул \ $\varphi
(\overrightarrow{a},\underline{l}) $ \  с кванторным префиксом,
состоящим из  \ $n$ \ максимальных блоков одноименных кванторов и
начинающимся с  \ $\exists $, \ а также всех эквивалентных им
формул обозначается через $\Sigma _{n}( \overrightarrow{a}) $.
\\
Двойственный класс обозначается через \ $\Pi _{n}(
\overrightarrow{a}) $, \ а их пересечение, класс  \ $\Sigma _{n}(
\overrightarrow{a} ) \cap \Pi _{n} ( \overrightarrow{a})$ \ -- \
через \ $\triangle_{n} \left( \vec{a}\right)$.
\\
В результате возникает элементарная иерархия Леви:
\[
    \{ \Sigma _{n}( \overrightarrow{a} ) ;\ \Pi _{n}
    ( \overrightarrow{a}) \} _{n\in \omega _{0}} .
\]
{\em 2)}\quad  Через \  $\Sigma _{n}^{\dashv \vdash }(
\overrightarrow{a}) $ \  обозначается класс все формул,
генерически эквивалентных формулам из \ $\Sigma _{n}(
\overrightarrow{a})$.
\\
Двойственный класс обозначается через \ $\Pi_{n}^{\dashv \vdash }(
\overrightarrow{a}) $, \ а класс \ $\Sigma _{n}^{\dashv \vdash }
(\overrightarrow{a}) \cap \Pi _{n}^{\dashv \vdash }
(\overrightarrow{a}) $ \ -- \ через \  $\Delta _{n}^{\dashv \vdash
}( \overrightarrow{a}) $.
\\
\quad \\
Символы \  $Q_{n} (\overrightarrow{a})$, $ Q_{n}^{\dashv \vdash}
(\overrightarrow{a}) $ \ служат общими обозначениями
соответственно классов
\[
    \Sigma _{n}( \overrightarrow{a}), ~
    \Pi _{n}( \overrightarrow{a}), ~~~
    \Sigma _{n}^{\dashv \vdash}( \overrightarrow{a}), ~
    \Pi _{n}^{\dashv \vdash }(\overrightarrow{a}).
\]
Натуральный индекс  \ $n$ \ называется их уровнем и уровнем их
формул.
\\
Рассматривая фиксированное расширение \ $L_{k}[ l] $, \ множества,
которые определены в нём формулами из \ $Q_{n}^{\dashv \vdash }
(\overrightarrow{a}) $, будут называться \ $Q_{n}^{\dashv \vdash
}( \overrightarrow{a}) $-множествами.
\\
В результате возникают иерархия формул
\[
    \{ \Sigma _{n}^{\dashv \vdash } ( \overrightarrow{a} ) ;\
    \Pi_{n}^{\dashv \vdash } ( \overrightarrow{a} ) \}_{n \in \omega_0}
\]
и соответствующая иерархия множеств.\\
Во всех этих обозначениях последовательность \
$\overrightarrow{a}$ \ будет опускаться, если она произвольна или
оговорена в контексте.
\end{definition}

\vspace{12pt}

Следуя традиции будем полагать, чтов всякая рассматриваемая
формула предполагается преобразованной в эквивалентную пренексную
форму
 \textit{минимального уровня} \ --- \
эквивалентную в $ZFC^{-}+V=L[\underline{l}]$ или коструктивно или
генерически эквивалентную, в зависимости от структур её
интерпретации \ --- \ разумеется, если контекст не имеет ввиду
другую ситуацию.
\\
Будем также полагать, что если фомулы рассматриваются как
интерпретируемые в структуре \  $ ( L_{k},\; \in ,\; = ) $, \ то
будут использоваться все эти обозначения и термины, но без
верхнего индекса \ $\dashv \vdash $.
\\

В дальнейшем рассматриваются классы \ $Q_{n}$, $Q_{n}^{\dashv
\vdash }$ \ и формулы и множества этих классов некоторого
фиксированного уровня \ $n>3$ \ (если контекст не влечёт
обратное). Это соглашение принимается для того, чтобы использовать
субнедостижимые средства достаточно большого уровня (см. для
примера лемму 3.5 ниже), а также использовать некоторые
вспомогательные формулы, термы, отношения и множества, определимые
в \ $L_{k}$ \ или в \ $L_{k}[l] $, \ непосредственно как
дополнительные константы в обозначениях формул без увеличения их
уровня. Очевидно, таким образом могут быть использованы множества
\ $P$, \ $B$, отношения и операции на них, упомянутые выше, а
также следующие формулы и отношения:
\\
\quad \\
1)\quad $On ( x ) $ \  - ограниченная формула, означающая, что \
$x$ \ есть ординал:
\[
    \quad \forall y\in x ( \forall z\in x ( y\in
    z\vee z\in y ) \wedge \forall z_{1}\in y ( z_{1}\in x )  ) .
\]
Переменные и константы, ограниченные этой формулой, будут также
обозначаться малыми буквами конца греческого алфавита: \ $\alpha
$, $\beta $, $\gamma $, ... , опуская саму эту формулу. С помощью
этой формулы легко можно определить все натуральные числа,
ординалы \ $\omega_{0}, \ \ \omega_{0}+1 .... $ \ и так далее
соответствующими ограниченными формулами, так что можно
использовать эти ординалы как дополнительные константы нашего
языка.\\
2)\quad $F ( x,y ) $ \  --  \  $\Delta _{1}$-формула, сотавляющая
хорошо известное рекурсивное определение в \ $L_{k}$ \ гёделевой
конструктивной функции \ $F$~~\cite{Godel}, то есть для всяких \
$\alpha \in k, \: a\in L_{k} \: :$
\[
    a=F\mid  ( \alpha +1 ) \longleftrightarrow L_{k}\vDash F
    (\alpha, a).
\]
3)\quad $\lessdot $, \ $\underline{\lessdot }$ \quad -- отношения
вполе упорядочения на \quad $L_{k}$:
$$ \quad a\lessdot b\longleftrightarrow Od ( a ) <Od ( b ) \quad ;\quad a\mathbf{\
}\underline{\lessdot }b\longleftrightarrow a\lessdot b\vee a=b .
$$
4)\quad $\vartriangleleft $, \ $\trianglelefteq $ \  --
соответствующие отношения на \ $L_{k}\times k$:
\[
    a \vartriangleleft \beta \longleftrightarrow Od ( a ) < \beta
    \wedge On ( \beta ) \quad ; \quad a \trianglelefteq \beta
    \longleftrightarrow a\vartriangleleft \beta \vee Od ( a ) = \beta.
\]
5) \quad Нетрудно таким же образом использовать аналогичную \ $
\Delta _{1} $ -формулу \ $F ( x,y,\underline{l} )$, \ составляющую
рекурсивное определение в \ $L_{k}\left[ l\right] $ \
конструктивной функции Гёделя \ $F^{l}$ \ относительно \ $l$ \ и
получить функцию ординального номера \ $Od^{l}$, \ а также
отношения \ $\lessdot^{l}$, $\vartriangleleft ^{l}$:

\[
    ~Od^{l} ( a ) =\min \{ \alpha :F^{l} ( \alpha ) =a\} \ \ ;
\]
\[
    a \lessdot ^{l}b\longleftrightarrow Od^{l} ( a ) < Od^{l} ( b )~;
    \quad \quad\quad a\underline{\lessdot }^{l} b\longleftrightarrow
    a \lessdot ^{l}b\vee a=b ~;
\]
\[
    a\vartriangleleft ^{l}\beta\longleftrightarrow Od^{l} ( a )
    <\beta\wedge On ( \beta ) ~; \quad a\trianglelefteq^{l}\beta
    \longleftrightarrow a\vartriangleleft ^{l}\beta\vee
    Od^{l} ( a )=\beta ~.
\]

\vspace{12pt}

Нетрудно определить все эти упомянутые выше функции и отношения \
$\Delta _{1}$-формулами ссответственно в \ $L_{k}$, $L_{k}\left[
l\right] $ \ абсолютно по отношению к этим структурам, используя
формулы \ $F( x,y )$, $F(x,y,\underline{l} )$. Мы будем обозначать
эти формулы теми же символами, что и определяемые ими функции и
отношения, но заменяя функцию \  $l$ \ её именем \
$\underline{l}$. \ Ранее отношение \ $\lessdot $ \ было
использовано Аддисоном~\cite{Addison} над континуумом и
Когаловским~\cite{Kogalovski} над произвольной бесконечной
структурой произвольного уровня.
\\
С помощью обычной техники сколемовских функций легко
устанавливается

\begin{lemma} \label{2.2.} \hfill {} \\
\hspace*{1em} Пусть формула \  $\varphi $ \ содержится в классе \
$ Q_{n}$, \  $n\geq 1$, \ тогда формулы
\[
    \forall x\in y~\varphi, \quad \exists x\in y~\varphi
\]
снова содержатся в том же классе.
\\
Аналогично для класса  \ $Q_{n}^{\dashv \vdash }$, \ $n\geq 1$, \
заменяя ограничивающую формулу \  $ x\in y$ \ на формулы \
$x\lessdot ^{\underline{l}}y$, \ $x\vartriangleleft
^{\underline{l}}y$.   \hspace*{\fill} $\dashv$
\end{lemma}

Теперь обратимся к понятию спектра. Для большей прозрачности мы
введём его только для утверждений вида  \ $\varphi (
\overrightarrow{a},\underline{l} ) $ \  с кортежем индивидных
констант \ $\overrightarrow{a}= ( a_{1},...a_{m} )$, \ состоящим
из ординалов (если контекст не указывает на другой случай). Можно
обойтись без этого соглашения, заменяя вхождения каждого \ $a_{i}$
\ на вхождение терма \ $F^{\underline{l}} ( \alpha _{i} ) $ \ для
соответствующей ординальной константы \ $\alpha _{i}$.
\\
Будем также полагать, что каждый кортеж \ $\overrightarrow{a}= (
\alpha _{1},...,\alpha _{m} ) $ \ ординалов \ $ < k$ \
идентифицируется с ординалом, который является его образом при
каноническом порядковом изоморфизме \  $ ^{m}k$ \ на \ $k$. \
Изоморфизм \ $^{2}k$ \ на \  $k$ \  этого вида будет называться
\textit{функцией пары}.
\\
Следующее понятие играет в дальнейшем важную роль:

\begin{definition} \label{2.3.} \hfill {} \\
\hspace*{1em} Для каждой формулы \ $\varphi $ \  и ординала \
$\alpha _{1}\leq k$ \ через \ $\varphi ^{\vartriangleleft \alpha
_{1}}$  \ обозначается формула, полученная из \ $\varphi$ \
посредством \ $\vartriangleleft ^{\underline{l}}$~-ограничения
всех её кванторов ординалом \ $\alpha _{1}$, \ то есть заменой
всех вхождений таких кванторов \ $\exists x$, $\forall x$ \
соответствующими вхождениями
\[
    \exists x~ ( x\vartriangleleft ^{\underline{l}} \alpha _{1}
    \wedge ... ) , ~~~ \forall x~ ( x\vartriangleleft ^{\underline{l}}
    \alpha _{1}\longrightarrow ... ).
\]
Кроме того, если \  $\alpha _{1}<k$, \  то будем говорить, что \
$\varphi $ \ ограничена ординалом \ $\alpha _{1}$, \ или
релятивизирована к ординалу \ $\alpha _{1}$; \ если при этом
утверждение \ $\varphi ^{\vartriangleleft \alpha _{1}}$ \
выполняется, то будем говорить, что \ $\varphi $ \ выполняется
ниже \ $\alpha _{1}$, \ или что \ $\varphi $ \ сохраняется при
ограничении или релятивизации к \ $\alpha _{1}$.

Та же терминология будет применяться ко всем рассуждениям и
конструкциям с индивидными переменными и константами \
$\vartriangleleft ^{\underline{l}}$- ограниченными ординалом \
$\alpha _{1}$.
\\
Мы будем рассматривать  \ $\vartriangleleft$-ограничения вместо \
$\vartriangleleft ^{\underline{l}} $~-огра\-ни\-че\-ний во всех
подобных обозначениях и понятиях, если формулы, конструкции и
рассуждения интпретируюся в  \ $L_{k}$.
\\
Во все родобных случаях  \ $\alpha_{1}$ \ называется
соответственно \ $\vartriangleleft
^{\underline{l}}$~-ограничивающим или \ $\vartriangleleft
$~-ограничивающим ординалом.
\\
Если \  $\alpha _{1}=k$, \  то верхние индексы \ $\vartriangleleft
\alpha _{1}$ \  опускаются и такие формулы, рассуждения и
конструкции называются неограниченными или нерелятивизированными.
\hspace*{\fill} $\dashv$
\end{definition}

\vspace{-6pt}

\begin{definition}
\label{2.4.}\ \\
{\em 1)}\quad Пусть \  $\varphi (\overrightarrow{a},\underline{l}
) $ \  это утверждение \  $\exists x~\varphi _{1}
(x,\overrightarrow{a}, \underline{l} ) $ и \ $\alpha _{1}\leq k$.
Для каждого \ $\alpha <\alpha _{1}$ \  вводятся следующие булевы
значения:

\vspace{12pt}
\[
    A_{\varphi }^{\vartriangleleft \alpha _{1}} ( \alpha ,
    \overrightarrow{a}) = \left\| \exists x\trianglelefteq ^{\underline{l}}
    \alpha ~\varphi_{1}^{\vartriangleleft \alpha _{1}}
    ( x,\overrightarrow{a},\underline{l} ) \right\| ;
\]

\[
    \Delta _{\varphi }^{\vartriangleleft \alpha _{1}} ( \alpha ,
    \overrightarrow{a} ) =A_{\varphi }^{\vartriangleleft \alpha _{1}}
    ( \alpha ,\overrightarrow{a} ) - \sum_{\alpha ^{\prime }<\alpha}
    A_{\varphi }^{\vartriangleleft \alpha _{1}} ( \alpha ^{\prime },
    \overrightarrow{a} ) .
\]
\vspace{12pt}

\noindent {\em 2)}\quad Мы называем следующую функцию  \
$\mathbf{S}_{\varphi }^{\vartriangleleft \alpha _{1}} (
\overrightarrow{a} ) $ спектром утверждения \  $\varphi $ \ на
кортеже \ $\overrightarrow{a}$ \  ниже \  $\alpha _{1}$:
\vspace{12pt}
\[
    \mathbf{S}_{\varphi }^{\vartriangleleft \alpha _{1}}
    ( \overrightarrow{a}) =\{  ( \alpha ,
    \Delta_{\varphi }^{\vartriangleleft \alpha_{1}}
    ( \alpha ,\overrightarrow{a} )  ) :\alpha <\alpha_{1}
    \wedge \Delta _{\varphi }^{\vartriangleleft \alpha _{1}}
    ( \alpha , \overrightarrow{a} ) >0\}.
\]
\vspace{0pt}

\noindent Проекции
\[
    dom \left( \mathbf{S}_{\varphi }^{\vartriangleleft \alpha _{1}}
    ( \overrightarrow{a} ) \right), \quad  rng \left (
    \mathbf{S}_{\varphi }^{\vartriangleleft \alpha _{1}}
    ( \overrightarrow{a} ) \right )
\]
\vspace{-3pt}

\noindent называются соответственно ординальным и булевым
спектрами \ $\varphi $ \ на кортеже \ $\overrightarrow{a}$ \ ниже
\ $\alpha _{1}$.

\noindent {\em 3)}\quad  Если \ $ ( \alpha ,\Delta  ) \in
\mathbf{S}_{\varphi }^{\vartriangleleft \alpha _{1}} (
\overrightarrow{a} ) $, \ то \  $\alpha $ \  называется ординалом
скачка этих формулы и спектров, а \ $\Delta $ \ называется их
булевым значением на кортеже \ $\overrightarrow{a}$ \ ниже \
$\alpha _{1}$.
\\
\noindent {\em 4)}\quad Ординал \ $\alpha _{1}$ \ называется
носителем этих спектров.  \hspace*{\fill} $\dashv$
\end{definition}
Если кортеж \  $\overrightarrow{a}$ \  пустой, то мы опускаем его
в обозначениях и опускаем все упоминания о нём.
\\

Для анализа различных утверждений естественно использовать их
спектры. Можно осуществить более тонкий анализ, используя их
двухмерные спектры, трёхмерные спектры и так далее. \label{c2}
\endnote{
\ стр. \pageref{c2}. \ Двухмерные спектры строятся следующим
образом. Рассмотрим утверждение:
\[
    \varphi  ( \overrightarrow{a},\underline{l} ) =\exists x\ \varphi_{1}
    ( x,\overrightarrow{a},\underline{l} ) =\exists x\forall y\
    \varphi_{2} ( x,y,\overrightarrow{a},\underline{l} ) .
\]
К каждому ординалу скачка \ $\alpha $ \ формулы \ $\varphi $\ для
кортежа \ $\overrightarrow{a}$\  ниже \ $\alpha _{1}$\ следует
присоединить спектр утверждения \ $\varphi _{1}^{\prime }=\exists
y\neg \varphi _{2} ( F^{\underline{l}} ( \alpha  ) ,y,
\overrightarrow{a},\underline{l} ) $ \  ниже  \ $\alpha _{1}$,\ но
уже на кортеже \ $ ( \alpha ,\overrightarrow{a} ) $.\ В результате
получается функция, которая является двухмерным спектром:
\\
\vspace{6pt}
\[
    \mathbf{S}_{\varphi ,2}^{\vartriangleleft \alpha _{1}}
    ( \overrightarrow{a} ) = \{  ( \alpha ,\beta ,
    \Delta _{\varphi}^{\vartriangleleft \alpha _{1}}
    ( \alpha ,\overrightarrow{a} ) ,
    \Delta _{\varphi _{1}^{\prime }}^{\vartriangleleft \alpha _{1}}
    (\alpha ,\beta ,\overrightarrow{a} )  ) :\alpha ,
    \beta<\alpha_{1}\wedge \\
\]
\[
    \qquad\qquad \qquad\qquad\wedge
    \Delta _{\varphi}^{\vartriangleleft \alpha _{1}}
    ( \alpha ,\overrightarrow{a} )>0\wedge
    \Delta _{\varphi _{1}^{\prime }}^{\vartriangleleft \alpha _{1}}
    ( \alpha ,\beta ,\overrightarrow{a} ) >0 \}.
\]
\vspace{6pt}

\noindent Его первая и вторая проекция составляет двухмерный
ординальный спектр, в то время как третья и четвёртая проекция
составляет двухмерный булев спектр утверждений \ $\varphi $ \  на
кортеже \ $ \overrightarrow{a}$ \  ниже \ $\alpha _{1}$. \
Трёхмерный спектр получается аналогичным образом последующим
присоединением спетров \ $\varphi _{2}$ \  и так далее. Однако
применения многомерных спектров лежат за пределами этой работы.
\\
\quad \\
} 

Все введённые спектры обладают следующими простыми свойствами:
\begin{lemma} \label{2.5.} \quad \\
\hspace*{1em} Пусть \ $\varphi $ \ это утверждение
\[
    \exists x\ \varphi _{1} ( x,\overrightarrow{a},l ), \quad
    \quad \varphi_{1}\in \Pi _{n-1}^{\dashv \vdash }, \quad
    \alpha _{1}\leq k,
\]
тогда:
\\
{\em 1)}\quad $\sup dom \left (
\mathbf{S}_{\varphi}^{\vartriangleleft \alpha_{1}} (
\overrightarrow{a} ) \right ) <k $;
\hfill {} \\
\hfill {} \\
{\em 2)} Спектры \quad $\mathbf{S}_{\varphi }^{\vartriangleleft
\alpha _{1}} ( \overrightarrow{a} ) $, \ $dom \left (
\mathbf{S}_{\varphi }^{\vartriangleleft \alpha _{1}} (
\overrightarrow{a} ) \right) $ \     \ $\Delta _{n}$-определимы, в
то время как спектр
\\
$rng \left (\mathbf{S} _{\varphi }^{\vartriangleleft \alpha _{1}}
( \overrightarrow{a} ) \right )$ \   \ $\Sigma _{n}$-определим в \
$L_{k}$ \ для \ $\alpha _{1}=k$.
\\
Для \ $\alpha _{1}<k $ \ все эти спектры \ $\Delta_{1}$-определимы
в \ $L_{k}$;
\hfill {} \\
\hfill {} \\
{\em 3)}\quad $\alpha \in dom \left (\mathbf{S}_{\varphi
}^{\vartriangleleft \alpha _{1}} ( \overrightarrow{a} ) \right)$ \
в точности тогда, когда существует \ $\mathfrak{M}$-генерическая
функция
\[
    l\stackrel{\ast }{\in }\Delta _{\varphi }^{
    \vartriangleleft \alpha _{1}} ( \alpha , \overrightarrow{a}
    ).\
\]

\end{lemma}
Здесь важное утверждение 1) следует непосредственно из леммы
\ref{1.1.}; это утверждение будет часто использоваться в
дальнейшем.  \hspace*{\fill} $\dashv$
\\

Среди всех спектров особую роль играют так называемые
\textit{универсальные спектры}.
\\
Хорошо известно, что класс \ $\Sigma_{n} ( \overrightarrow{a})$ \
уровня \  $n>0$ \  содержит формулу, универсальную для этого
класса (Тарский~\cite{Tarski}, см. Аддисон \cite{Addison}); мы
будем обозначать её через \ $U_{n}^{\Sigma } ( \mathfrak{n},
\overrightarrow{a},\underline{l} ) $. \ Следовательно, она
универсальна также и для класса \
$\Sigma_{n}^{\dashv\vdash}(\overrightarrow{a}) $. \ Её
универсальность означает, что для любой \ $\Sigma
_{n}^{\dashv\vdash} ( \overrightarrow{a} ) $-формулы \ $ \varphi (
\overrightarrow{a},\underline{l} ) $ \ существует натуральное \
$\mathfrak{n}$ \ такое, что
\[
    \varphi  ( \overrightarrow{a},\underline{l} ) \dashv \vdash
    U_{n}^{\Sigma } ( \mathfrak{n},\overrightarrow{a},\underline{l} ) ;
\]
это \ $\mathfrak{n}$ \ называется гёделевым номером формулы \
$\varphi $. \ Двойственная формула, универсальная для
$\Pi_{n}^{\dashv\vdash} ( \overrightarrow{a} ) $, \ обозначается
через $ U_{n}^{\Pi } ( \mathfrak{n},\overrightarrow{a},\underline{
l} ) $. \ Для некоторого удобства мы будем использовать \
$U_{n}^{\Sigma }$ \ в форме \ $\exists x~U_{n-1}^{\Pi } (
\mathfrak{n},x,\overrightarrow{a}, \underline{l} ) $. \ В этих
обозначениях верхние индексы \ $^{\Sigma} $,~$^{\Pi} $ \ будут
опускаться, если они могут быть восстановлены из контекста либо
произвольны.
\\
Напомним, что универсальные формулы возникают на следующем пути.
Рассмотрим формулу \ $ \varphi\left( \overrightarrow{{a}},
\underline{l}\right)$ \ и её эквивалентную пренексную форму
\textit{минимального} уровня $e$:
\[
    Q_{1}x_{1}Q_{2}x_{2}...,Q_{i}x_{i}\varphi_{1}\left(
    x_{1},x_{2}...,x_{i},\overrightarrow{{a}},\underline{l}\right),
\]
где \ $Q_{1}x_{1}Q_{2}x_{2}...,Q_{i}x_{i}$ \ это её кванторный
префикс и \ $ \varphi_{1}$ \ содержит только ограниченные
кванторы; предположим, что \ $Q_{i}=\exists, \quad e>0$. \
Передвигая эти последние ограниченные кванторы влево можно
получить эквивалентную формулу вида
\begin{eqnarray*}
    Q_{1}x_{1}Q_{2}x_{2}...,Q_{i}x_{i}
    \quad \qquad \qquad \qquad \qquad \qquad \qquad \qquad \qquad
    \qquad \qquad
\\
    \quad \forall y_{1}\in z_{1} \exists v_{1}\in w_{1}\forall y_{2}
    \in z_{2}\exists v_{2} \in w_{2}...,\forall y_{j}\in z_{j}
    \exists v_{j}\in w_{j}\qquad \qquad
\\
    \varphi_{2}\left(x_{1},x_{2}...,x_{i},y_{1},z_{1},
    v_{1},w_{1},y_{2},z_{2},v_{2},w_{2}...,y_{j},z_{j},v_{j},
    w_{j},\overrightarrow{{a}},\underline{l}\right),
\end{eqnarray*}
\vspace{0pt}

\noindent где каждая из переменных \ $z_{1},w_{1}...,z_{j},w_{j}$
\ ограничена другими или термами из кортежа \
$x_{1}...,x_{i},\overrightarrow{{a}}$ \ и \ $\varphi_{2}$ \ не
содержит кванторов. Используя технику сколемовских функций можно
легко преобразовать её в эквивалентную формулу$ \
U_{\varphi}\left( \overrightarrow{{a}},\underline{l}\right)$ \
того же уровня \ $e$:
\[
    Q_{1}x_{1}Q_{2}x_{2}...,Q_{i}x_{i}\exists f_{1} \exists
    f_{2}...,\exists f_{j}\quad
    \varphi_{3}\left(x_{1},x_{2}...,x_{i},f_{1},f_{2}...,f_{j},
    \overrightarrow{{a}},\underline{l}\right),
\]
где \ $f_{1}...,f_{j}$ \ - это переменные сколемовские функции и \
$\varphi_{3}$ \ содержит только ограниченные кванторы стандартного
вида и расположения.
\\
Случай \ \mbox{$Q_{i}=\forall$} \ следует рассмотреть двойственным
образом.

Мы будем называть эту формулу \ $U_{\varphi}\left(
\overrightarrow{{a}},\underline{l}\right)$
\textit{предуниверсальной формой} рассматриваемой формулы \
$\varphi \left(\overrightarrow{{a}}, \underline{l}\right)$. \
После этого легко получить эквивалентную формулу с каждым блоком
одноименных кванторов, заменённым на один такой квантор,
посредством функции пары. Осталось применить к любой такой формуле
рекурсивную нумерацию её подформул, следующих за кванторным
префиксом.

Положим, что в дальнейшем все вводимые формулы мы будем
рассматривать в их \textit{предуниверсальной форме} и будем
использовать для них те же обозначения (если только контекст не
имеет ввиду другую ситуацию). Это соглашение сохраняет
эквивалентность формулы \ $\varphi (a,\underline{l}) \in
\Sigma_{e}^{\dashv\vdash} $ \ с гёделевым номером  \
$\mathfrak{n}$ \ универсальной формуле при \ $\vartriangleleft^{
\underline{l}}$~-ограничении:
\[
    \varphi^{\triangleleft \alpha_{1}}(\overrightarrow{a},
    \underline{l}) \longleftrightarrow U_{e}^{\Sigma \triangleleft
    \alpha_{1}}( \mathfrak{n},\overrightarrow{a},\underline{l})
\]
для ординалов  \  $\alpha_{1} $ \  многих видов (например, для
кардиналов). Аналогично для \ $\Pi_{e}^{\dashv\vdash}$ \ и для
генерической эквивалентности.

\begin{definition}
\label{2.6.}
\quad \\
{\em 1)}\quad Мы называем формулой спектрально универсальной для
класса \ $\Sigma_{n}^{\dashv \vdash }$ \  уровня \ $n$ \ формулу \
$u_{n}^{\Sigma } (\overrightarrow{a}, \underline{l} ) $, \
полученную из универсальной формулы \ $U_{n}^{\Sigma } (
\mathfrak{n},\overrightarrow{a},\underline{l} ) $ \  заменой всех
вхождений переменной \ $\mathfrak{n}$ \ вхождениями терма \
$\underline{l} ( \omega _{0} ) $.
\\
Спектрально универсальная для класса $\Pi_{n}^{\dashv \vdash }$
формула \ $u_{n}^{\Pi } ( \overrightarrow{a},\underline{l} )$
вводится двойственным образом. Таким образом, будет использоваться
формула
\[
     u_{n}^{\Sigma } ( \overrightarrow{a}, \underline{l} ) =
     \exists x~u_{n-1}^{\Pi } ( x,\overrightarrow{a},\underline{l}
     ),
\]
где \ $u_{n-1}^{\Pi} ( x, \overrightarrow{a}, \underline{l} )$ \
это спектрально универсальная для класса \ $\Pi_{n-1}^{\dashv
\vdash }$ \ формула.
\quad \\
{\em 2)}\quad Булевы значения
\[
    A_{\varphi}^{\vartriangleleft \alpha _{1}} ( \alpha ,
    \overrightarrow{a} ), \Delta _{\varphi }^{\vartriangleleft
    \alpha _{1}} ( \alpha,\overrightarrow{a} )
    \mbox{\it \ и спектр \ }\mathbf{S}_{\varphi}^{\vartriangleleft
    \alpha _{1}} ( \overrightarrow{a} )
\]
формулы \ $\varphi =u_{n}^{\Sigma } (\overrightarrow{a},
\underline{l} )$ \ и его проекции называются универсальными
булевыми значениями и спектрами уровня \ $n$ \  на кортеже \
$\overrightarrow{a}$ \  ниже \ $\alpha _{1}$ \ и в их обозначениях
индекс \ $u_{n}^{\Sigma }$  \ заменяется на \ $n$, \ то есть они
обозначается соответственно через
\[
    A_{n}^{\vartriangleleft \alpha _{1}} ( \alpha ,
    \overrightarrow{a} ), \quad \Delta _{n}^{\vartriangleleft
    \alpha_{1}} ( \alpha ,\overrightarrow{a} ) ,\quad
    \mathbf{S}_{n}^{\vartriangleleft \alpha _{1}}
    ( \overrightarrow{a} ).
\]
\hfill {} \\
{\em 3)}\quad Если \ $( \alpha ,\Delta  ) \in
\mathbf{S}_{n}^{\vartriangleleft \alpha _{1}} ( \overrightarrow{a}
) $, \ то \ $\alpha $ \ называется ординалом скачка этих формулы и
спектров, а \ $\Delta $ \ называется их булевым значением на
кортеже \ $\overrightarrow{a}$ \ ниже \ $\alpha _{1}$.
\\
{\em 4)}\quad Ординал \  $\alpha _{1}$ \ называется носителем этих
спектров. \hspace*{\fill} $\dashv$
\end{definition}

Везде в дальнейшем  \  $\vartriangleleft ^{\underline{l}}$~- \ или
\ $\vartriangleleft $-ограничивающие ординалы \ $\alpha_1 $ \ это
предельные кардиналы  \ $<k $ \ (в \ $L_{k} $) или \ $\alpha_1 =k
$ \ (если контекст не указывает на другое).
\\
\hfill {} \\
Здесь термин ``универсальный'' оправдывает следующая

\begin{lemma} \label{2.7.} \hfill {} \\
\hspace*{1em} Для любого утверждения \ $\varphi =\exists x\
\varphi _{1} ( x, \overrightarrow{a},\underline{l} ) ,\quad
\varphi _{1}\in \Pi _{n-1}^{\dashv \vdash }$:
\[
    dom \left (
    \mathbf{S}_{\varphi }^{\vartriangleleft \alpha _{1}} (
    \overrightarrow{a} )\right ) \subseteq dom \left (
    \mathbf{S}_{n}^{\vartriangleleft \alpha _{1}}
    (\overrightarrow{a} ) \right) .
\]
\end{lemma}
\vspace{-12pt}

\noindent \textit{Доказательство.} \ Для функции  \  $l\in {}^{k}k
$ \ пусть \ $l_0 $  \  обозначает любую функцию $\in {}^{k}k $, \
принимающую значения:
\\
\begin{equation} \label{e2.1}
l_0 \left( \alpha\right) = \left\{
\begin{array}{ll}
l(\alpha) & ,  \ \alpha< \omega_{0} \vee \alpha > \omega_{0}+1; \\
(l(\omega_{0}), l(\omega_{0}+1)) & , \ \alpha = \omega_{0}+1;
\end{array} \right.
\end{equation}
\hfill {} \\
значение  \ $l_0(\omega_{0}) $ \ здесь произвольно. Очевидно,\
$l_0$ \ это \ $\mathfrak{M}$-генерическая функция, если \ $l $  \
обладает тем же свойством. Пусть \ $\varphi_{0} \left
(\overrightarrow{a}, l \right )$  \  обозначает формулу,
полученную последовательной заменой каждой субформулы \
$\underline{l} \left (t_{1} \right ) = t_{2} $ \ формулы \
$\varphi $  \ на ограниченную формулу:
\\
\begin{eqnarray} \label{e2.2}
& \nonumber \exists y_{1},y_{2} < \omega_{0} ( \underline{l}
\left( \omega_{0}+1 \right) = (y_{1},y_{2})\wedge
\\
& \wedge \left ( \left ( t_{1} < \omega_{0} \vee t_{1}
> \omega_{0} + 1 \right) \longrightarrow \underline{l} (t_{1}) =
t_{2}  \right) \wedge
\\
& \nonumber \wedge \left ( t_{1} = \omega_{0} \longrightarrow
y_{1} = t_{2} \right )\wedge \left ( t_{1} = \omega_{0} +1
\longrightarrow y_{2} = t_{2} \right ) ).
\end{eqnarray}
\hfill {} \\
Подформулы вида  \  $ \underline{l} (t_{1}) \in t_{2}, \quad t_{2}
\in \underline{l} (t_{1}) $  \  заменяются аналогично. Пусть \ $
\mathfrak{n}_{0} $  \  это гёделев ноомер \ $\varphi_{0}$. \ Для
произвольной  \ $\mathfrak{M}$-генерической функции \ $l $ \
положим \ $l_0(\omega_{0}) = \mathfrak{n}_{0} $. Очевидно, для
каждого  \ $\alpha < \alpha_{1} $

\[
    l \stackrel{\ast }{\in }~ A_{\varphi}^{\triangleleft \alpha_{1}}
    (\alpha, \overrightarrow{a}) \longleftrightarrow l_0
    \stackrel{\ast}{\in }~ A_{\varphi_{0}}^{\triangleleft \alpha_{1}}
    (\alpha, \overrightarrow{a}) \longleftrightarrow l_0
    \stackrel{\ast }{\in }~ A_{n}^{\triangleleft \alpha_{1}}
    (\alpha, \overrightarrow{a})\quad.
\]
\vspace{0pt}

Осталось применить теперь лемму 2.5 \ 3). \hfill $\dashv$ \\

Отсюда можно видеть, что универсальные спектры аккумулируются
когда их ординальные константы возрастают; для \
$\overrightarrow{a}= ( \alpha_{1}...,\alpha _{m} )$ \ пусть \
$\max \overrightarrow{a} =\max \{ \alpha _{1}...,\alpha _{m}\} $:
\\

\begin{lemma}\label{2.8.} \hfill {} \\
\hspace*{1em} Для всяких \ $\overrightarrow{a}_{1},
\overrightarrow{a}_{2}<\alpha _{1}$:
\hfill {} \\
\[
    \max \overrightarrow{a}_{1}<\max \overrightarrow{a}_{2}
    \longrightarrow dom \left ( \mathbf{S}_{n}^{\vartriangleleft
    \alpha _{1}} (\overrightarrow{a}_{1} )\right )
    \subseteq dom \left ( \mathbf{S}_{n}^{\vartriangleleft
    \alpha _{1}} (\overrightarrow{a}_{2} )\right) .
\]
\end{lemma}
Доказательство этого может быть получено так называемым
расщепляющим методом ( см. доказательство леммы 4.6 ниже для
примера), но эта лемма, хотя и проясняющая спектральные свойства,
не применяется далее и поэтому её доказательство здесь опускается.
\hspace*{\fill} $\dashv$

\newpage
\quad 

\newpage

\section{Субнедостижимые кардиналы}
\setcounter{equation}{0}

Здесь разрабатывается теория субнедостижимости в её базовых
аспектах.

Дальнейшие рассуждения проводятся в \ $L_{k}$ \ (или в  \
$\mathfrak{M} $, \  если контекст не имеет ввиду некоторую другую
ситуацию).
\\
Введём центральное понятие субнедостижимости -- недостижимость
средствами нашего языка. ``Значение'' утверждений заключается в их
спектрах и поэтому весьма естественно определить такую
недостижимость посредством спектров заданного уровня:

\begin{definition} \label{3.1.} \hfill {} \\
\hspace*{1em} Пусть \  $\alpha _{1} \leq k $.  \\
Мы называем ординал  \ $\alpha <\alpha _{1}$ \ субнедостижимым
уровня \ $n$ \ ниже \ $\alpha _{1}$, \ если выполняется следующая
формула, обозначаемая через \ $SIN_{n}^{<\alpha _{1}} ( \alpha )$:
\hfill {} \\
\[
    \forall \overrightarrow{a} <\alpha \quad dom
    \left ( \mathbf{S}_{n}^{\vartriangleleft \alpha_{1}}
    ( \overrightarrow{a} ) \right ) \subseteq \alpha .
\]
Множество \
\[
    \{ \alpha <\alpha _{1}:SIN_{n}^{<\alpha _{1}} ( \alpha ) \}
\]
всех таких ординалов также обозначается через \ $SIN_{n}^{<\alpha
_{1}}$ \ и они называются \ $SIN_{n}^{<\alpha _{1}}$-ординалами.
\\
Как обычно, для \ $\alpha _{1}<k$ \  мы говорим при этом, что
субнедостижимость ординала \ $\alpha $ \ ограничена ординалом \
$\alpha _{1}$ \ или релятивизирована к \ $\alpha _{1}$. \ Для \
$\alpha_1 = k$ \ верхние индексы \ $< \alpha_1$, $\vartriangleleft
\alpha_1$ \ опускаются. \hfill $\dashv$
\end{definition}

Очевидно, сам кардинал  \ $k$ \ является \textit{субнедостижимым}
любого уровня, если мы определяем это понятие для \ $\alpha =
\alpha_1 = k$.
\\
Таким образом, сравнение понятий недостижимости и
субнедостижимости естественно возникает на следующем пути:
\\
Кардинал \ $k$ \ (слабо) недостижим, так как он несчётен и не
может быть достигнут средствами меньших мощностей в том смысле,
что он обладает двумя свойствами: 1) он регулярен и 2) он замкнут
относительно операции перехода к следующей мощности: \ $\forall
\alpha < k ~~~ \alpha^+ < k$.
\\
Обращаясь к \textit{субнедостижимости} ординала \ $\alpha < k$ \
уровня \ $n$ \ (в её нерелятивизированной форме, для некоторой
краткости), можно видеть, что свойство регулярности теперь
опущено, но \ $\alpha$ \ всё-таки не может быть достигнут, но на
сей раз более сильными средствами: второе условие усилено и теперь
\ $\alpha$ \ закнут относительно более сильной операции перехода к
ординалам скачка универсального спектра:
\[
    \forall \overrightarrow{a} <\alpha \quad \forall \alpha^{\prime}
    \in dom \mathbf{S}_{n} ( \overrightarrow{a} )
    \quad \alpha^{\prime} < \alpha,
\]
то есть средствами ординальных спектров \textit{всех} утверждений
уровня \ $n$ \ (см. лемму 2.7 выше).
\\
Это означает замкыкание \ $\alpha$ \ \textit{относительно всех \
$\Pi_{n-1}^{\dashv\vdash}$-функций во всех генерических
расширениях} \ $L_k$, \ а не только относительно операции перехода
к следующей мощности в  \ $L_k$ \ (см. лемму 3.5 ниже).
\\

Очевидно, что действуя в \ $L_{k} $ \ следует рассматривать
формулу \ $SIN_{n}^{< \alpha_{1}}(\alpha) $ \ фактически как две
формулы: одну из них без константы \ $\alpha_{1} $ \ когда \
$\alpha_{1}=k$, \ и другую содержащую \ $\alpha_{1} $ \ когда \
$\alpha_{1} < k$; \ такое же замечание относится ко всем формулам,
конструкциям и рассуждениям, содержащим параметр \
$\alpha_{1} \leq k$.\\
\quad \\

Из определения 3.1 и леммы 2.7 очевидно следует

\begin{lemma}
\label{3.2.} {\em (Об ограничении)} \hfill {} \\
\hspace*{1em} Пусть \  $\alpha <\alpha _{1}\leq k$, \ $\alpha \in
SIN_{n}^{<\alpha _{1}}$  \ и  \ $\Sigma _{n}^{\dashv \vdash }$--\
утверждение \  $ \exists x~\varphi  ( x,\overrightarrow{a},
\underline{l} ) $,\   \ $\varphi \in \Pi _{n-1}^{\dashv \vdash }$,
\ имеет кортеж \  $ \overrightarrow{a}<\alpha $, \ тогда для
каждой \ $\mathfrak{M}$-генеричекой функции \ $l$
\\
\[
    L_{k}\left[ l\right] \vDash \left( \exists x\vartriangleleft ^{l}
    \alpha _{1}\mathit{\ }\varphi ^{\vartriangleleft \alpha _{1}}
    ( x,\overrightarrow{a}, l ) \longrightarrow \exists x
    \vartriangleleft ^{l}\alpha \ \varphi ^{\vartriangleleft
    \alpha _{1}} ( x,\overrightarrow{a},l) \right) \ ,
\]
\hfill {} \\
В таком случае будем говорить, что ниже \  $\alpha _{1}$ \ ординал
\  $\alpha $ \ ограничивает или релятивизирует утверждение \
$\exists x~\varphi $.
\\
Рассматривая то же самое в обращённой форме для \ $\varphi \in
\Sigma _{n-1}^{\dashv \vdash }$:
\\
\[
    L_{k}\left[ l\right] \vDash \left( \forall x\vartriangleleft ^{l}
    \alpha \ \varphi ^{\vartriangleleft \alpha _{1}} ( x,
    \overrightarrow{a},l ) \longrightarrow \forall x
    \vartriangleleft ^{l}\alpha _{1}\varphi ^{\vartriangleleft
    \alpha _{1}} ( x,\overrightarrow{a},l ) \right) \ ,
\]
\hfill {} \\
будем говорить, что ниже \  $\alpha _{1}$ \  ординал \ $\alpha $ \
расширяет или продолжает утверждение \ $\forall x\ \varphi $ \ до
\ $\alpha _{1}$. \hspace*{\fill} $\dashv$
\end{lemma}

Разумеется, лемма~3.2 представляет более сильное утверждение --
критерий \ $SIN_{n}^{<\alpha _{1}}$-субнедостижимости. \label{c3}
\endnote{
\ стр. \pageref{c3}. \ Очевидно, понятие субнедостижимости можно
сформулировать в терминах элементарных подмоделей, но только в
следующей существенно более сложной и искусственой форме:

Ординал $\alpha$ \ является субнедостижимым уровня \ $n$ \ ниже \
$\alpha_1$, \ если  $L_{\alpha}[l | \alpha_1]$ \ это \
$\Sigma_n$-элементарная подмодель \ $L_{\alpha_1}[l | \alpha_1]$ \
\textit{для каждой \ \mbox{$\mathfrak{M}$-генерической} функции} \
$l$.
\\
Именно это делается в лемме 3.2. Здесь несколько слов следует
сказать по этому поводу.
\\
Это понятие не может быть естественно редуцировано к элементарной
эквивалентности конструктивных сегментов \ $L_{\alpha}$, \
$L_{\alpha_1}$ \ только, но требует вовлечения всех их расширений
указанного вида.
\\
Таким образом, используя подобную терминологию, когда параметры \
$\alpha$, $\alpha_1$ \ варьируются одновременно с другими также
субнедостижимыми кардиналами, мы получаем некоторое многоуровневое
и громоздкое описание понятия, которое особенно неестественно,
поскольку на самом деле это понятие работает только внутри  \
$L_{k} $.
\\
Более того, \ и понятие диссеминатора (\S 6), которое является
простым обобщением понятия субнедостижимости, становится
неестественно усложнённым в терминах элементарных подмоделей.
Поэтому требуется некоторое более удобное описание этого понятия,
указывающее непосредственно на сущность этого явления, просто
определимое и поэтому пригодное для исследования поставленной
проблемы -- то есть субнедостижимость посредством формульных
спектров, введённая выше.
\\
\quad \\
} 
\\

Теперь следующие леммы \ref{3.3.} -- \ref{3.8.} можно просто
вывести из определения 3.1 и леммы \ref{2.5.}:

\begin{lemma} \label{3.3.} \hfill {} \\
\hspace*{1em} Формула \  $SIN_{n}^{<\alpha _{1}} ( \alpha ) $ \
принадлежит классу \ $\Pi _{n} $ \ для \ $\alpha _{1}=k$ \ и
классу \ $\Delta _{1} $ \ для \ $\alpha _{1}<k$.
\end{lemma}

\vspace{-6pt}

\begin{lemma} \label{3.4.} \hfill {} \\
\hspace*{1em} Для каждого \ $n>0$:\\
\quad \\
{\em 1)}\quad множество \  $SIN_{n}^{<\alpha _{1}}$ \  замкнуто в
 \ $\alpha _{1}$, \ то есть для каждого \  $\alpha <\alpha
_{1}$
\[
    \sup  ( \alpha \cap SIN_{n}^{<\alpha _{1}} ) \in
    SIN_{n}^{<\alpha_{1}};
\]
{\em 2)}\quad множество \  $SIN_{n}$ \  неограничено в \ $k$,
то есть выполняется \ $\sup SIN_{n}=k$;\\
\quad \\
{\em 3)}\quad \quad \quad \quad \quad \quad \quad \quad \quad
$SIN_{n}^{<\alpha _{1}} ( \alpha  ) \longleftrightarrow
SIN_{n}^{\vartriangleleft \alpha _{1}} ( \alpha )$.
\end{lemma}

\vspace{-6pt}

\begin{lemma} \label{3.5.} \hfill {} \\
\hspace*{1em} Пусть \  $\alpha \in SIN_{n}^{<\alpha _{1}}$ \ и
функция \ $f\subset \alpha _{1}\times \alpha _{1}$ \ определена в
\ $L_{k}[l] $  \  формулой \  $\varphi ^{\vartriangleleft \alpha
_{1}} ( \beta ,\gamma ,l ) $, \ \ $\varphi \in \Pi _{n-1}^{\dashv
\vdash}$,
\\
тогда \  $\alpha $ \  замкнут относительно \  $f$.
\\
В частности, для каждого \  $n \geq 2$
\[
    \mbox{\it \quad если \quad} \alpha \in SIN_{n} ,\mbox{\it \quad то \quad}
    \alpha =\omega _{\alpha } \quad (\mbox{\it в \ } L_{k} ).
\]
\end{lemma}
Следующие леммы 3.6-3.8 представляют важные технические средства
исследования субнедостижимости.

\begin{lemma} \label{3.6.} \hfill {} \\
\hspace*{1em} Для всякого \ $m<n$
\\
1) \hspace*{3.5cm} $SIN_{n}^{<\alpha _{1}} \subset
SIN_{m}^{<\alpha _{1}}$;
\\
\\
2) более того, каждый \ $\alpha \in SIN_{n}^{<\alpha _{1}}$ \ это
предельный ординал в \ $SIN_{m}^{<\alpha _{1}}$~:
\[
    \sup \left ( \alpha  \cap SIN_{m}^{<\alpha_{1}} \right ) = \alpha
\]
\end{lemma}
\textit{Доказательство.} \ Утверждение 1) очевидно, так как каждая
\ $\Sigma_m$-формула есть в то же время \ $\Sigma_n$-формула.
\\
Обращаясь к 2) рассмотрим \ $\Sigma_n$-утверждение
\[
    \exists\gamma(\beta<\gamma\wedge SIN_m(\gamma))
\]
с произвольной константой \ $\beta<\alpha$. \ Это утверждение
выполняется ниже \ $\alpha_1$, потому что благодаря 1) сам ординал
\ $\alpha$ \ может быть использован как \ $\gamma$. \ После этого
\ $SIN_{n}^{<\alpha_{1}}$-ординал \ $\alpha$ \ ограничивает это
утверждение и некоторый \ $SIN_{m}^{<\alpha_{1}}$-ординал \
$\gamma>\beta$ появляется ниже \ $\alpha$. \hfill $\dashv$
\\

\noindent Очевидно, что обратное утвеждение неверно:
субнедостижимые кардиналы иногда теряют это свойство на следующем
уровне -- например, все наследники в данном классе \
$SIN_{m}^{<\alpha _{1}}$.

\begin{lemma} \label{3.7.} \hfill {} \\
\hspace*{1em} Пусть
\[
    \overrightarrow{a}<\alpha _{2}<\alpha_{1}\leq k ,\quad
    \alpha _{2}\in SIN_{n}^{<\alpha _{1}}
\]
тогда для всякого \ $ Q_{n}^{\dashv\vdash}$-утверждения \ $\varphi
( \overrightarrow{a},\underline{l} )$
\\
\[
     \varphi ^{\vartriangleleft \alpha _{1}} ( \overrightarrow{a},
     \underline{l} ) \dashv \vdash \varphi ^{\vartriangleleft
     \alpha_{2}} ( \overrightarrow{a},\underline{l} ).
\]
\end{lemma}
\textit{Доказательство.} \ Рассмотрим утверждение \ $\varphi =
\exists x \ \ \varphi_{1} \left (x,\overrightarrow{a}
,\underline{l} \right ) $,  \ $\varphi_{1}\in \Pi_{n-1} $, \ и
ординал
\[
    \alpha_{0} = min \left \{ \alpha : L_{k} [l] \vDash
    \varphi_{1}^{\triangleleft \alpha_{1}} \left (F^{l}(\alpha),
    \overrightarrow{a}, l \right ) \right \}
\]
для  \  $\mathfrak{M} $-генерической  \  $l $. \ Согласно леммам
2.5 \ 3), 2.7 \ $\alpha_{0}\in dom \left
(\mathbf{S}_{n}^{\triangleleft \alpha_{1}} (\overrightarrow{a})
\right ) $. \ Так как  \ $\alpha_{2}\in SIN_{n}^{< \alpha_{1} }$,\
то это влечёт  \ $\alpha_{0} < \alpha_{2} $ \  и поэтому \
$\vartriangleleft^{l}$-ограничение ординалом \  $\alpha_{1}$ \
утверждения \  $ \varphi_{1}^{\triangleleft \alpha_{1}} \left
(F^{l}(\alpha_{0}), \overrightarrow{a}, l \right ) $ \  может быть
заменено на  \ $\vartriangleleft^{l}$-ограничение ординалом \
$\alpha_{2}$. Следовательно
\[
    L_{k}[l]\vDash \left( \varphi^{\triangleleft \alpha_{1}}
    (\overrightarrow{a},l) \longrightarrow
    \varphi^{\triangleleft \alpha_{2}}(\overrightarrow{a},l) \right) \ .
\]
Осталось обратить это рассуждение. \hfill $\dashv$

\begin{lemma} \label{3.8.} \hfill {} \\
\hspace*{1em} Пусть \ $\alpha _{2}<\alpha _{1}\leq k$, \ тогда:\\
\quad \\
{\em 1)}\quad Если \  $\alpha _{2}\in SIN_{n-1}^{<\alpha _{1}}$,
то множество \  $SIN_{n}^{<\alpha _{1}}\cap \alpha _{2}$ \
составляет начальный сегмент множества \ $SIN_{n}^{<\alpha_{2}}$,
\  то есть:
\\

 \quad \

\qquad (i)\quad $SIN_{n}^{<\alpha _{1}}\cap \alpha
_{2}\subseteq SIN_{n}^{<\alpha _{2}}$;\\

\qquad (ii)\quad $SIN_{n}^{<\alpha_{2}} \cap \sup \left (
SIN_{n}^{<\alpha_{1}} \cap \alpha_{2} \right )
\subseteq SIN_{n}^{<\alpha_{1}} $.\\
\quad \\
\hfill {} \\
{\em 2)}\quad Если \  $\alpha _{2}\in SIN_{n}^{<\alpha _{1}}$, то
для каждого \ $m \leq n$ \
\[
    SIN_{m}^{<\alpha _{1}}\cap \alpha_{2}=SIN_{m}^{<\alpha _{2}} .
\]
\end{lemma}
\textit{Доказательство.} \ Утверждения этой леммы доказываются
аналогичным образом, поэтому продемонстрируем это рассуждение для
 $1.(i).$
\\
Пусть \ $\alpha \in SIN_{n}^{<\alpha _{1}}\cap \alpha_{2}$; \
будет установлено, что для каждого \ $\overrightarrow{a} < \alpha
$
\[
    dom \left ( \mathbf{S}_{n}^{\triangleleft \alpha_{2}}
    (\overrightarrow{a}) \right )  \subseteq \alpha,
\]
для этого рассмотрим любое \ $\beta \in dom \left (
\mathbf{S}_{n}^{\triangleleft \alpha_{2}} (\overrightarrow{a})
\right ) $. \ Теперь имеется  \ $\beta < \alpha_{2} ,
\overrightarrow{a} < \alpha_{2} $ \  и \  $\alpha_{2} \in
SIN_{n-1}^{<\alpha_{1}} $, поэтому для каждого  \  $\beta^{'} \leq
\beta $ \  можно заменить  \ $\vartriangleleft^{\underline{l}}
$~-ограничение ординалом \ $\alpha_{2} $ \ утверждения \
$u_{n-1}^{\Pi \triangleleft \alpha_{2}} \left (F^{\underline{l}}
(\beta^{'}), \overrightarrow{a}, \underline{l} \right ) $ \ на  \
$\vartriangleleft^{\underline{l}} $~-ограничение ординалом \
$\alpha_{1} $. \ Отсюда и из определения 3.1 следует \ $\beta \in
dom \left ( \mathbf{S}_{n}^{\triangleleft \alpha_{1}}
(\overrightarrow{a}) \right ) $ \  и тогда  \  $ \alpha \in
SIN_{n}^{<\alpha_{1}} $ \ влечёт  \  $\beta < \alpha $. \ Это
означает, что  \  $dom \left ( \mathbf{S}_{n}^{\triangleleft
\alpha_{2}} (\overrightarrow{a}) \right ) \subseteq \alpha $ \ и
поэтому  \ $\alpha \in
SIN_{n}^{<\alpha_{2}} $.  \hfill $\dashv$\\
\quad \\

Когда утверждения эквивалентно преобразуются, то их спектры могут
меняться. Можно использовать это явление для анализа
субнедостижимых кардиналов. С этой целью мы введём универсальные
формулы, чьи ординальные спектры состоят только из субнедостижимых
кардиналов меньшего уровня. Для большей ясности всех конструкций
будут рассматриваться формулы без индивидных констант. Начнём со
спектрально универсальной формулы для класса \ $\Sigma
_{n}^{\dashv \vdash }$. \  Верхние индексы \ $^{\Sigma} $, $
^{\Pi}$ \ будут как правило опускаться (там где это не повлечёт
недоразумений). \\
В дальнейшем достаточно использовать ограничивающие ординалы
только из класса \ $SIN_{n-2}$, \ поэтому везде далее в качестве \
$\vartriangleleft^{\underline{l}}$\;- \ или \ $\vartriangleleft
$-ограничивающих ординалов \ $\alpha $ \ рассматриваются только \
$SIN_{n-2}$-ординалы или \ $\alpha =k$ \ (если контекст не
означает иное). \\
Следовательно, все такие ординалы \ $\alpha \leq k$ \ являются
кардиналами вида $\alpha=\omega_{\alpha}$ благодаря лемме 3.5.
\hfill {} \\
\begin{definition}
\label{3.9.} \hfill {} \\
{\em 1)}\quad Мы называем монотонной спектрально универсальной для
класса \ $\Sigma _{n}^{\dashv \vdash }$ \ формулой уровня \ $n$ \
следующую \ $\Sigma _{n}$-формулу
\[
    \widetilde{u}_{n} ( \underline{l} ) =\exists x~\widetilde{u}_{n-1}
    ( x, \underline{l} )
\]
где \  $\widetilde{u}_{n-1} ( \underline{l} ) \in \Pi _{n-1}$ \ и
\[
    \widetilde{u}_{n-1} ( x,\underline{l} ) \dashv \vdash
    \exists x^{\prime}\vartriangleleft
    ^{\underline{l}}x~u_{n-1}^{\Pi}
    ( x^{\prime },\underline{l} ) .
\]
\hfill {} \\
{\em 2)}\quad Мы называем субнедостижимо универсальной для класса
\ $\Sigma _{n}^{\dashv \vdash }$ \ формулой уровня \ $n$ \
следующую \ $\Sigma _{n}$-формулу
\[
    \widetilde{u}_{n}^{\sin } ( \underline{l} ) =
    \exists x~\widetilde{u}_{n-1}^{\sin} ( x,\underline{l} )
\]
где \ $\widetilde{u}_{n-1}^{\sin }\in \Pi _{n-1} $  \  и
\[
    \widetilde{u}_{n-1}^{\sin } ( x,\underline{l} )
    \dashv \vdash SIN_{n-1} ( x )
    \wedge \widetilde{u}_{n-1} ( x ,
    \underline{l} ) .
\]
\hfill {} \\
Монотонная и субнедостижимо универсальная для класса \
$\Pi_{n}^{\dashv \vdash} $ \ формулы вводятся двойственным
образом.
\\
\hfill {} \\
{\em 3)}\quad Булевы значения
\[
    A_{\varphi}^{\vartriangleleft \alpha _{1}} ( \alpha  ), \
    \Delta_{\varphi }^{\vartriangleleft \alpha _{1}} ( \alpha  )
    \mbox{\it \ и спектр \ }
    \mathbf{S}_{\varphi }^{\vartriangleleft \alpha_{1}}
\]
формулы  \ $\varphi = \widetilde{u}_{n}^{\sin }$ \ и его проекции
(см. определение 2.4, где \ $\trianglelefteq^{\underline{l}}$ \
следует заменить на \ $\leq$) будут называться субнедостижимо
универсальными уровня \ $n$ \ ниже \ $\alpha_{1}$ \ и обозначаться
соответственно через

\[
    \widetilde{A}_{n}^{\sin \vartriangleleft \alpha _{1}} ( \alpha ) ,
    \quad \widetilde{\Delta }_{n}^{\sin \vartriangleleft \alpha _{1}}
    ( \alpha  ) ,\quad \widetilde{\mathbf{S}}_{n}^{\sin \vartriangleleft
    \alpha _{1}}.
\]
\hfill {} \\
{\em 4)}\quad Если \ $( \alpha ,\Delta  ) \in
\widetilde{\mathbf{S}}_{n}^{\sin \vartriangleleft \alpha _{1}}$, \
то \ $\alpha $ \ называется ординалом скачка этих формулы и
спектров, а \ $\Delta $ \ называется их булевым значением ниже \
$\alpha _{1}$.
\hfill {} \\
{\em 5)}\quad Кардинал \ $\alpha_{1}$ \ называется носителем этих
спектров. \hspace*{\fill} $\dashv$
\end{definition}
Из этого определения и лемм 3.4~~2) ( для \ $n-1$ \ вместо \ $n$),
2.5~~3), 2.7 и 3.8 непосредственно вытекают следующие простые
леммы:

\begin{lemma}
\label{3.10.} \quad \\
{\em 1)}\qquad \qquad \qquad \qquad \qquad $u_{n}^{\Sigma} (
\underline{l} ) \dashv \vdash \widetilde{u}_{n}^{\sin } (
\underline{l} )$~.
\\
\\
{\em2)} \qquad   $ dom \ \bigl ( \widetilde{\mathbf{S}}_{n}^{\sin
\vartriangleleft \alpha _{1}} \bigr ) \subseteq
SIN_{n-1}^{<\alpha_{1}} \cap dom  \bigl
(\mathbf{S}_{n}^{\triangleleft \alpha_{1}} \bigr )$~.
\\
{\em 3)} Пусть
\[
    \alpha \in dom \bigl ( \mathbf{S}_{n}^{\triangleleft \alpha_{1}} \bigr )
    \mbox{\it \quad и \quad} \alpha^{\prime} = min  \left \{ \alpha^{\prime \prime} \in
    SIN_{n-1}^{< \alpha_{1}} : \alpha^{\prime \prime} > \alpha
    \right \}~,
\]
тогда
\[
    \alpha^{\prime} \in dom \bigl (
    \widetilde{\mathbf{S}}_{n}^{\sin \triangleleft \alpha_{1}}
    \bigr )~.
\]
{\em 4)} Пусть
\[
    \alpha \leq \alpha_{1} \mbox{\it \quad - предельный кардинал в \quad}
    SIN_{n-1}^{< \alpha_{1}}~,
\]
тогда
\[
    \sup dom \bigl ( \widetilde{\mathbf{S}}_{n}^{
    \sin \triangleleft \alpha_{1}} | \alpha \bigr ) =
    \sup dom \bigl ( \mathbf{S}_n^{\triangleleft
    \alpha_{1}} | \alpha \bigr ).
\]
\end{lemma}

\vspace{-6pt}

\begin{lemma} \label{3.11.} \hfill {} \\
\hspace*{1em} Пусть
\[
    \alpha_{2} \in SIN_{n-2}^{< \alpha_{1}} \quad \mbox{\it и}\quad
    \alpha_{0} =\sup \bigl ( SIN_{n-1}^{<\alpha _{1}}\cap \alpha_{2} \bigr ),
\]
тогда
\[
    \widetilde{\mathbf{S}}_{n}^{\sin \vartriangleleft
    \alpha _{2}}  | \alpha_{0} =\widetilde{\mathbf{S}}_{n}^{\sin
    \vartriangleleft \alpha _{1}} | \alpha_{0} \quad .
\]
\end{lemma}
\textit{Доказательство.} \ Пусть  \  $\alpha<\alpha_{0},(\alpha,
\Delta ) \in \widetilde{\mathbf{S}}_{n}^{\sin \triangleleft
\alpha_{2} } $ \ и поэтому
\[
    \Delta = \widetilde{\Delta}_{n}^{\sin
    \triangleleft \alpha_{2} } (\alpha)>0.
\]
Для каждого  \  $\alpha^{\prime} \leq \alpha $ \  благодаря лемме
3.10~2) и лемме 3.8 (использованной для  \  $n-1 $ \ вместо \ $n
$)
\[
    SIN_{n-1}^{< \alpha_{2}} (\alpha^{\prime})
    \longleftrightarrow SIN_{n-1}^{< \alpha_{1}}
    (\alpha^{\prime})
\]
и поэтому  \ $\vartriangleleft^{\underline{l}}$~-ограничение
ординалом \ $\alpha_{2} $ \  в утверждении  \
$\widetilde{u}_{n-1}^{ \triangleleft \alpha_{2}} \left
(F^{\underline{l}} (\alpha^{\prime}), \underline{l}) \right )$ \
можно заменить на \
$\vartriangleleft^{\underline{l}}$~-ограничение ординалом \
$\alpha_{1} $. \ Это влечёт
\[
    \Delta= \widetilde{\Delta}_{n}^{\sin
    \triangleleft \alpha_{1} } (\alpha)
    \quad \mbox{\it и} \quad (\alpha, \Delta ) \in
    \widetilde{\mathbf{S}}_{n}^{\sin \triangleleft \alpha_{1} }.
\]
Обратное рассуждение завершает доказательство. \hfill $\dashv$\\

\quad \\

Наша цель -- отыскать способ ``сравнения'' универсальных спектров
друг с другом на \textit{различных} носителях \ $\alpha_1$, \
расположенных конфинально \ $k$, \ чтобы построить
\textit{монотонную} матричную функцию. С этой целью естественно
использовать значения гёделевой функции \ $Od$ \ на таких
спектрах.
\\
Также естественно найти некоторые оценки ``информационной
сложности'' таких спектров посредством оценок их порядковых типов.
Однако из лемм \ref{2.7.}, \ref{2.8.}, \ref{3.10.} \ следует, что,
например, спектры \ $dom \left (
\widetilde{\mathbf{S}}_{n-1}^{\sin \triangleleft \alpha_{1}} (
\alpha  ) \right ) $ \ накапливаются, когда \ $\alpha $, \
$\alpha_{1}$ \ возрастают и поэтому их порядковые типы возрастают
до \ $k $; \ кроме того они замкнуты относительно \ $\Pi
_{n-2}$-функций и т.д. Спектры \ $\widetilde{\mathbf{S}}_{n}^{\sin
\triangleleft \alpha_{1}}$ \ обладают аналогичными свойствами.
\\
Поэтому требуемое сравнение таких спектров весьма сложно получить
естественным должным образом, так как они \textit{``слишком
отличаются''} друг от друга на произвольно больших носителях \
$\alpha_1$.
\\
Значит, ничего не остаётся кроме использования в дальнейшем
\textit{спектров, редуцированных к некоторому фиксированному
кардиналу}, и, далее, \textit{редуцированных матриц}.

\newpage
\quad 

\newpage
\section{Редуцированные спектры}
\setcounter{equation}{0}

Здесь мы начинаем формировать основной материал для построения
матричных функций -- редуцированные матрицы.
\\
С этой целью мы рассмотрим необходимые предварительные конструкции
-- редуцированные спектры.
\\

Для ординала \ $\chi \leq k$ \ пусть \ $P_{\chi }$ \ обозначает
множество \ $\{ p\in P:dom ( p ) \subseteq \chi \} $ \ и \
$B_{\chi }$ \ обозначает подалгебру \ $B$, \ порождённую \
$P_{\chi}$ \ в \  $L_{k}$. \  Для каждого \ $A\in B$ \ введём
множество
\[
    A \lceil \chi =\{ p\in P_{\chi }:\exists q~ ( p=\left. q\right|
    \chi \wedge q\leq A ) \}
\]
которое называется значением \ $A$, \ редуцированным к \ $\chi $.
\ Хорошо известно (см.~\cite{Jech}), что
\[
    B_{\chi }=\left \{ \sum X:X\subseteq P_{\chi } \right \},
\]
и поэтому каждое \ $A\in B_{\chi }$  \ совпадает с \ $ \sum
A\lceil \chi $. \ Поэтому будем отождествлять каждое \ $A\in
B_{\chi }$ \ с его редуцированным значением \ $A\lceil \chi $; \ и
снова следует отметить, что именно поэтому каждое значение \ $A
\in B_\chi$ \ это множество в \ $L_k$, не класс, и \ $B_\chi$ \
рассматривается как множество таких значений.

\begin{definition} \label{4.1.} \hfill {} \\
\hspace*{1em} Пусть \  $\chi \leq k$,\quad $\alpha_{1}\leq k$.
\hfill {} \\
{\em 1)} \quad Для всякого \  $\alpha <\alpha _{1}$  \ введём
булевы значения и спектр:
\hfill {} \\
\[
    \mathbf{\ }\widetilde{A}_{n}^{\sin \vartriangleleft \alpha _{1}}
    (\alpha  ) \lceil \chi \mathbf{;~}\widetilde{\Delta }_{n}^{\sin
    \vartriangleleft \alpha _{1}} ( \alpha  ) \overline{\lceil}\chi =
    \widetilde{A}_{n}^{\sin \vartriangleleft \alpha _{1}} ( \alpha  )
    \lceil\chi \mathbf{-}\sum_{\alpha ^{\prime }<\alpha }
    \widetilde{A}_{n}^{\sin \vartriangleleft \alpha _{1}}
    ( \alpha ^{\prime } ) \lceil \chi ;
\]
\[
    \widetilde{\mathbf{S}}_{n}^{\sin \vartriangleleft \alpha _{1}}
    \overline{\overline{\lceil}}\chi =\{
    ( \alpha ;~\widetilde{\Delta }_{n}^{\sin \vartriangleleft \alpha _{1}}
    ( \alpha  ) \overline{\lceil}\chi  ) :\alpha <\alpha _{1}\wedge
    \widetilde{\Delta }_{n}^{\sin \vartriangleleft \alpha _{1}} ( \alpha  )
    \overline{\lceil}\chi >0\}  .
\]
\hfill {} \\
{\em 2)} \quad Эти значения, спектр и его первая и вторая проекции
называются субнедостижимо универсальными редуцированными к \ $\chi
$ \  уровня \ $n$\quad ниже \ $\alpha _{1}$.
\hfill {} \\
{\em 3)}\quad Если \ $( \alpha ,\Delta  ) \in
\widetilde{\mathbf{S}}_{n}^{\sin \vartriangleleft \alpha _{1}}
\overline{\overline{\lceil}} \chi$, \ то \ $\alpha $ \ называется
ординалом скачка этих спектров, а \ $\Delta $ \ называется их
булевым значением редуцированным к \ $\chi $ \ ниже \
$\alpha_{1}$.
\hfill {} \\
{\em 4)}\quad Кардинал \  $\alpha _{1}$ \ называется носителем
этих спектров.  \hspace*{\fill} $\dashv$
\end{definition}
Аналогичным образом можно ввести многомерные редуцированные
спектры. \label{c4}
\endnote{
\ стр. \pageref{c4}. \ С этой целью следует все булевы значения
многомерных спектров редуцировать к данному ординалу \ $\chi$. \
Например, двухмерные спектры (см. примечание 2) ) трансформируются
в их редуцированные формы:
\[
    \mathbf{S}_{\varphi, 2 }^{\vartriangleleft \alpha _{1}}
    ( \overrightarrow{a}) \overline{\overline{\lceil}} \chi = \{
    ( \alpha , \beta, \Delta_{\varphi}^{\vartriangleleft
    \alpha_{1}}  ( \alpha ,\overrightarrow{a} ) \overline{\lceil} \chi,
    \Delta_{\varphi_1^{\prime}}^{\vartriangleleft \alpha_{1}}
    (\alpha , \beta, \overrightarrow{a} ) \overline{\lceil} \chi  ):
\]
\[
    : \alpha,\beta <\alpha_{1} \wedge \ \Delta_{\varphi }^{\vartriangleleft \alpha_1}
    (\alpha , \overrightarrow{a}) \overline{\lceil} \chi > 0 \
    \wedge \
    \Delta_{\varphi_1^{\prime} }^{\vartriangleleft \alpha_1}
    (\alpha, \beta, \overrightarrow{a}) \overline{\lceil} \chi > 0  \}.
\]
\vspace{6pt}
} 

\noindent Далее всегда предполагается, что \ $\chi $ \ замкнут
относительно функции пары; если $\chi =k$, \ то все упоминания о \
$\chi $ \ будут опускаться.
\\
Следующие две леммы нетрудно получить из определений \ref{3.9.},
\ref{4.1.} и лемм \ref{3.10.}, \ref{3.11.}:

\begin{lemma}
\label{4.2.}
\[
    dom \left ( \widetilde{\mathbf{S}}_{n}^{\sin \vartriangleleft
    \alpha _{1}}\overline{\overline{\lceil}}\chi \right )
    \subseteq dom  \left (  \widetilde{
    \mathbf{S}}_{n}^{\sin \vartriangleleft \alpha _{1}} \right )
    \subseteq SIN_{n-1}^{<\alpha _{1}} \cap
    dom \Bigl( \mathbf{S}_{n}^{\vartriangleleft \alpha_1 } \Bigr).
\]
\end{lemma}

\vspace{-18pt}

\begin{lemma} \label{4.3.} \hfill {} \\
\hspace*{1em} Пусть
\[
    \alpha _{2}\in SIN_{n-2}^{<\alpha _{1}} \mbox{\it \quad и \quad}
    \alpha_{0}=\sup  ( SIN_{n-1}^{<\alpha _{1}}\cap \alpha_{2} ),
\]
тогда
\[
    ( \widetilde{\mathbf{S}}_{n}^{\sin \vartriangleleft
    \alpha _{2}} \overline{\overline{\lceil }}\chi  )
    | \alpha_{0} = ( \widetilde{\mathbf{S}}_{n}^{\sin
    \vartriangleleft \alpha _{1}}\overline{\overline{\lceil }}\chi)
    | \alpha_{0} \quad.
\]
\hspace*{\fill} $\dashv$

\end{lemma}
Следующая лемма аналогична лемме \ref{2.5.} и следует из
определений:

\begin{lemma} \label{4.4.} \hfill {} \\
\hspace*{1em} Пусть \ $\alpha <\alpha _{1}$,\quad $ \chi \leq k$, \  тогда: \\
\hfill {} \\
{\em 1)}\quad $ \sup dom  \left (
\widetilde{\mathbf{S}}_{n}^{\sin \vartriangleleft
\alpha_{1}}\overline{\overline{\lceil }}\chi \right) < k$.
\\
\quad \\
\medskip
{\em 2)}\ $\widetilde{\mathbf{S}}_{n}^{\sin \vartriangleleft
\alpha_{1}}\overline{\overline{\lceil }}\chi $, \ $dom  \left (
\widetilde{ \mathbf{S}}_{n}^{\sin \vartriangleleft \alpha
_{1}}\overline{\overline{\lceil}}\chi  \right ) $
$\Delta_{n}$-определимы, а
\\
$rng  \left (\widetilde{\mathbf{S}}_{n}^{\sin \vartriangleleft
\alpha _{1}}\overline{\overline{\lceil}}\chi   \right )$ \ \ $
\Sigma _{n}$-определима \medskip в  \  $L_{k}$ \  для \ $\alpha
_{1}=k$.
\\
Для \  $\alpha _{1}<k$ \ эти спектры $\Delta_{1}$-определимы;
\\
\quad \\
{\em 3)}\quad $\alpha \in dom  \left (
\widetilde{\mathbf{S}}_{n}^{\sin \vartriangleleft \alpha
_{1}}\overline{\overline{\lceil }}\chi  \right ) $ \ если
существует \ $\mathfrak{M}$-генерическая функция
\[
    l\stackrel{\ast }{\in }\widetilde{\Delta }_{n}^{\sin
    \vartriangleleft \alpha _{1}} (\alpha) \overline{\lceil}\chi ;
\]
\hfill {} \\
{\em 4)}\quad \quad \quad \quad \quad \quad \quad  $\left\|
\widetilde{u}_{n}^{\sin \vartriangleleft \alpha_{1}}
( \underline{l} ) \right\| \lceil \chi =\sum rng  \left ( \widetilde{%
\mathbf{S}}_{n}^{\sin \vartriangleleft \alpha _{1}}\overline{\overline{%
\lceil }}\chi  \right ) $. \hspace*{\fill} $\dashv$
\end{lemma}
Как и в лемме \ref{2.5.}~3) утверждение 3) здесь делает возможным
обнаружение ординалов скачка \ $\alpha $ \ с помощью генерических
функций \ $l$; \ эта техника используется ниже в доказательстве
леммы \ref{4.6.}. Можно обойтись без этого, используя вместо \ $l$
\ условия \ $p\subset l$ \ с достаточно длинной областью
определения.
\\

Теперь обратимся к обсуждению спектральных порядковых типов. Если
\ $X$ \ это вполне упорядоченное множество, то его порядковый тип
обозначается через \ $OT ( X ) $; \
\\
если \  $X$ \  это функция со вполне упорядоченной областью
определения, то мы полагаем, что его
порядковый тип это ординал \  $OT ( X ) =OT ( dom ( X )  ) $.\\

Грубая верхняя оценка спектральных типов непосредственно вытекает
из леммы \ref{1.1.}, \ $\left| P_{\chi }\right| =\left| \chi
\right| $ \ и \  $GCH$ \ в \ $L_{k}$~:
\begin{lemma} \label{4.5.}
\[
    OT ( \widetilde{\mathbf{S}}_{n}^{\sin \vartriangleleft
    \alpha_{1}}\overline{\overline{\lceil }}\chi  ) < \chi ^{+}\ .
\]
\end{lemma}

Теперь обсудим оценки снизу таких типов. Здесь возникает лемма,
весьма существенная для доказательства основной теоремы. Она
показывает, что как только ординал \ $ \delta < \chi ^{+}$ \
определим через некоторые ординалы скачка субнедостижимо
универсального редуцированного к \ $\chi$ \ спектра, то порядковый
тип этого спектра превосходит \ $\delta$ \ при некоторых
естественных условиях.

Мы будем использовать здесь и в дальнейшем метод рассуждений,
который может быть назван \textit{расщепляющим методом}. В своём
простейшем варианте он состоит в расщеплении некоторого
рассматриваемого значения
\[
    \Delta = \widetilde{\mathbf{S}}_{n}^{\sin \triangleleft
    \alpha_1} \overline{\overline{\lceil}} \chi (\alpha)
\]
(или нескольких таких значений) на последовательность его частей,
которые затем приписываются последовательно расположенным \
$SIN_{n-1}^{< \alpha_1} $-кардиналам соответственно, после
некоторых их незначительных трансформаций. После этого такие
кардиналы становятся \textit{кардиналами предскачка} этого
спектра. Предварительно для этого кардинал \ $\alpha $ \ должен
быть \textit{фиксирован} некоторой функцией  \  $l $, \ то есть \
$\mathfrak{M} $-генерической функцией \ $l\stackrel{\ast }{\in }
\Delta $.

Для некоторого удобства будут использоваться традиционные
обозначения ординальных промежутков с концами \ $\alpha
_{1}<\alpha _{2}$:
\\
$\left[ \alpha _{1},\alpha _{2}\right[ =\alpha _{2}-\alpha _{1}$;
\quad $\left] \alpha _{1},\alpha _{2}\right[ =\alpha _{2}- (
\alpha _{1}+1 ) $; \quad $\left[ \alpha _{1},\alpha _{2}\right] =
( \alpha _{2}+1 ) -\alpha _{1}$; \quad $\left] \alpha _{1},\alpha
_{2}\right] = ( \alpha _{2}+1 ) - (\alpha _{1} + 1)$ \quad (здесь
и далее \ $\alpha _{1}, \alpha _{2}$ \ это множества меньших
ординалов).
\begin{lemma}
\label{4.6.} {\em (О спектральном типе)} \hfill {} \\
\hspace*{1em} Пусть ординалы \ $ \overline{\delta
}$,~$\overline{\chi }$,~$\overline{\alpha
}_{0}$,~$\overline{\alpha}_1$ \ такие, что:
\\
\quad \\
(i) \quad $\overline{\delta }<\overline{\chi
}^{+}<\overline{\alpha}_{0}<\overline{\alpha}_1\leq k$~; \\
\quad \\
(ii) \quad $SIN_{n-2} ( \overline{\alpha}_1 ) \wedge OT (
SIN_{n-1}^{<\overline{\alpha}_1 } ) =\overline{\alpha}_1$~;
\\
\quad \\
(iii) \quad $SIN_{n-1}^{<\overline{\alpha}_1 } ( \overline{\chi }
) \wedge \sup dom  \left ( \widetilde{\mathbf{S}}_{n}^{\sin
\vartriangleleft
\overline{\chi }}  \right ) =\overline{\chi }$~; \\
\quad \\
(iv) \quad $\sum rng  \left ( \widetilde{\mathbf{S}}_{n}^{\sin
\vartriangleleft \overline{\chi }}  \right ) \in
B_{\overline{\chi }}$~;\\
\quad \\
(v) \quad $\overline{\alpha }_{0}\in dom  \left (
\widetilde{\mathbf{S}} _{n}^{\sin \vartriangleleft
\overline{\alpha}_1}\overline{\overline{\lceil }}
\overline{\chi}  \right ) $~;\\
\quad \\
(vi) \quad $\overline{\delta }$ \ определён в  \ $ L_{k} $ \ через
ординалы \  $\overline{\alpha }_{0}$, $\overline{\chi }$ \
формулой класса \  $\Sigma
_{n-2}\cup \Pi _{n-2}$~.\\
\quad \\
Тогда $\quad \quad \quad \quad \quad \quad \quad \quad \quad
\quad\overline{\delta }<OT ( \widetilde{\mathbf{S}}_{n}^{\sin
\vartriangleleft \overline{\alpha}_1}\overline{\overline{\lceil }}%
\overline{\chi}  ) \quad.$
\end{lemma}
\textit{Доказательство.} \ Введём следующие формулы, описывающие
существенные аспекты ситуации ниже \ $\overline{\alpha}_1$.
\\
По условию $(vi)$ существует  \ $\Sigma_{n-2}\cup
\Pi_{n-2}$-формула \  $\psi _{0} ( \alpha _{0},\chi ,\delta ) $, \
определяющая \  $\overline{\delta } $ \ через \ $\alpha _{0}
=\overline{\alpha }_{0}$, $\chi =\overline{\chi }$, \ то есть \
$\delta =\overline{\delta }$ \  это единственный ординал,
выполняющий \ $\psi _{0} ( \overline{\alpha }_{0},\overline{\chi
},\delta  ) $ \ в \ $L_{k}$.
\\
Благодаря условим  $(i)$,  $(iii)$, лемме 3.5 и минимальности
кардинала \ $k$, \ $\overline{\chi }$ \ это сингулярный кардинал,
а тогда по лемме~1.3~2) для каждой \ $\mathfrak{M}$-генерической
функции \ $\left. l\right| \overline{\chi } $ \ на \ $\overline{
\chi }$ ординалы \ $\overline{\chi }$, \ $\overline{\delta }$ \
счётны в \ $L_{k}\left[ \left. l\right| \overline{\chi }\right] $.
\ Обозначим через\quad $ ( f:\omega _{0}\longrightarrow \delta +1
) $ \ формулу
\begin{center}
{``$f$ \ отображает \ $\omega _{0}$ \  на \ $ \delta +1$''~~.}
\end{center}
В алгебре \ $B_{\overline{\chi }}$ \ имеется булево значение
\[
    \left\| \exists f\  ( f:\omega _{0}\longrightarrow
    \delta +1 ) \right\|_{\overline{\chi }}=1
\]
и поэтому сществует имя \ $\underline{f }\in
L_{k}^{B_{\overline{\chi }}}$ \ для которого
\[
    \| \underline{f}:\omega _{0}\longrightarrow
    \delta +1 \| _{\overline{\chi } }=1
\]
в \ $B_{\overline{\chi}}$ \ (см.~\cite{Jech}). Поэтому существует
ординал \ $\overline{\beta }$, \ определимый через \ $\alpha _{0}=
\overline{\alpha }_{0}$, \ $\chi =\overline{\chi }$ \ следующей
формулой, которую обозначим через \ $\psi_{1} ( \alpha _{0},\chi
,\beta  ) $:
\\
\quad \\
\hspace*{1em} $ \exists \delta<\alpha _{0}( ^{^{{}}}\psi _{0} (
\alpha _{0},\chi ,\delta ) \wedge$
$$
\qquad \wedge \beta =\min \{ \beta ^{\prime }:F ( \beta ^{\prime }
) \in L^{B_{\chi }}\wedge \left\| F ( \beta ^{\prime } ) :\omega
_{0}\longrightarrow \delta +1\right\| _{\chi }=1\}  ) .
$$
Из условий $(ii)$, $(v)$ и лемм~\ref{3.8.},~\ref{4.2.} \ следует
$\overline{\alpha }_{0} \in SIN_{n-2}$,\ а тогда по условию $(i)$
и лемме 3.5 выполняется \ $\overline{\beta
}<\overline{\alpha}_{0}$. \ Далее, из $(ii)$, $(v)$ и той же леммы
~\ref{4.2.} следует \ $ \overline{\alpha }_{0}\in dom \left (
\widetilde{\mathbf{S}}_{n}^{\sin \vartriangleleft \overline{\alpha
}_1}  \right ) $ \  и поэтому можно определить ординал скачка
\[
    \overline{\alpha}_0^\prime = \min  ( \left[ \overline{\chi },
    \overline{\alpha }_{0}\right] \cap dom  \left
    ( \widetilde{\mathbf{S}}_{n}^{\sin \vartriangleleft
    \overline{\alpha}_1}  \right ) ) .
\]
Так как \ $\overline{\chi },\overline{\alpha}_0^\prime \in
SIN_{n-1}^{<\overline{\alpha}_1}$, \ то из $(iii)$, $(iv)$ следует
\[
    \overline{\alpha}_0^\prime \in dom  \left
    ( \widetilde{\mathbf{S}}_{n}^{\sin \vartriangleleft
    \overline{\alpha}_1} \overline{\overline{\lceil }}\overline{\chi}
    \right )
\]
и кардинал \  $\overline{\chi }$ \  определяется через \
$\alpha_0^\prime = \overline{\alpha}_0^\prime$ \ формулой
\[
    \chi =\sup dom  \left ( \widetilde{\mathbf{S}}_{n}^{\sin
    \vartriangleleft \alpha_0^\prime } \right ),
\]
которую обозначим через \ $\psi _{2} ( \alpha_0^\prime , \chi)$.
\\
Будем полагать, что формулы \ $\psi _{i}$,\quad
$i=\overline{0,2}$, \  преобразованы в их предуниверсальные формы;
далее они будут использованы в генерических расширениях \ $L_{k}$
\ и в подобных случаях они будут релятивизированы к
конструктивному классу и обозначены через \ $\psi _{i}^{L}$,
$i=\overline{0,2}$.

Теперь начинает действовать \textit{расщепляющий метод}. Для этого
необходимо использовать прямые произведения функций. Пусть \
$l_{0},...,l_{m}\in {}^{k}k$; \ на \ $k$ \ определяется функция \
$l=l_{0}\oplus ...\oplus l_{m}$ \ следующим образам: для каждого
ординала \  $\alpha <k, \quad \alpha =\alpha _{0}+ ( m+1 ) i+j$, \
где ординал \ $\alpha _{0}$ \ предельный и \ $i\in \omega _{0}$,
\quad $j \leq m $, \ она имеет значение \ $l ( \alpha ) =l_{j} (
\alpha _{0}+i ) $. \ Обозначим через \ $ ( , ) _{j}^{m}$ \
операцию, реконструирующую \ $l_{j}$ \ по \ $l$, то есть \ $l_{j}=
( l ) _{j}^{m}$. \ Известно (см. Соловай ~\cite{Solovay}), что \
$l$ \ это \ $\mathfrak{M}$-генерическая функция на \ $k$ \
\textit{если} \ $l_{0}$ \  это \ $\mathfrak{M}$-генерическая
функция и \ $l_{j}$ \ -- \ $\mathfrak{M}\left[ l_{0}\oplus
...\oplus l_{j-1}\right] $-генерическая функция на \ $k$, \
$j=\overline{1,m}.$
\\
Мы возьмём здесь генерические \ $l_0$, \ $l_1$, \ чтобы
зафиксировать булевы значения спектра на кардиналах скачка \
$\overline{\alpha}_0$, \ $\overline{\alpha}_0^\prime$, \ то есть
возьмём \ $l_0$, \ $l_1$ \ такие, что
\[
    l_0\stackrel{\ast }{\in }\widetilde{\Delta }_{n}^{\sin
    \vartriangleleft \overline{\alpha}_1} (\overline{\alpha}_0)
    \overline{\lceil}\chi, \quad \quad
    l_1\stackrel{\ast }{\in }\widetilde{\Delta }_{n}^{\sin
    \vartriangleleft \overline{\alpha}_1} (\overline{\alpha}_0^\prime)
    \overline{\lceil}\chi;
\]
их существование следует из леммы~4.4~3).
\\
Теперь рассмотрим формулу \ $\varphi =\exists \alpha \ \varphi
_{1} ( \alpha ,\underline{l} ) $, \ собирающую всю полученную
информацию о ситуации ниже \ $\overline{\alpha}_1$ \ и задающую
расщеплённые булевы значения на ординалах скачка \
$\overline{\alpha}_0$, \ $\overline{\alpha}_0^\prime$, \ где \
$\varphi _{1} $ \ это следующая формула:
\begin{eqnarray*}
    SIN_{n-1} ( \alpha  ) \wedge \quad \qquad \qquad \qquad \qquad
    \qquad \qquad \qquad \qquad \qquad  \qquad \qquad
\\
    \wedge \exists \alpha_0 <\alpha \exists \alpha_0^\prime
    \leq \alpha _{0}\exists \chi \leq \alpha_0^\prime \exists \beta
    <\alpha _{0}\exists i\in \omega _{0}\exists \delta <\alpha
    _{0}\exists \delta _{i}\leq \delta \exists p\in P_{\chi }
\\
    \left[ SIN_{n-1}^{<\alpha } ( \alpha _{0} ) \wedge
    SIN_{n-1}^{<\alpha } ( \alpha_0^\prime ) \wedge \widetilde{u}%
    _{n-1}^{\vartriangleleft \alpha } ( \alpha _{0,} ( \underline{l})
    _{0}^{3} ) \wedge \widetilde{u}_{n-1}^{\vartriangleleft \alpha}
    ( \alpha_0^\prime, ( \underline{l} ) _{1}^{3} ) \right. \wedge
\\
    \wedge \psi _{0}^{L} ( \alpha _{0},\chi ,\delta  ) \wedge
    \psi_{1}^{L} ( \alpha_0,\chi ,\beta  ) \wedge \psi _{2}^{L}
    (\alpha _{0}^\prime,\chi  ) \wedge
\\
    \wedge ( \underline{l} )
    _{2}^{3} ( \omega _{0} ) =i\wedge p\subseteq  ( \underline{l}%
     ) _{3}^{3}\wedge p\leq \left\| F ( \beta  )  ( i )
     =\delta _{i}\right\| _{\chi }\wedge
\\
    \left . \wedge OT\{ \alpha ^{\prime }<\alpha :SIN_{n-1}^{<\alpha }
    ( \alpha ^{\prime } ) \} = \alpha_{0} + \delta _{i}^{^{{}}} \right].
\end{eqnarray*}
\hfill {} \\
 Эта формула содержится
в \  $\Pi _{n-1}$, \ так как все её переменные в квадратных
скобках ограничены \ $SIN_{n-1}$-переменной \ $ \alpha $. \ Для
каждого \ $\delta \leq \overline{\delta }$ \ обозначим через \
$\alpha _{\delta }$ \ кардинал \ $\alpha < \overline{\alpha }_1$ \
такой, что \vspace{6pt}
\begin{equation}
\label{e4.1} SIN_{n-1}^{<\overline{\alpha}_1 } ( \alpha  ) \wedge
OT\{ \alpha ^{\prime }<\alpha :SIN_{n-1}^{<\alpha } ( \alpha
^{\prime } ) \} = \overline{\alpha}_{0} + \delta ~;
\end{equation}
\vspace{0pt}

\noindent такой кардинал существует благодаря условию $(ii)$ и
лемме~3.5. Покажем, что каждый \ $\alpha _{\delta }$ \ это ординал
скачка \ $\varphi $ \ ниже \ $ \overline{\alpha }_1$, \
посредством \textit{расщепления} упомянутых булевых значений. С
этой целью рассмотрим \ $\mathfrak{M}$-генерическую функцию на
 \ $k$ \ \ $l=l_{0} \oplus l_{1}
\oplus l_{2} \oplus l_{3}$ \ такую, что

\vspace{6pt}
\begin{equation}
\label{e4.2} L_{k}\left[ l_{0}\right] \vDash
\widetilde{u}_{n-1}^{\vartriangleleft \overline{\alpha}_1} (
\overline{\alpha }_{0,}l_{0} )\quad ;
\end{equation}

\begin{equation}
\label{e4.3}
l_{0}\stackrel{\ast }{\notin }\left\| \exists \alpha <\overline{\alpha }%
_{0} ( SIN_{n-1}^{\vartriangleleft \overline{\alpha}_1} ( \alpha )
\wedge \widetilde{u}_{n-1}^{\vartriangleleft \overline{\alpha}_1 }
( \alpha _{,}\underline{l} )  ) \right\| \lceil \overline{\chi }
\quad;
\end{equation}

\begin{equation}
\label{e4.4} L_{k}\left[ l_{1}\right] \vDash
\widetilde{u}_{n-1}^{\vartriangleleft \overline{\alpha}_1} (
\overline{\alpha}_0^\prime, l_{1} )\quad ;
\end{equation}

\begin{equation}
\label{e4.5} l_{1}\stackrel{\ast }{\notin }\left\| \exists \alpha
<\overline{\alpha}_0^\prime ( SIN_{n-1}^{\vartriangleleft
\overline{\alpha}_1} ( \alpha ) \wedge
\widetilde{u}_{n-1}^{\vartriangleleft \overline{\alpha}_1 } (
\alpha _{,}\underline{l} )  ) \right\| \lceil \overline{\chi
}\quad;
\end{equation}

\begin{equation}
\label{e4.6} l_{2} ( \omega _{0} ) =i_{\delta }\quad;
\end{equation}

\begin{equation}
\label{e4.7} \exists p\in P_{\overline{\chi }} ( p\subset
l_{3}\wedge p\leq \left\| F ( \overline{\beta } )  ( i_{\delta }
) =\delta \right\| _{\overline{\chi }} )
\end{equation}
\vspace{0pt}

\noindent для соответствующего \ $i_{\delta } \in \omega_0$; \
здесь существование \ $l_{2}$, \ $l_{3}$ \ очевидно. Из
(\ref{e4.2}), (\ref{e4.4}), (\ref{e4.6}), (\ref{e4.7}) и
(\ref{e4.1}) следует

\vspace{6pt}
\begin{equation}
\label{e4.8} l\stackrel{\ast }{\in }\left\| \varphi
_{1}^{\vartriangleleft \overline{\alpha}_1} ( \alpha _{\delta
},\underline{l} ) \right\| \lceil \overline{\chi }\quad.
\end{equation}
\vspace{0pt}

\noindent Тогда из (\ref{e4.3}), (\ref{e4.5}) нетрудно заключить,
что

\[
    l\stackrel{\ast }{\notin }\left\| \exists \alpha <
    \alpha_{\delta }\ \varphi_{1}^{\vartriangleleft
    \overline{\alpha}_1} ( \alpha ,\underline{l} )
     \right\| \lceil \overline{\chi }\quad.
\]
\vspace{0pt}

\noindent Вместе с (\ref{e4.8}) это влечёт \quad \(
    \alpha_{\delta } \in dom \left (
    \mathbf{S}_{\varphi }^{\vartriangleleft \overline{\alpha}_1}
    \overline{\overline{\lceil} }\overline{\chi } \right )
\). \ После этого с помощью условия~$(ii)$ нетрудно видеть, что
для каждого \ \(
    \delta \leq \overline{\delta }\quad  dom \left (
    \widetilde{\mathbf{S}}_{n}^{\sin \vartriangleleft
    \overline{\alpha}_1} \overline{\overline{\lceil }}
    \overline{\chi }\right )
\) \ содержит кардинал, следующий за \ $\alpha _{\delta}$ \ в \
$SIN_{n-1}^{<\overline{\alpha}_1}$. \hfill $\dashv$
\\

Незначительно усложняя это рассуждение можно доказать эту лемму
для случаев, когда ординал \  $\overline{\delta } $ \ определён
через \ $\overline{\chi} $ \ и несколько кардиналов скачка
\[
    \overline{\alpha}_{1}, \dots, \overline{\alpha}_{m} \in dom  \left
    (\widetilde{\mathbf{S}}^{\sin \triangleleft \overline{\alpha}_1}
    \overline{\overline{\lceil}} \overline{\chi} \right ) ,
\]
или кардиналов предскачка. Легко также видеть, что эта лемма
справедлива и в том случае, когда ординал \  $\overline{\delta } $
\ определим формулой \emph{любого уровня,} но ограниченной одним
из этих кардиналов.
\\
\noindent Лемма \ref{4.6.} допускает разнообразные усиления и
версии, однако они безразличны для дальнейшего и поэтому
опускаются. \label{c5}
\endnote{
\ стр. \pageref{c5}. \ Например, можно потребовать определимость \
$\overline{\delta }$ \ не в \ $ L_{k} $, \ а в некотором \
$L_{k}\left[ l\right] $; \ можно ослабить условие (ii) до условия

\[
    SIN_{n-2} ( \overline{\alpha}_1 ) \wedge OT ( \;
    ]\overline{\alpha }_0, \overline{\alpha}_1 [ \; \cap \;
    SIN_{n-1}^{< \overline{\alpha}_1 } ) \geq \overline{\chi}^{+} ,
\]

\noindent также можно существенно ослабить требования, наложенные
на \ $\overline{\chi}$ \ и т.д.
\\
\quad \\
} 
\hfill {} \\
\hfill {} \\
В этом доказательстве интерпретируемые в \ $L_{k}$ \ формулы
снабжались верхним индексом \ $L$ \ при переходе к расширениям \
$L_{k}$. \ В дальнейшем этот индекс будет опускаться, если
соответствующий переход подразумевается контекстом.
\hfill {} \\
\hfill {} \\

Лемма~4.6 о спектральном типе указывает на важнейшее свойство
редуцированных спектров наряду с их информативными свойствами (см.
лемму~5.1 ниже).

\noindent Однако всё ещё остаётся следующее существенное
препятствие: такие спектры на их различных носителях не всегда
могут быть сравнимы друг с другом относительно их базовых свойств,
потому что их области определения -- ординальные спектры -- могут
содержать произвольно большие кардиналы \ $<k$ когда их носители
возрастают до \ $k$. Чтобы преодолеть это препятствие  мы
преобразуем их в редуцированные матрицы в следующем параграфе.

\newpage
\quad 

\newpage

\section{Редуцированные матрицы}
\setcounter{equation}{0}

Мы начинаем формировать матричные функции. Для этого, обладая
редуцированными спектрами, мы переходим к их простым
трансформациям -- редуцированным матрицам, служащим значениями
таких функций. Такие матрицы получаются из редуцированных спектров
простым изоморфным перечислением их областей определения:

\begin{definition}
\label{5.1.}
\ \\
{\em 1)}\quad Мы называем матрицей редуцированной к ординалу
 \ $\chi $ \  любое
отношение \ $M$, \ выполняющее следующую формулу, обозначаемую
через \ $\mu  ( M,\chi  )$:
\[
    (M \mbox{\it \ это функция})\wedge (dom ( M )
    \mbox{\it \ это ординал})\wedge rng ( M ) \subseteq B_{\chi }~.
\]

\noindent {\em 2)}\quad Пусть \  $M$ \  это матрица и \
$M_{1}\subset k\times B$. Мы называем наложением матрицы\ $M$ \ на
\ $M_{1}$ \ функцию \ $f$, \ выполняющую следующую формулу,
обозначаемую через \ $f:M\Rightarrow M_{1}$~:
\begin{eqnarray*}
&(f \ \mbox{\it это порядковый изоморфизм}\ dom (M )\mbox{\it \ на
}\ \ dom ( M_{1}))\wedge \qquad\qquad
\\
& \wedge \forall \alpha ,\alpha ^{\prime }\forall \Delta ,\Delta
^{\prime } ( f ( \alpha  ) =\alpha ^{\prime }\wedge  ( \alpha
,\Delta  ) \in M\wedge  ( \alpha ^{\prime },\Delta ^{\prime } )
\in M_{1}\longrightarrow \qquad
\\
& \qquad\qquad\qquad\qquad\qquad\qquad\qquad\qquad\qquad
\longrightarrow \Delta =\Delta ^{\prime } ).
\end{eqnarray*} Если наложение существует, то говорим, что \ $M$
 \  накладывается на \  $M_{1}$ \  и пишем \
$M\Rightarrow M_{1} $.

\noindent {\em 3)}\quad Если матрица \  $M$ \ накладывается на
спектр  \ $\widetilde{\mathbf{S}}_{n}^{\sin \vartriangleleft
\alpha }\overline{\overline{\lceil }}\chi $, то \ $M$ \ называется
матрицей этого спектра на \ $\alpha $, \ или субнедостижимо
универсальной матрицей уровня \  $n$ \ редуцированной к \ $\chi$ \
на \ $\alpha $.
\\
{\em 4)}\quad В этом случае если \ $(\alpha^{\prime} ,\Delta ) \in
\widetilde{\mathbf{S}}_{n}^{\sin \vartriangleleft \alpha
}\overline{\overline{\lceil }}\chi $, \ то  \ $\alpha^{\prime}$ \
называется кардиналом скачка матрицы \ $M$, \ а \ $\Delta$ \
называется их булевым значением на \ $\alpha$.
\\
{\em 5)}\quad В этом случае кардинал \ $\alpha$ \ называется
носителем матрицы \ $M$.  \hspace*{\fill} $\dashv$
\end{definition}

\noindent По этому определению если \ $M\Rightarrow
\widetilde{\mathbf{S}}_{n}^{\sin \vartriangleleft \alpha
}\overline{ \overline{\lceil }}\chi$, \ то
\[
    rng (M) = rng (\widetilde{\mathbf{S}}_{n}^{\sin
    \vartriangleleft \alpha }\overline{ \overline{\lceil }}\chi),
\]
поэтому начиная с этого места мы будем рассматривать булев спектр
на  \ $\alpha$ \ также как \ $rng (M)$.
\\
\hfill {} \\
Введённые выше матрицы могут быть названы одномерными; следуя
этому определению можно ввести многомерные матрицы, накладываемые
на многомерные редуцированные спектры той же размерности ( см.
комментарий 4) для более тонкого анализа утверждений нашего языка.
Многомерные матрицы этого типа использовались автором долгое время
как основной инструмент исследования субнедостижимости.
\\
В дальнейшем носители \  $\alpha $ \ и редуцирующие каординалы \
$\chi $ \ всегда будут $SIN_{n-2}$-кардиналами, \ $\chi <\alpha
\leq k$ \ (если только не рассматривается некоторый другой
случай).
\\
\hfill {} \\
Из этого определения и лемм \ref{4.4.}, \ref{4.5.} вытекают
следующие две леммы:

\begin{lemma}
\label{5.2.} \hfill {} \\
\hspace*{1em} Формулы
\[
    f:M\Rightarrow \widetilde{\mathbf{S}}_{n}^{\sin
    \vartriangleleft \alpha }\overline{\overline{\lceil }}\chi
    \quad , \ M\Rightarrow \widetilde{\mathbf{S}}_{n}^{\sin
    \vartriangleleft \alpha }\overline{\overline{\lceil }}\chi
\]
принадлежат соответственно классам \ $\Pi_{n},\Sigma _{n+1}$ \ для
\ $\alpha =k$ \ и классу  \ $\Delta _{1}$ \ для \ $\alpha <k$.
\end{lemma}

\vspace{-6pt}
\quad \\
\begin{lemma}
\label{5.3.} \hfill {} \\
\hspace*{1em} Пусть
\[
    M\Rightarrow \widetilde{\mathbf{S}}_{n}^{\sin
    \vartriangleleft \alpha } \overline{\overline{\lceil }}\chi ,
\]
тогда \medskip

\noindent {\em 1)} \ $\left\| \widetilde{u}_{n}^{\sin
\vartriangleleft \alpha } ( \underline{l} ) \right\| \lceil \chi
=\sum rng ( M )
$; \\

\noindent {\em 2)} \ $OT ( M ) =dom ( M ) \leq Od ( M ) <\chi
^{+}.$   \hspace*{\fill} $\dashv$
\end{lemma}

\noindent Последнее утверждение показывает, что теперь
редуцированные матрицы могут быть сравнены друг с другом (в смысле
гёделевой функции $Od$) \textit{только внутри \ $L_{\chi^{+}}$ \ }
и это обстоятельство далее сделает возможным определение матричных
функций с нужными свойствами.
\\
Основную роль далее играют матрицы и спектры, редуцированные к так
называемым полным кардиналам; их существование следует из леммы
\ref{2.5.}~1) (для \ $\widetilde{u}_n^{\sin}$, \ $k$ \ как \
$\varphi$, \ $\alpha_1$ \ соответственно):

\begin{definition}
\label{5.4.} \hfill {} \\
\hspace*{1em} Мы называем полным ординалом уровня  \ $n$
 \  всякий ординал \  $\chi $ \  такой, что
\[
    \exists x~\widetilde{u}_{n-1}^{\sin } ( x,\underline{l} )
    \dashv \vdash \exists x<\chi
    ~\widetilde{u}_{n-1}^{\sin}( x,\underline{l} ) .
\]
Наименьший из таких ординалов обозначается через \ $\chi ^{\ast
}$, \ в то время как значение \ $\left\| \widetilde{u}_{n}^{\sin }
(\underline{l} ) \right\| $ \ обозначается через \ $A^{\ast}$~.
\end{definition}

\vspace{-6pt}

\begin{lemma}
\label{5.5.} \

\quad

\noindent {\em 1)}\ $\chi ^{\ast }=\sup dom \left (
\widetilde{\mathbf{S}}_{n}^{\sin } \right ) =\sup dom \left (
\mathbf{S}_{n}\right ) <k$ \quad;

\quad

\noindent  {\em 2)}\ $SIN_{n-1} ( \chi ^{\ast } ) ~ , ~ \chi
^{\ast }=\omega _{\chi ^{\ast }}$\quad;

\quad

\noindent {\em 3)}\ $OT ( \chi ^{\ast }\cap SIN_{n-1} ) =OT (
dom  \left ( \widetilde{ \mathbf{S}}_{n}^{\sin } \right ) ) =OT (
dom \left ( \mathbf{S}_{n} \right )) =\chi ^{\ast }$;

\quad

\noindent {\em 4)}\ $\widetilde{\mathbf{S}}_{n}^{\sin
}=\widetilde{\mathbf{S}}_{n}^{\sin \triangleleft \chi^{\ast}} =
\widetilde{\mathbf{S}}_{n}^{\sin }\overline{\overline{\lceil
}}\chi ^{\ast }; \ \ \mathbf{S}_{n}= \mathbf{S}_{n}^{\triangleleft
\chi^{\ast}} =\mathbf{S}_{n} \overline{\overline{\lceil }}\chi
^{\ast }; \ \ A^{\ast }=A ^{\ast } \lceil \chi ^{\ast }$,

\quad \\
и аналогично для редукции этих спектров к любому \
$SIN_{n-1}$-кардиналу \ \ $\geq \chi ^{\ast }$.
\end{lemma}

\noindent \textit{Доказательство} 1) следует из лемм 3.10~2), 2.7,
1.1; утверждение 2) следует из лемм~\ref{3.10.}~2), \ref{3.4.}~1),
\ref{3.5.} (для уровня \ $n-1$) когда \ $\alpha _{1}=k$. После
этого 3),~4) устанавливаются посредством расщепляющего метода
доказательства леммы ~\ref{4.6.}.~\hfill~$\dashv$

Можно видеть, что лемма 5.5 2) наилучшая возможная и
\ $\chi^{\ast}
\notin SIN_{n} $.
\\
\hfill {} \\
Из этой леммы следует ещё одно важное свойство полного кардинала,
усиливающее лемму 3.2 (для \ $\alpha_{1}=k , \alpha=\chi^{\ast}
$), потому что здесь \ $\chi^{\ast} $ \  это только  \  $
SIN_{n-1} $-кардинал:

\begin{lemma}
\label{5.6.} \hfill {} \\
\hspace*{1em} Пусть \  $\forall x~\varphi $ \  это \ $\Pi
_{n}^{\dashv \vdash }$-утверждение, \  $\varphi \in \Sigma
_{n-1}^{\dashv \vdash }$, \ тогда:
\\
\\
{\em1)}\quad если \  $\varphi $ \  содержит индивидные константы
только из \ $L_{\chi ^{\ast }}^{B_{\chi ^{\ast }}}$, то для каждой
\ $\mathfrak{M}$-генерической функции \ $l$
\[
    L_{k} [ l] \vDash \left( \forall x\vartriangleleft ^{\underline{l}}
    \chi ^{\ast} ~\varphi \longleftrightarrow \forall
    x~\varphi \right)~;
\]
{\em 2)}\quad если \  $\varphi $ \  содержит индивидные константы
$\vartriangleleft \chi^{\ast} $ \ только из \ $L_{k} $ \ и не
содержит \ $\underline{l}$, \ то
\[
    \quad L_{k}\vDash \left( \forall x\vartriangleleft
    \chi ^{\ast }~\varphi \longleftrightarrow \forall
    x~\varphi \right)~;
\]
{\em 3)}\quad пусть \  $\omega _{0}^{\ast }=\sup dom \left (
\widetilde{\mathbf{S}} _{n}^{\sin }\overline{\overline{\lceil }} (
\omega _{0}+1 ) \right ) $  \  и \  $\varphi $ \  не содержит
индивидных констант и \ $\underline{l}$, \ тогда
\[
    \qquad L_{k}\vDash \left( \forall x\vartriangleleft
    \omega _{0}^{\ast}~\varphi \longleftrightarrow
    \forall x~\varphi \right)~.
\]
\end{lemma}
\textit{Доказательства} 1) - 3) аналогичны и сводятся к тому, что
все ординалы скачка утверждения \ $\exists x~\neg \varphi $ \
(если они существуют) меньше чем \ $\chi ^{\ast }$, \ а в 3) даже
меньше чем \ $\omega _{0}^{\ast }$; \ поэтому продемонстрируем их
сначала для~2).
\\
Предварительно нужно сделать следующее простое замечание: для
всяких \ $x \in L_{k}$, \ $m \geq 2 $, \ $\alpha\in SIN_{m}$ \ и \
$\mathfrak{M} $-генерической \ $l $\ выполняется \quad
$x\vartriangleleft \alpha \longleftrightarrow x
\vartriangleleft^{l} \alpha $.
\\
Поэтому в таком случае ограничение  \ $x \vartriangleleft \alpha $
\ следует рассматривать как \ $Od(x)<\alpha $ \  в  \ $L_{k} $ \ и
как  \ $Od^{l}(x)<\alpha $ \ в  \  $L_{k}[l]$.
\\
Рассмотрим утверждение \ $\forall x \ \ \varphi (x,\alpha_{0}) $,
\ $\varphi \in \Sigma_{n-1} $, \ только с одной индивидной
ординальной константой \ $\alpha_{0} \vartriangleleft \chi^{\ast }
$ \ (для большей ясности) и пусть \ $\alpha_{0} \in dom \left
(\widetilde{\mathbf{S}}_{n}^{\sin} \right ) $ \ (иначе можно
использовать лемму 5.5.~3)~). Сначала костанта \ $\alpha_{0} $ \
должна быть зафиксирована как в доказательстве леммы 4.6 некоторой
\ $\mathfrak{M} $-генерической \ $l \stackrel{\ast }{\in }
\widetilde{\bigtriangleup}_{n}^{\sin} (\alpha_{0}) $. Предположим,
что
\[
    L_{k} \vDash \exists x \neg \varphi (x, \alpha_{0}),
\]
тогда
\[
    L_{k} [l] \vDash \exists \alpha \ \  \varphi_{1} (\alpha,l)
\]
\\
где  \  $\varphi_{1} (\alpha, \underline{l}) $ \  это следующая  \
$\Pi_{n-1}$-формула:
\[
    SIN_{n-1}(\alpha) \wedge \exists \alpha^{\prime} <
    \alpha  \left ( \underline{l} \stackrel{\ast }{\in }
    \widetilde{\bigtriangleup}_{n}^{\sin \triangleleft \alpha}
    (\alpha^{\prime}) \wedge \exists x \vartriangleleft \alpha
    \neg \varphi^{L} (x,\alpha^{\prime}) \right ) .
\]
Утвеждение   \  $\varphi_{2}(\underline{l}) = \exists \alpha \ \
\varphi_{1} (\alpha,\underline{l})  $ \ уже не имеет индивидных
констант и выполняется \ $l \stackrel{\ast }{\in }
\bigtriangleup_{\varphi_{2}}(\alpha) $ \ для некоторого его
ординала скачка \ $\alpha $. \ По леммам~2.7, \ 5.5~1) \ $ \alpha
< \chi^{\ast} $ \ и для некоторого ординала \ $\alpha^{\prime} <
\alpha $
\[
    L_{k} [l] \vDash  \left ( l \stackrel{\ast }{\in }
    \widetilde{\bigtriangleup}_{n}^{\sin \triangleleft \alpha}
    (\alpha^{\prime}) \wedge \exists x \vartriangleleft \alpha \neg
    \varphi^{L} (x,\alpha^{\prime}) \right ).
\]
Так как  \  $\alpha \in SIN_{n-1} $ \  и  \ $\alpha^{\prime} <
\alpha $, \ то здесь можно опустить ограничение \ $
\vartriangleleft^{l} \alpha $ \  и поэтому  \ $ l \stackrel{\ast
}{\in } \widetilde{\bigtriangleup}_{n}^{\sin} (\alpha^{\prime}) $.
\ Это влечёт   \  $\alpha^{\prime} = \alpha_{0} $ \  и
\[
    L_{k} [l] \vDash   \exists x \vartriangleleft \alpha \neg
    \varphi^{L} (x,\alpha_{0})
\]
и, наконец,
\[
    L_{k} \vDash   \exists x \triangleleft \chi^{\ast}
    \neg \varphi (x,\alpha_{0}) .
\]
\\
Обращаясь к 3) предположим, что \ $L_{k} \vDash \exists x \neg
\varphi (x) $, \ тогда для каждой
 \ \mbox{$\mathfrak{M}$-генерической} \ $l $
\[
    L_{k} [l] \vDash \exists \alpha \ \ \varphi_{1}^{L} (\alpha)
\]
где
\[
    \ \varphi_{1} (\alpha) = \exists x \vartriangleleft
    \alpha \neg \varphi(x) .
\]
\\
Из леммы~1.4 ясно, что
\[
    \left \| \exists \alpha \ \ \varphi_{1}^{L} (\alpha)  \right \| = 1
\]
и для каждого \ $\alpha$
\[
    \quad \left \|
    \varphi_{1}^{L} (\alpha) \right \| \in  \left \{ 0;1 \right \} .
\]
Пусть  \  $\mathfrak{n}$ \  это гёделев номер формулы \ $\exists
\alpha \ \ \varphi_{1}^{L} (\alpha) $ \ и  \ \mbox{$p_{
\mathfrak{n}}= \left \{ (\omega_{0}, \mathfrak{n}) \right \} $}. \
Отсюда и из определений~2.6, 3.9 можно видеть, что

\[
    p_{ \mathfrak{n}} \leq  \left \| \exists \alpha \ \
    \widetilde{u}_{n-1} (\alpha,\underline{l}) \right \|
\]
и для каждого   \  $\alpha $
\[
    p_{\mathfrak{n}} \cdot \left \| \widetilde{u}_{n-1}
    (\alpha,\underline{l}) \right \| = 0
    \mbox{\it \quad или \quad} p_{\mathfrak{n}}  \leq \left \| \widetilde{u}_{n-1}
    (\alpha,\underline{l}) \right \|~.
\]

Эти утверждения сохраняются при редуцировании к \ $ \chi =
\omega_{0}+1 $ \ значений  \ $\left \| \widetilde{u}_{n-1}
(\alpha,\underline{l}) \right \| $. \ По определению~4.1 это
влечёт существование некоторого кардинала скачка \
$\alpha_{\mathfrak{n}} $ \ такого, что
\[
    p_{\mathfrak{n}}\cdot \widetilde{\bigtriangleup}_{n}^{\sin}
    (\alpha_{\mathfrak{n}}) \overline{\lceil}(\omega_{0}+1)>0 .
\]
Следовательно
\[
    \alpha_{\mathfrak{n}} \in dom  \left ( \widetilde{\mathbf{S}}_{n}^{\sin}
    \overline{\overline{\lceil}} (\omega_{0}+1)\right )
    \mbox{\it \quad и \quad} \alpha_{\mathfrak{n}}<\omega_{0}^{\ast} .
\]
После этого по определениям~3.9, 2.6 \
\[
    L_{k} [l] \vDash \exists \alpha^{\prime} < \alpha_{\mathfrak{n}} \ \
    \varphi_{1}^{L} (\alpha^{\prime})
\]
и, наконец,
\[
    L_{k} \vDash \exists x \vartriangleleft
    \omega_{0}^{\ast} \neg \varphi (x).
\]
\hfill $\dashv$

Используя здесь аргументы из доказательства~2) можно допустить
индивидные константы \ $ \in dom \left (
\mathbf{\widetilde{S}}_{n}^{\sin} \overline{\overline{\lceil}}
(\omega_{0}+1) \right ) $ \ в
 \ $\forall x \varphi $,  \ но в дальнейшем это не потребуется.
\\
\hfill {} \\
\hfill {} \\
В дальнейшем очень особую роль играют так называемые
\textit{сингулярные матрицы}:
\begin{definition}
\label{5.7.} \hfill {} \\
\hspace*{1em} Через \  $\sigma  ( \chi ,\alpha  ) $ \ обозначается
конъюнкция следующих формул:
\\

{\em 1)}\quad $SIN_{n-2} ( \alpha  ) \wedge  (\chi
\mbox{\it \ предельный кардинал } <\alpha )$~; \\

{\em 2)}\quad $OT ( SIN_{n-1}^{<\alpha } ) =\alpha $~;\\

{\em 3)}\quad $\sup dom \left ( \widetilde{\mathbf{S}}_{n}^{\sin
\vartriangleleft \alpha }\overline{\overline{\lceil}}\chi\right )
=\alpha $~.
\\
\\
И пусть \  $\sigma  ( \chi ,\alpha ,M ) $ \ обозначает формулу:
\[
    \sigma  ( \chi ,\alpha  ) \wedge  ( M\Rightarrow
    \widetilde{\mathbf{S}}_{n}^{\sin \vartriangleleft \alpha }
    \overline{\overline{\lceil }}\chi  ) .
\]
Матрица \  $M$  \  и спектр\ $\widetilde{\mathbf{S}}_{n}^{\sin
\vartriangleleft \alpha} \overline{\overline{\lceil }}\chi $ \
редуцированные к \ $\chi $ \ называются сингулярными на носителе \
$\alpha $ \ (на интервале \ $[\alpha _{1},\alpha _{2}[\ $) \
\textit{если} выполняется \ $ \sigma ( \chi ,\alpha ,M ) $ \
 (для некоторого \ $\alpha \in \lbrack \alpha
_{1},\alpha _{2}[\ $).
\\
Символ \ $S$  \  далее используется для общего обозначения
сингулярных матриц.  \hspace*{\fill} $\dashv$
\end{definition}

В дальнейшем все матрицы и спектры \textit{будут редуцированы к
определённому кардиналу \ $\chi$ \ и будут сингулярными на их
рассматриваемых носителях}; все рассуждения будут проводиться в \
$L_{k}$ \ ( или в \ $\mathfrak{M}$, \ если контекст не указывает
на другой случай).
\\
Из определения~5.7 и леммы~5.2 непосредственно следует

\begin{lemma}
\label{5.8.} \hfill {} \\
{\em 1)}\quad Формулы \ $\sigma  ( \chi ,\alpha  ) $, \ $\sigma  (
\chi ,\alpha ,M ) $ \ принадлежат классу \ $\Pi _{n-2}$.
\hfill {} \\
{\em 2)}\quad Если \quad $\sigma  ( \chi ,\alpha  ) $,\quad то
\[
    \widetilde{A}_{n}^{\sin \vartriangleleft \alpha }
    (\chi  ) \lceil \chi < \sum rng \left
    (\widetilde{\mathbf{S}}_{n}^{\sin \vartriangleleft \alpha }
    \overline{\overline{\lceil}}\chi  \right ) =\left\|
    \widetilde{u}_{n}^{\sin \vartriangleleft \alpha }
    ( \underline{l} ) \right\| \lceil\chi .
\]
\end{lemma}

Благодаря лемме \ref{5.8.}~2) можно ввести следующие важные
кардиналы: \vspace{-6pt}

\begin{definition}
\label{5.9.} \hfill {} \\
\hspace*{1em} Пусть выполняется \ $\sigma  ( \chi ,\alpha, S ) $
,\ тогда мы называем  кардиналом скачка и кардиналом предскачка
после\ $\chi $ \ матрицы \ $S$ \ и спектра на носителе \ $\alpha
$, \ или, короче, самого кардинала \ $\alpha$ \ соответственно
следующие кардиналы:
\[
    \alpha _{\chi }^{\downarrow }=\min \{ \alpha ^{\prime}
    \in ]\chi ,\alpha \lbrack ~:~ \widetilde{A}_{n}^{\sin
    \vartriangleleft \alpha } ( \chi  ) \lceil \chi
    <\widetilde{A}_{n}^{\sin \vartriangleleft \alpha } (
    \alpha^{\prime } ) \lceil \chi \wedge
\]
\[
    \qquad \qquad \qquad \qquad \qquad \qquad \qquad \qquad \qquad \qquad
    \wedge SIN_{n-1}^{<\alpha } ( \alpha^{\prime } ) \} ;
\]
\[
    \alpha _{\chi }^{\Downarrow }=\sup \{ \alpha ^{\prime}
    <\alpha _{\chi }^{\downarrow }  ~:~
    \widetilde{A}_{n}^{\sin \vartriangleleft \alpha } ( \chi )
    \lceil \chi =\widetilde{A}_{n}^{\sin \vartriangleleft \alpha }
    (\alpha ^{\prime }) \lceil \chi \wedge
\]
\[
    \qquad \qquad \qquad \qquad \qquad \qquad \qquad \qquad \qquad \qquad
    \wedge SIN_{n-1}^{<\alpha } ( \alpha ^{\prime } ) \}.
\]
\end{definition}

В этих обозначениях и везде далее индекс \ $\chi $ \ может быть
опущен, если он восстанавливается из контекста (или произволен).

\begin{lemma}
\label{5.10.} \hfill {} \\
\hspace*{1em} Пусть выполняется \  $\sigma  ( \chi ,\alpha, S ) $,
\ тогда кардиналы \  $\alpha _{\chi }^{\downarrow }$, \ $\alpha
_{\chi }^{\Downarrow }$ \ существуют и

\quad

{\em 1)}~$~\alpha _{\chi }^{\downarrow }=\min \left\{ \alpha
^{\prime }>\chi :\alpha ^{\prime }\in dom   \left (
\widetilde{\mathbf{S}}_{n}^{\sin \vartriangleleft \alpha }
\overline{\overline{\lceil }}\chi \right ) \right\} \quad;$

\quad

{\em 2)}\quad $\alpha _{\chi }^{\Downarrow }<\alpha _{\chi
}^{\downarrow }<\alpha ;\quad ]\alpha _{\chi }^{\Downarrow
},\alpha _{\chi }^{\downarrow }[\cap SIN_{n-1}^{<\alpha
}=\varnothing \quad;$

\quad

{\em 3)}\quad $\alpha _{\chi }^{\downarrow },$ $\alpha _{\chi
}^{\Downarrow }\in SIN_{n-2} \quad$.
\end{lemma}
\textit{Доказательство} 1), 2) следует из определений и леммы~3.4,
а доказательство 3) -- из лемм~3.6, \ref{3.8.}~2) (для \
$k,\alpha,n-2 $ \ как \ $\alpha_{1},\alpha_{2},n $ \
соответственно). \hfill $\dashv$
\\

Следующая лемма вместе с леммой 4.6 о спектральном типе составляет
основные инструменты в обращении с матричными функциями и
представляет важное информативное свойство редуцированных матриц и
соответствующих спектров.\\
До сих пор, рассматривая редуцированные спектры, мы использовали
главным образом их первые проекции, то есть ординальные спектры.
Что касается их вторых проекций (то есть редуцированных булевых
спектров, или редуцированных матриц), то они играли
вспомогательную роль, служа инструментами в доказательствах
лемм~\ref{4.6.}, \ref{5.5.}, \ref{5.6.} и других. Однако эти
вторые прoекции (или редуцированные булевы спектры, или
редуцированные матрицы) обладают следующим важным
характеристическим свойством: они содержат информацию о начальных
частях универсума, ограниченных кардиналами скачка на их
носителях, так что когда такой спектр или матрица переносится с
одного носителя на другой, свойства таких начальных частей
универсума сохраняются. Более точно:

\begin{lemma}
\label{5.11.} {\em (О матричной информативности)} \hfill {} \\
\hspace*{1em} Пусть  \ $S$ \  это редуцированная к \ $\chi$ \
матрица на носителях \  $\alpha _{1},\alpha _{2}>\chi $, \
накладываемая на спектры
\[
    S \Rightarrow \widetilde{\mathbf{S}} _{n}^{\sin \vartriangleleft
    \alpha_{1}}\overline{\overline{\lceil }} \chi, \quad
    S \Rightarrow \widetilde{\mathbf{S}} _{n}^{\sin \vartriangleleft
    \alpha_{2}}\overline{\overline{\lceil }} \chi
\]
на этих носителях, и пусть
\[
    \overrightarrow{a_{1}}= ( \alpha _{10},\alpha_{11},...,\alpha_{1m}),
    \quad \overrightarrow{a_{2}}= ( \alpha _{20},\alpha_{21},...,\alpha_{2m})
\]
это кортежи кардиналов скачка из ординальных спектров
\[
    dom  \left ( \widetilde{\mathbf{S}}_{n}^{\sin \vartriangleleft
    \alpha _{1}} \overline{\overline{\lceil }}\chi  \right ),
    \quad
    dom  \left ( \widetilde{\mathbf{S}}_{n}^{\sin \vartriangleleft
    \alpha _{2}} \overline{\overline{\lceil }}\chi  \right )
\]
или кардиналов предскачка на этих носителях, соответствующие
равным булевым значениям:
\[
    \widetilde{\mathbf{S}}_{n}^{\sin \vartriangleleft
    \alpha _{1}}\overline{\overline{\lceil }}\chi  (
    \alpha _{1i} ) =\widetilde{\mathbf{S}}_{n}^{\sin
    \vartriangleleft \alpha_{2}}\overline{\overline{\lceil }}
    \chi ( \alpha _{2i} ) , \ \ i=\overline{0,m}\mathit{.}
\]
Пусть, наконец, \  $\psi  ( x_{1},x_{2},...,x_{m},\underline{%
l} )$ \ это произвольная формула произвольного уровня со
свободными переменными \ $ x_{1}$, $~x_{2}$,$...$, $x_{m}$ \ и без
индивидных констант. Тогда

\[
    \psi ^{\vartriangleleft \alpha _{10}} ( \alpha _{11},...,\alpha _{1m},
    \underline{l} ) \dashv \vdash \psi ^{\vartriangleleft \alpha_{20}}
    ( \alpha _{21},...,\alpha _{2m},\underline{l} ) .
\]
\end{lemma}
\textit{Доказательство} проводится методами, аналогичными
использованным в доказательстве леммы~\ref{4.6.}. Рассмотрим для
большей ясности случай, когда формула  \ $\psi $ \ не содержит
свободные переменные и имя \ $\underline{l} $ \ и \
$\overrightarrow{a_{1}}$, $\overrightarrow{a_{2}}$ \ состоят
только из кардиналов  скачка (или предскачка) после \  $\chi $ \
на $\alpha_{1}$, $\alpha _{2}$ \ соответственно:
\[
    \alpha _{10}=\alpha _{1\chi }^{\downarrow }, \
    \alpha_{20}=\alpha _{2\chi }^{\downarrow },
\]
так как именно этот случай потребуется в дальнейшем. Пусть \
$\alpha _{11}$, $\alpha _{21}$ \  это следующие за \ $\alpha
_{10}$, $\alpha _{20}$ \ кардиналы скачка в редуцированных
ординальных спектрах соответственно на носителях
 \ $ \alpha _{1}$, $\alpha _{2}$. \
 \\
Идея доказательства состоит в следующем:
\\
Выполнение утверждения \ $\psi ^{\vartriangleleft \alpha _{10}}$ \
означает, что значение

\[
    \Delta =\widetilde{\mathbf{S}}_{n}^{\sin \vartriangleleft
    \alpha_{1}}\overline{\overline{\lceil }}\chi
    ( \alpha _{11} )
\]
содержит соответствующее условие\ $ p\in P_{\chi }$ \ в котором
оно ``закодировано''. Когда \ $S$ \ переносится с носителя \
$\alpha _{1}$ \ на носитель \ $\alpha _{2}$ \ значение
\[
    \Delta =\widetilde{\mathbf{S}}_{n}^{\sin \vartriangleleft
    \alpha _{2}}\overline{\overline{\lceil }}\chi (\alpha_{21} )
\]
всё ещё содержит \  $p$ \  и это означает, что \
$\psi^{\vartriangleleft \alpha _{20}}$ \  тоже выполнено. Но более
удобно использовать вместо \ $p $ \ некоторую \  $ \mathfrak{M}
$-генерическую функцию \ $l $, \ включающую \ $p $. \ Нижний
индекс  \ $\chi $ \ в обозначениях будет опускаться для краткости.
\\
Итак, пусть  \  $\Delta$ \  это булево значение спектра \
$\widetilde{\mathbf{S}}_{n}^{\sin \triangleleft \alpha_{1}}
\overline{\overline{\lceil}} \chi $, \ соответствующее \
$\alpha_{10} = \alpha_{1}^{\downarrow} $, \ то есть
\[
    \Delta = \widetilde{\Delta}_{n}^{\sin \triangleleft
    \alpha_{1}}(\alpha_{10})\overline{\lceil} \chi,
\]
и пусть некоторая \  $\mathfrak{M} $-генерическая функция \ $l $ \
фиксирует это  \ $\Delta$, \ то есть \ $l\stackrel{\ast }{\in }
\Delta $, \ поэтому она фиксирует \ $\alpha_{10} $ \ при
рассмотрении носителя \ $\alpha_{1} $ \ и фиксирует \ $\alpha_{20}
$ \ при рассмотрении носителя \ $\alpha_{2} $. \ По определению
4.1

\begin{equation}
\label{e5.1} l \stackrel{\ast }{\in }  \left \|
\widetilde{u}_{n-1}^{\triangleleft
\alpha_{1}}(\alpha_{10},\underline{l}) \right \| \lceil \chi ~;
\end{equation}

\begin{equation}
\label{e5.2} l \stackrel{\ast }{\notin }  \left \| \exists
\alpha^{\prime} < \alpha_{10}  \left ( SIN_{n-1}^{\triangleleft
\alpha_{1}} (\alpha^{\prime}) \wedge
\widetilde{u}_{n-1}^{\triangleleft \alpha_{1}} (\alpha^{\prime},
\underline{l} ) \right ) \right \| \lceil \chi ~.
\end{equation}
\vspace{0pt}

\noindent Так как \ $\alpha_{10} \in SIN_{n-1}^{< \alpha_{1}} $, \
то ограничение \ $\vartriangleleft \alpha_{1} $ \ в (5.2) может
быть заменено на ограничение \ $\vartriangleleft \alpha_{10} $. \
После опускания редукции \ $\chi $ \ в (5.1), (5.2) отсюда следует

\[
    l_{1} \stackrel{\ast }{\in } \left \|
    \widetilde{u}_{n-1}^{\triangleleft \alpha_{1}}(\alpha_{10},
    \underline{l}) \right \| ~;
\]
\[
    l_{1} \stackrel{\ast }{\notin } \left \|\exists \alpha^{\prime} <
    \alpha_{10}  \left ( SIN_{n-1}^{\triangleleft \alpha_{10}}
    (\alpha^{\prime}) \wedge \widetilde{u}_{n-1}^{\triangleleft
    \alpha_{10}}(\alpha^{\prime}, \underline{l}) \right ) \right \| ~,
\]
\vspace{0pt}

\noindent где  \  $l_{1} $ \  это некоторая \ $\mathfrak{M}
$-генерическая функция, совпадающая с \ $l $ \ на \ $\chi $; \ но
для некоторой краткости будем использовать прежний символ \ $l$. \
Тогда это означает, что

\begin{equation}
\label{e5.3} L_{k}[l] \vDash \widetilde{u}_{n-1}^{\triangleleft
\alpha_{1}}(\alpha_{10},l) \wedge \neg \exists \alpha^{\prime} <
\alpha_{10}  \left ( SIN_{n-1}^{\triangleleft
\alpha_{10}}(\alpha^{\prime}) \wedge
\widetilde{u}_{n-1}^{\triangleleft \alpha_{10}}(\alpha^{\prime},l)
\right )~.
\end{equation}
\vspace{0pt}

\noindent Таким образом кардинал  \ $\alpha_{10} $ \ определён в
 \  $ L_{k}[l] $ \  через  \  $l $, \ $\alpha_{1}
$ \  и также аналогично \ $\alpha_{20} $ \ определён через \ $l $,
\ $\alpha_{2} $ \  в  \  $ L_{k}[l] $. \ Теперь рассмотрим
утверждение  \ $\varphi(\underline{l}) = \exists \alpha \ \
\psi_{1}(\alpha, \underline{l}) $, \  где  \ $\psi_{1}(\alpha,
\underline{l}) $ \ это следующая  \ $\Pi_{n-1}$-формула:

\vspace{6pt}
\[
    SIN_{n-1}(\alpha) \wedge \widetilde{u}_{n-1} (\alpha,
    \underline{l}) \wedge \neg \exists \alpha^{\prime} < \alpha
    \left( SIN_{n-1}^{< \alpha} (\alpha^{\prime}) \wedge
    \widetilde{u}_{n-1}^{\triangleleft \alpha}
    (\alpha^{\prime},\underline{l}) \right ) \wedge
    \psi^{\triangleleft \alpha} ~.
\]
\vspace{0pt}

\noindent Из (5.3) ясно, что \vspace{6pt}
\begin{equation}
\label{e5.4} L_{k}[l]\vDash \psi^{\triangleleft \alpha_{10}}
\longleftrightarrow L_{k}[l] \vDash \exists \alpha < \alpha_{1} ~
\psi_{1}^{\triangleleft \alpha_{1}} (\alpha, l) ~.
\end{equation}
\vspace{0pt}
\\
После этого следует повторить аргумент из доказательства леммы
2.7: функция  \ $l $ \ должна быть заменена на функцию \ $l_{0} $
\ как в (2.1) и в то же время формула \ $\varphi $ \ должна быть
трансформирована в формулу \  $\varphi_{2} $ \  буквально как в
этом доказательстве, то есть посредством замены его подформул вида
\ $\underline{l}(t_{1})=t_{2} $ \ подформулами (2.2) и так далее.
Пусть \ $\mathfrak{n}$ \ это гёделев номер формулы \ $\varphi_{2}$
\ и  \ $l_{0}(\omega_{0}) = \mathfrak{n}$. \ В результате мы
переходим от\ $\varphi $ \ к спектрально универсальной формуле и
(5.4) влечёт \vspace{6pt}
\[
    L_{k}[l]\vDash \psi^{\triangleleft \alpha_{10}}
    \longleftrightarrow L_{k}[l_{0}] \vDash \exists x
    \vartriangleleft \alpha_{1} ~ u_{n-1}^{\triangleleft
    \alpha_{1}}(x, l_{0}) ~.
\]
\vspace{0pt}

\noindent Используя здесь  \ $\widetilde{u}_{n-1}^{\sin} $ \
вместо  \  $u_{n-1}$ \ из определения 5.7 можно видеть, что
\vspace{6pt}
\begin{equation}
\label{e5.5} L_{k}[l]\vDash \psi^{\triangleleft \alpha_{10}}
\longleftrightarrow L_{k}[l_{0}] \vDash \exists x \vartriangleleft
\alpha_{1} \widetilde{u}_{n-1}^{ \sin \triangleleft \alpha_{1}} ~
(x, l_{0}) \ .
\end{equation}
\vspace{0pt}

\noindent Анализируя конструкцию  формулы \
$\widetilde{u}_{n}^{\sin} $ \  из (5.5), (5.3) нетрудно заключить,
что \vspace{6pt}
\[
    L_{k}[l]\vDash \psi^{\triangleleft \alpha_{10}}
    \longleftrightarrow l_{0} \stackrel{\ast }{\in }
    \widetilde{\Delta}_{n}^{\sin \triangleleft
    \alpha_{1}} (\alpha_{11}) ~;
\]
\vspace{0pt}

\noindent напомним, что здесь  \ $\alpha_{11} $ \ это наследник  \
$\alpha_{10} $ \ в редуцированном ординальном спектре на носителе
\ $\alpha_{1} $. \ В этом рассуждении значения \ $l_{0}(\alpha) $
\ для \ $\alpha \geq \chi $ \  не использовались и поэтому
\vspace{6pt}
\begin{equation}
\label{e5.6} L_{k}[l]\vDash \psi^{\triangleleft \alpha_{10}}
\longleftrightarrow l_{0} \stackrel{\ast }{\in }
\widetilde{\Delta}_{n}^{\sin \triangleleft \alpha_{1}}
(\alpha_{11}) \overline{\lceil} \chi \quad.
\end{equation}
\vspace{0pt}

\noindent Матрица  \  $S $ \  накладывается на \
$\widetilde{\mathbf{S}}_{n}^{\sin \triangleleft \alpha_{2}}
\overline{\overline{\lceil}} \chi $ \  и поэтому \vspace{6pt}
\begin{equation}
\label{e5.7} L_{k}[l]\vDash \psi^{\triangleleft \alpha_{10}}
\longleftrightarrow l_{0} \stackrel{\ast }{\in }
\widetilde{\Delta}_{n}^{\sin \triangleleft \alpha_{2}}
(\alpha_{21}) \overline{\lceil} \chi ~.
\end{equation}
\vspace{0pt}

\noindent Это же рассуждение, повторённое  вплоть до (5.6), но
только с \ $\alpha_{1}, \alpha_{10},\alpha_{11} $ \ заменёнными
соответственно на \ $\alpha_{2}, \alpha_{20},\alpha_{21} $, \
устанавливает, что \vspace{6pt}
\[
    L_{k}[l]\vDash \psi^{\triangleleft \alpha_{20}}
    \longleftrightarrow l_{0} \stackrel{\ast }{\in }
    \widetilde{\Delta}_{n}^{\sin \triangleleft
    \alpha_{2}} (\alpha_{21}) \overline{\lceil} \chi ~.
\]
\vspace{0pt}

\noindent Остаётся сравнить эту эквивалентность с (5.7).
\\
В случае, когда  \  $\psi $ \  содержит свободные переменные \
$x_{1}...,x_{m} $ \  и  \ $\underline{l} $, \ следует следует
провести аналогичное рассуждение, рассматривая несколько \
$\mathfrak{M} $-генерических функций
\[
    l_{i} \stackrel{\ast }{\in }
    \widetilde{\Delta}_{n}^{\sin \triangleleft
    \alpha_{1}} (\alpha_{1i}) \overline{\lceil} \chi ,
    \quad i=\overline{0,m} ,
\]
и используя их прямое произведение как в доказательстве леммы 4.6.
\hfill $\dashv$
                            \\            \hfill {} \\
\\
Незначительно изменяя это рассуждение несложно провести его для
кардиналов предскачка  \ $\alpha_{10}=\alpha_{1\chi}^{\Downarrow},
\ \alpha_{20}=\alpha_{2\chi}^{\Downarrow}$~.
\\

\noindent Основные инструменты доказательства основной теоремы --
это матричные функции, представляющие собой последовательности
 сингулярных редуцированных матриц специального вида. Следующаяя
лемма доставляет возможность построения таких функций:

\begin{lemma}
\label{5.12.}
\[
    \forall \chi \forall \alpha _{0} \bigl(  ( \chi
    \mbox{\it \ предельный кардинал }>\omega _{0} ) \rightarrow
    \exists \alpha _{1}>\alpha _{0}~\sigma  ( \chi ,\alpha _{1} )  \bigr) .
\]
\end{lemma}

\noindent \textit{Доказательство.} \ Это \  $\Pi _{n}$-утверждение
не содержит индивидных констант или \ $\underline{l}$ \ и поэтому
по лемме \ref{5.6.} 3) достаточно доказать, что оно выполняется,
когда переменные \ $\chi $, $ \alpha _{0}$ \ ограничены кардиналом
\ $\omega _{0}^{\ast }$.
\\
Итак, пусть \  $\chi $, $\alpha _{0}<\omega _{0}^{\ast }$; \
рассмотрим кардинал
\[
    \alpha _{1}=\sup dom  \left (
    \widetilde{\mathbf{S}}_{n}^{\sin }\overline{
    \overline{\lceil }}\chi \right ).
\]
Можно видеть, что \ $\alpha _{1}\geq \omega _{0}^{\ast }$, \ и
поэтому \  $\chi $, $\alpha _{0}<\alpha _{1}$. Из лемм \ref{4.2.},
\ref{3.4.} (для \ $\alpha _{1}=k$) \  следует \ $\alpha _{1}\in
SIN_{n-1}$, \ поэтому условия 1), 3) определения \ref{5.7.}
выполняются; условие 2) может быть установлено методом,
составляющим упрощенную версию расщепляющего метода доказательства
леммы~\ref{4.6.}. И действительно, допустим, что напротив,
\vspace{6pt}
\[
    \alpha_{2} = OT  \left ( SIN_{n-1}^{< \alpha_{1}} \right )
    < \alpha_{1} ~,
\]

\noindent тогда существуют кардинал
\[
\alpha_{3} \in [ \alpha_{2}, \alpha_{1} [ \cap dom
    \left ( \widetilde{\mathbf{S}}_{n}^{\sin}
    \overline{\overline{\lceil}} \chi \right )
\]
\noindent и некоторая \  $\mathfrak{M} $-генерическая функция \ $l
\stackrel{\ast }{\in } \widetilde{\Delta}_{n}^{\sin} (\alpha_{3})
\overline{\lceil} \chi $. \ Теперь рассмотрим утверждение \
$\varphi(\underline{l}) = \exists \alpha \ \ \varphi_{1}(\alpha,
\underline{l}) $, \  где \ $\varphi_{1}(\alpha, \underline{l} ) $
\ это  \ $\Pi_{n-1} $\- формула:

\begin{eqnarray*}
&   SIN_{n-1}(\alpha)\wedge
\qquad\qquad\qquad\qquad\qquad\qquad\qquad\qquad\qquad\qquad
\\
& \wedge \exists
\alpha^{\prime},\alpha^{\prime\prime}<\alpha \ \ (
\alpha^{\prime}<\alpha^{\prime\prime} \wedge
SIN_{n-1}^{<\alpha}(\alpha^{\prime\prime})\wedge
\widetilde{u}_{n-1}^{\sin \triangleleft \alpha^{\prime\prime}}
(\alpha^{\prime},\underline{l})\wedge
\\
&   \wedge \neg \exists \alpha^{\prime\prime\prime} <
\alpha^{\prime} \ \ \widetilde{u}_{n-1}^{\sin \triangleleft
\alpha^{\prime\prime}}
(\alpha^{\prime\prime\prime},\underline{l}) \wedge
\\
& \qquad\qquad\qquad\qquad \wedge OT \left (
\left \{ \alpha^{\prime\prime\prime}< \alpha^{\prime\prime} :
SIN_{n-1}^{<\alpha^{\prime\prime}}(\alpha^{\prime\prime\prime})
\right \} \right ) = \alpha^{\prime} ) \ .
\end{eqnarray*}
\vspace{0pt}

\noindent Нетрудно видеть, что  \  $ L_{k} [l] \vDash \exists
\alpha \ \ \varphi_{1}(\alpha,l) $ \ и для каждого \ $\alpha $

\begin{eqnarray*}
L_{k}[l] \vDash \varphi_{1}(\alpha,l) \longleftrightarrow
SIN_{n-1}(\alpha)\wedge \exists
\alpha^{\prime\prime} < \alpha (
SIN_{n-1}(\alpha^{\prime\prime}) \wedge
\\
\qquad\qquad\qquad\qquad\qquad\qquad\qquad\qquad
\wedge OT \left (\alpha^{\prime\prime} \cap SIN_{n-1} \right ) =
\alpha_{3} ).
\end{eqnarray*}
\vspace{0pt}

\noindent Поэтому \ $ dom  \left
(\mathbf{S}_{\varphi}\overline{\overline{\lceil}} \chi  \right )$
\ содержит кардинал  \ $  > \alpha_{1}  $ \ и с помощью лемм~2.7,
3.10 можно видеть, что \ $dom \left (
\widetilde{\mathbf{S}}_{n}^{\sin}\overline{\overline{\lceil}} \chi
\right ) $ \ тоже содержит такие кардиналы в противоречии с
предположениями.\hfill $\dashv$
\vspace{12pt} \\

\noindent Построение матричных функций базируется на следующем
перечислении (в~$L_{k}$) \  субнедостижимых кардиналов:

\begin{definition}
\label{5.13.} \hfill {} \\
\hspace*{1em} Пусть \  $\alpha _{1}\leq k$.
\\
Мы определяем следующую функцию \ $\gamma _{f}^{<\alpha _{1}}= (
\gamma _{\tau }^{<\alpha _{1}} ) _{\tau }$~ ниже \  $\alpha _{1}$
рекурсией по \  $ \tau <\alpha _{1}$: \

\[
    \gamma _{0}^{<\alpha _{1}}=0 ~; \mbox{\it \quad для \quad} \tau>0
    \qquad \qquad\qquad\qquad\qquad\qquad\qquad
\]

\vspace{-6pt}

\[
    \gamma _{\tau }^{<\alpha _{1}}=\min \{ \gamma <\alpha_{1}:
    SIN_{n-1}^{<\alpha _{1}} ( \gamma  ) \wedge \forall
    \tau^{\prime }<\tau \quad
    \gamma_{\tau ^{\prime }}^{<\alpha_{1}}<\gamma \} \ .
\]
\vspace{0pt}

\noindent Определяется также обратная функция \ $\tau
_{f}^{<\alpha _{1}}= ( \tau_{\gamma }^{<\alpha _{1}} ) _{\gamma }$
:

\[
    \tau =\tau _{\gamma }^{<\alpha _{1}}\longleftrightarrow
    \gamma =\gamma _{\tau }^{<\alpha _{1}}.
\]
\end{definition}

Доказательство основной теоремы состоит в получении противоречия:
в создании в \ $L_{k}$ \ специальной матричной функции, обладающей
несовместными свойствами монотонности и немонотонности; эта
функция возникает в результате последовательного усложнения её
следующей простейшей версии:
\begin{definition}
\label{5.14.} \hfill {} \\
\hspace*{1em} Мы \medskip называем матричной функцией уровня \ $n$
\ ниже $ \alpha _{1}$, редуцированной к $\chi $, следующую функцию
\ \ \ $S_{\chi f}^{<\alpha _{1}}= ( S_{\chi \tau }^{<\alpha _{1}}
) _{\tau }$, \ принимающую значения: \vspace{6pt}
\[
    S_{\chi \tau }^{<\alpha _{1}}=\min_{\underline{\lessdot }}\{
    S:\exists \alpha < \alpha_{1}  \left ( \gamma_{\tau}^{<
    \alpha_{1}}< \alpha \wedge \sigma ^{\vartriangleleft
    \alpha _{1}} ( \chi ,\alpha ,S )\right ) \}  .
\]
\end{definition}
\vspace{0pt}

\noindent Таким образом эти значения являются матрицами \ $S$, \
редуцированными к \ $\chi$ \ и сингулярными на их носителях \
$\alpha$ ниже $ \alpha _{1}$.
\\

\noindent Как обычно, если \  $\alpha _{1}<k$, \ то все введённые
функции называются ограниченными или релятивизированными к \
$\alpha _{1}$; \ если  \ $\alpha_{1} = k $, \ то все упоминания о
\ $\alpha_{1} $ \ опускаются.
\\
Напомним, что все ограничивающие ординалы \ $\alpha_{1} $ \
полагаются \ $SIN_{n-2}$-кардиналами или \ $\alpha_{1}=k $.

\begin{lemma}
\label{5.15.} {\em (О абсолютности матричной функции).} \\
\hspace*{1em} Пусть \ $\chi < \gamma_{\tau+1}^{<\alpha_1} <
\alpha_2 < \alpha_1 \leq k$ \ и \ $\alpha_2 \in
SIN_{n-2}^{<\alpha_1}$, \ тогда:
\\

\noindent {\em 1)} функции \ \ $\gamma _{f}^{<\alpha_2}$, \
$\gamma _{f}^{<\alpha_1}$ \ \ тождественно совпадают на множестве
\ $\{ \tau^\prime: \gamma _{\tau^\prime} ^{<\alpha_2} \leq
\gamma_{\tau+1} ^{<\alpha_1} \}$:

\[
    \gamma _{\tau^\prime} ^{<\alpha_2} \leq
    \gamma_{\tau+1} ^{<\alpha_1}  \longrightarrow
    \gamma _{\tau^\prime} ^{<\alpha_2} =
    \gamma _{\tau^\prime} ^{<\alpha_1};
\]
\vspace{0pt}

\noindent {\em 2)} функции \ \ $S_{\chi f}^{<\alpha_2}$, \
$S_{\chi f}^{<\alpha _{1}} $ \ \ тождественно совпадают на
множестве \ $\{ \tau^\prime: \chi \leq \gamma _{\tau^\prime}
^{<\alpha_2} \leq \gamma_\tau ^{<\alpha_1} \}$:

\[
    \chi \leq
    \gamma _{\tau^\prime} ^{<\alpha_2} \leq
    \gamma_\tau ^{<\alpha_1}  \longrightarrow
    S _{\chi \tau^\prime} ^{<\alpha_2} =
    S _{\chi \tau^\prime} ^{<\alpha_1}.
\]
\end{lemma}
\textit{Доказательство} следует из лемы~\ref{3.8.} 1) (где \ $n$ \
заменяется на \ $n-1$) \ и леммы 5.17~2) $(ii)$ ниже. \hfill $\dashv$ \\

Из леммы~\ref{3.3.} (где \ $n$ \ также заменяется на \ $n-1$ \ )
следует
\begin{lemma}
\label{5.16.} \hfill {} \\
\hspace*{1em} Для\ $\alpha_{1}<k $ \ функции

\[
    \gamma =\gamma _{\tau }^{<\alpha _{1}}, \quad S=S_{\chi
    \tau}^{<\alpha _{1}}
\]
\vspace{0pt}

\noindent  \  $\Delta _{1}$-определимы через \ $\chi $, \ $\alpha
_{1}$.
\\
Для \  $\alpha_{1}=k $ \  эти функции  \ $\Pi_{n-1} $-определима,
\ $\Delta_{n}$-определима соответственно.  \hspace*{\fill}
$\dashv$
\end{lemma}

Следующая лемма представляет собой ``зародыш'' всех дальнейших
рассуждений: она устанавливает, что матричная функция обладает
свойством $\underline{\lessdot }$-монотонности -- и в дальнейшем
мы будем модифицировать эту функцию таким образом, чтобы и
сохранить и исключить это свойство одновременно.
\\
Поэтому мы часто будем возвращаться к идее этой леммы и её
доказательства в различных формах:

\begin{lemma}
\label{5.17.}
\quad \medskip \\
{\em 1)}\quad Функция \  $S_{\chi f}^{<\alpha _{1}}$ \ является
$\underline{\lessdot }$-монотонной, то есть для всяких \ $\tau
_{1}$, $ \tau _{2}\in dom \left( S_{\chi f}^{<\alpha _{1}} \right)
$
\[
    \quad \tau _{1}<\tau _{2}\longrightarrow S_{\chi
    \tau _{1}}^{<\alpha _{1}} \underline{\lessdot }
    S_{\chi \tau _{2}}^{<\alpha _{1}} \ .
\]
{\em 2)}\quad Пусть \ $ \tau \in dom \left ( S_{\chi f}^{<
\alpha_{1}} \right ) $, \ тогда:
\\
\begin{itemize}
\item[(i)] \ $ \{ \tau^{\prime}: \chi \leq
\gamma_{\tau^{\prime}}^{< \alpha_{1}} \leq \gamma_{\tau}^{<
\alpha_{1}} \} \subseteq dom \left ( S_{\chi f}^{< \alpha_{1}}
\right ) $~;
\\
\item[(ii)] \ если \ $\gamma_{\tau+1}^{<\alpha_{1}} $ \ и матрица \ $S_{\chi \tau}^{< \alpha_{1}}$ \ существуют, \ то
эта матрица обладает носителем \ $\alpha \in \;
]\gamma_{\tau}^{<\alpha_{1}}, \gamma_{\tau+1}^{<\alpha_{1}}[ $~.
\end{itemize}
\end{lemma}
\textit{Доказательство.} \ Утверждения 1, 2 $(i)$ следуют из
определения~5.14 непосредственно. Доказательство 2~$(ii) $
представляет типичное приложение леммы~3.2 об ограничении. Так как
существует матрица  \ $ S=S_{\chi \tau}^{< \alpha_{1}} $ \ на
некотором носителе
\[
    \alpha \in ]\gamma_{\tau}^{<\alpha_{1}}, \alpha_{1} [,
\]
то следующее \  $\Sigma_{n-1} $-утверждение \ $\varphi
(\chi,\gamma_{\tau}^{<\alpha_{1}}, S ) $~
\[
    \exists \alpha \ \ (\gamma_{\tau}^{< \alpha_{1}}< \alpha
    \wedge \sigma (\chi, \alpha, S))
\]
выполняется ниже\  $\alpha_{1} $, \ то есть выполняется \
$\varphi^{\triangleleft
\alpha_{1}}(\chi,\gamma_{\tau}^{<\alpha_{1}}, S) $.

Это утверждение содержит индивидные константы
\[
    \chi < \gamma_{\tau +1}^{<\alpha _{1}}, \quad
    \gamma _{\tau }^{<\alpha_{1}}<
    \gamma _{\tau +1}^{<\alpha _{1}},\quad S\vartriangleleft
    \gamma _{\tau +1}^{<\alpha _{1}}
\]
согласно лемме~\ref{5.3.}~2), и поэтому по лемме~\ref{3.2.} (где \
$n$ \ заменяется на \ $n-1$)  \ $ SIN_{n-1}^{<\alpha
_{1}}$-кардинал \ $\gamma _{\tau +1}^{<\alpha_{1}} $ \
ограничивает утверждение \ $\varphi $, \ то есть выполняется
утверждение
\[
    \exists \alpha \in \left[ \gamma _{\tau }^{<\alpha _{1}},
    \gamma_{\tau +1}^{<\alpha _{1}}\right[ ~\sigma ^{\vartriangleleft
    \alpha _{1}} ( \chi ,\alpha ,S ).
\]
\hfill $\dashv$

 Отсюда, из леммы 5.12 (для \ $\chi = \chi^{\ast}$) и
леммы~\ref{5.3.}~2) \ (для \ $\alpha_1 = k$) \ сразу же следует
\quad \\
\quad \\
\begin{lemma}
\label{5.18.} \quad \\
{\em 1)} Нерелятивизированная функция \ $S_{\chi^{\ast} f}$ \
определена и монотонна на заключительном сегменте недостижимого
кардинала $k$:

\[
    dom \left( S_{\chi ^{\ast } f } \right)=
    \left \{ \tau : \chi^{\ast} \leq \gamma_{\tau} < k \right \};
\]
\vspace{0pt}

\noindent {\em 2)}\quad эта функция стабилизируется, то есть
существует ординал \ $\tau ^{\ast }>\chi ^{\ast }$ \ такой, что
для каждого \ $\tau \geq \tau ^{\ast }$ существует
\[
    S_{\chi ^{\ast }\tau }=S_{\chi ^{\ast }\tau ^{\ast }};
\]
поэтому для каждого \ $\gamma_{\tau} \geq \chi ^{\ast }$ \
\[
    S_{\chi ^{\ast }\tau } \ \underline{\lessdot}
    \ S_{\chi ^{\ast }\tau ^{\ast }}.
\]
\end{lemma}

Теперь приведём эскиз первого приближения к идее доказательства
основной теоремы.
\\
Мы попытаемся получить требуемое противоречие на следующем пути.
Нижний индекс  \ $\chi ^{\ast } $ \ будет опускаться для некоторой
краткости.
\\
Рассмотрим матричную функцию \ $S_f$ на стадии её стабилизации, то
есть рассмотрим \ произвольный достаточно большой ординал \ $\tau
_{0}>\tau ^{\ast }$ \ и соответствующее значение этой функции, то
есть матрицу \ $ S_{\tau _{0}}$ \ на некотором её носителе \
$\alpha _{0}\in \left] \gamma_{\tau _{0}},\gamma _{\tau
_{0}+1}\right[ $, \ её кардинал предскачка \ $\alpha^0=\alpha
_{0}^{\Downarrow }$ \ и функцию \  $ S_{f}^{<\alpha^0}$, но уже
ниже этого кардинала.
\\
Предстоящие рассуждения будут более удобными и прозрачными, если
будет  использоваться следующее описание ситуации, когда наше
рассмотрение, ``опираясь'' или ``вставая''  \ на данный ординал \
$\alpha> \chi^{\ast}$, обозревает начальную часть универсума
\emph{ниже этого ординала.}
\\
\quad \\
Именно, мы будем говорить, что, \textit{ стоя на ординале \
$\alpha> \chi^{\ast}$, мы обнаруживаем (или видим) некоторую
информацию или некоторый объект ниже этого \ $\alpha$, \ если эта
информация или объект могут быть определены некоторой формулой,
ограниченной этим \ $\alpha$, -- причём определены как некоторый
ординал\ $< \chi^{\ast +}$.} \ Если этого сделать нельзя, то будем
говорить, что эту информацию или объект \textit{ нельзя обнаружить
(или увидеть) } \ ниже \ $\alpha$.
\\
\quad \\
Итак, идея использования матричной функции в её базовой форме
будет состоять в следующем:
\\
\hfill {} \\
\textit{Допустим, что матрица \ $S_{\tau_0}$ \  -- значение этой
матричной функции -- рассматривается на её носителе  \  $ \alpha
_{0}$. \ Стоя на кардинале предскачка \ $\alpha^0=\alpha
_{0}^{\Downarrow }$ \ после \ $\chi^{\ast}$ \ этой матрицы,
следует рассмотреть ситуацию ниже \ $\alpha^0$ и, главное,
поведение этой самой функции, но в её релятивизированной к \
$\alpha^0$ \ форме, с целью обнаружения её свойств}.
\\
\hfill {} \\
Тогда информация об этой функции, которую можно обнаружить и
выразить в виде некоторых ординалов \ $< \chi^{\ast +}$ -- эта
информация может доставить некоторые противоречия.
\\
Например, она может вызвать возрастание порядкового типа \
$S_{\tau_0}$ \ согласно лемме~4.6, и, следовательно, может
нарушить стабилизацию всей функции \ $S_f$ \ в противоречии с
леммой 5.18.
\\
Более подробно в этом контексте:
\\
По лемме~\ref{5.15.}~2) функция \ $S_f^{< \alpha^0}$ совпадает с \
$S_{f}$ \ на ординале \ $\tau _{0}$ \  и также \
$\underline{\lessdot }$-монотонна. Применим лемму~\ref{4.6.} к
этой ситуации, рассматривая
\[
    \overline{\delta }= \sup{}_{\tau} Od
    (S_{\tau }^{<\alpha^0} ),  \quad \overline{\chi }=
    \chi^{\ast }, \quad \overline{\alpha}_{0}=
    \alpha _{0}^{\downarrow }, \quad \overline{\alpha}_1 =
    \alpha_{0}.
\]
Предположим, что \ $\overline{\delta }<\chi ^{\ast +}$, \ тогда
можно определить \ $\overline{\delta}$ \ стоя на кардинале \
$\alpha^0$. \ Но теперь из определения~\ref{5.7.} и
леммы~\ref{5.5.} следует, что все условия леммы~\ref{4.6.}
выполнены и поэтому она влечёт противоречие:
\[
    Od ( S_{\tau ^{\ast }} ) \leq \overline{\delta }<OT
    ( S_{\tau_{0}} ) \leq Od ( S_{\tau ^{\ast }} )  .
\]
Следовательно, в действительности \ $\overline{\delta }=\chi
^{\ast +}$ \ и функция \ $S_{f}^{<\alpha^0}$ \  \
$\underline{\lessdot }$-неубывает до \ $\chi ^{\ast +}$.

Это означает, что функция \  $ ( S_{\tau } ) _{\tau <\tau _{0}}$,
\ оставаясь \ $\underline{\vartriangleleft }$~-огра\-ни\-чен\-ной
ординалом \ $Od \left ( S_{\tau^{\ast}} \right ) < \chi^{\ast +} $
\ (в смысле гёделевой функции \emph{Od}), теряет эту
ограниченность после её продолжения на множество \ $\{ \tau
:\gamma _{\tau}^{<\alpha^0 }<\alpha^0 \} $.
\\
Можно видеть, что это становится возможным только потому, что все
свойства субнедостижимости всех уровней \ $\geq n$ \ теряются
всеми кардиналами \ $\leq \gamma _{\tau _{0}}$ \  после их
релятивизации к \ $\alpha^0 $; \ в то же самое время возникают
некоторые \ $SIN_{n-1}$-кардиналы (релятивизированно к \
$\alpha^0$), не обладающие этим свойством ранее в
нерелятивизированном универсуме (Киселёв~\cite{Kiselev3},
\cite{Kiselev4}).
\\

Эти выводы означают, что многие важные свойства нижних уровней
универсума не продолжаются до некоторых релятивизирующих
кардиналов -- именно, до кардиналов скачка (а точнее предскачка)
редуцированных матриц на их носителях, являющихся значениями
матричных функций.

\textit{И для того, чтобы исключить это явление, мы введём
специальные кардиналы, так называемые диссеминаторы, которые
распространяют некоторые такие свойства без искажений до таких
кардиналов}.
\\
Требуемые распространения такого рода иногда производятся
субнедостижимыми кардиналами, но для дальнейшего требуется
существенно расширить и уточнить это понятие.

\newpage
\quad 

\newpage

\section{Диссеминаторы}
\setcounter{equation}{0}

Напомним, что все рассуждения проводятся и все формулы
интерпретируются в
 \ $ L_{k} $ \  (или в \ $\mathfrak{M}$ \ если,
конечно, контекст не имеет в виду другой случай). Также все \
$\vartriangleleft$-ограничивающие ординалы содержатся в \
$SIN_{n-2}$.
\\
Мы введём понятие диссеминатора в конструктивной структуре \ $
L_{k} $ \ для большей ясности рассуждений, хотя это может быть
сделано без всякой потери общности. Далее будут рассматриваться
классы формул \  $\Sigma _{m}$, \ $\Pi _{m}$ \ произвольного
фиксированного уровня \ $m>3$ (если только не оговаривается иное).
\begin{definition}
\label{6.1.} \hfill {} \\
\hspace*{1em} Пусть
\[
    0<\alpha <\alpha _{1} \leq k,\quad X\subseteq \alpha, \quad
    X\neq \varnothing .
\]
Мы будем называть ординал \  $\alpha $ \ диссеминатором уровня \
$m$ \  с базой данных \ $X$ \ ниже\ $\alpha _{1}$, \ если для
всяких \ $\mathfrak{n} \in \omega_0$ \ и кортежа \
$\overrightarrow{a}$ \ ординалов \ $\in X$ в  \ $ L_{k} $ \
выполняется
\[
    U_{m}^{ \Sigma \vartriangleleft \alpha _{1}} ( \mathfrak{n},
    \overrightarrow{a} ) \longrightarrow  U_{m}^{\Sigma
    \vartriangleleft \alpha } ( \mathfrak{n}, \overrightarrow{a} ) \quad.
\]
\vspace{0pt}

\noindent Формула, определяющая в \ $L_{k}$ \ класс всех
диссеминаторов с базой данных \ $X$ \ уровня \ $m $ \ ниже \
$\alpha _{1}$, \ будет обозначаться через \ $ SIN_{m}^{<\alpha
_{1}}\left[ X\right] ( \alpha ) $; \ само это множество также
будет обозначаться через \ $SIN_{m}^{<\alpha _{1}}\left[ X \right]
$, \ в то время как его диссеминаторы будут обозначаться общим
символом \ $\delta ^{X}$ \ или, коротко, символом \ $\delta $, \
указывая на \ $\alpha_1$ \ в контексте. \hspace*{\fill} $\dashv$
\end{definition}
Напомним, что символ \ $U_{m}^{\Sigma}$ \ обозначает здесь \
$\Sigma _{m}$-формулу, универсальную для класса \ $\Sigma
_{m}$-формул, но без всяких вхождений имени \ $\underline{l}$ \
(см. замечание после леммы 2.5).
\\
Можно получить определение понятия диссеминатора в более широком
смысле, если допустить вхождения имени \ $ \underline{l}$:
\\
\textit{для каждого кортежа \ $\overrightarrow{a}$ \ ординалов \
$\in X$}
\[
    \left\| u_{m}^{\Sigma \vartriangleleft \alpha _{1}}
    ( \overrightarrow{a},\underline{l} ) \right\| \leq
    \left\| u_{m}^{\Sigma \vartriangleleft \alpha }
    (\overrightarrow{a},\underline{l} ) \right\| ,
\]
и все дальнейшие рассуждения можно провести и в этом случае.
\\
\hfill {} \\
Следующие две леммы совершенно аналогичны леммам~3.3, 3.2 и могут
быть установлены тем же образом:

\begin{lemma}
\label{6.2.} \hfill {} \\
\hspace*{1em} Формула \ $SIN_{m}^{<\alpha _{1}}\left[
X\right] ( \alpha  ) $ \ принадлежит классу \ $\Pi _{m}$ \  для \  $%
\alpha _{1}=k$ \  и классу $ \Delta _{1}$ \ для \ $\alpha _{1}<k$.
\hspace*{\fill} $\dashv$
\end{lemma}

Следующая очевидная лемма оправдывает термин ``диссеминатор''
(распространитель), так как она показывает, что такой ординал
действительно распространяет, продолжает \ $\Pi_{m}$-свойства
(содержащие константы из его базы) с нижних уровней универсума
вплоть до релятивизирующих кардиналов:

\begin{lemma}
\label{6.3.} {\em (О продолжении)} \hfill {} \\
\hspace*{1em} Пусть
\[
    X\subseteq \alpha, \quad \alpha <\alpha_{1}, \quad
    \alpha \in SIN_{m}^{<\alpha_{1}}\left[ X\right]
\]
и утверждение\ $\forall x~\varphi ( x,\overrightarrow{a} ) $ \
имеет кортеж \ $\overrightarrow{a}$ \ констант \ $\in X$, \
$\varphi \in \Sigma _{m-1}$, \ тогда:

\vspace{6pt}
\[
    \forall x\vartriangleleft \alpha ~\varphi ^{\vartriangleleft \alpha }
    (x,\overrightarrow{a} ) \longrightarrow \forall x\vartriangleleft
    \alpha_{1}~\varphi ^{\vartriangleleft \alpha _{1}}
    ( x,\overrightarrow{a} )~.
\]
\vspace{0pt}

\noindent В этом случае мы будем говорить как и раньше, что ниже \
$\alpha _{1} $ \ ординал  $ \ \alpha $ \ расширяет или продолжает
\ $\forall x~\varphi $ \ до \ $\alpha _{1}$.
\\
Рассматривая то же самое в обращённой форме для \ $\varphi \in \Pi
_{m-1}$:

\vspace{6pt}
\[
    \exists x\vartriangleleft \alpha _{1}~\varphi ^{\vartriangleleft
    \alpha _{1}} ( x,\overrightarrow{a} ) \longrightarrow \exists
    x\vartriangleleft \alpha ~\varphi ^{\vartriangleleft \alpha }
    ( x,\overrightarrow{a} )~,
\]
\vspace{0pt}

\noindent мы будем говорить, что ниже \ $\alpha _{1}$ \ ординал \
$\alpha $ \ ограничивает или релятивизирует утверждение \ $\exists
x~\varphi $. \hspace*{\fill} $\dashv$
\end{lemma}
Свойства диссеминаторов становятся более прозрачными, если они
обладают дополнительными свойствами субнедостижимости меньших
уровней. Например, нетрудно видеть, что в этом случае класс
диссеминаторов  \ $\alpha \in SIN_{m}^{<\alpha _{1}} \left[
X\right]$ \  в \textit{широком смысле} ниже \ $\alpha _{1}$ \
уровня \ $m$ \ и с базой \ $X=\alpha $ \ совпадает с \ $
SIN_{m}^{<\alpha _{1}}$. В более общем случае справедлива

\begin{lemma}
\label{6.4.} \hfill {} \\
\hspace*{1em} Пусть
\[
    X\subseteq \alpha, \quad \alpha <\alpha _{1}, \quad
    \alpha \in SIN_{m-1}^{<\alpha _{1}},
\]
тогда следующие утвеждения равносильны:
\\
\\
{\em 1)}\quad $\alpha \in SIN_{m}^{<\alpha _{1}}\left[ X\right]$~;
\\
\\
{\em 2)}\quad для каждого \ $\mathfrak{n} \in \omega_0$ \ и
каждого кортежа \ $\overrightarrow{a}$ \ констант \ $\in X $
\[
    L_{k} \vDash \exists x\vartriangleleft
    \alpha _{1}~U_{m-1}^{\Pi \vartriangleleft \alpha _{1}}
    (\mathfrak{n}, x,\overrightarrow{a} ) \longleftrightarrow \exists
    x\vartriangleleft \alpha ~U_{m-1}^{\Pi \vartriangleleft
    \alpha _{1}} (\mathfrak{n}, x,\overrightarrow{a} )~;
\]

\noindent {\em 3)}\quad для каждой \ $\Pi_{m-1}$-формулы \ $
\varphi ( x,\overrightarrow{a} ) $ \ и каждого кортежа\
$\overrightarrow{a} $ \ констант~~$\in X$
\[
    L_{k} \vDash \exists x\vartriangleleft \alpha _{1}~
    \varphi^{\vartriangleleft \alpha _{1}} ( x,\overrightarrow{a} )
    \longleftrightarrow \exists x\vartriangleleft \alpha ~
    \varphi^{\vartriangleleft \alpha _{1}} ( x,\overrightarrow{a} )~.
\]
\end{lemma}

Эта лемма подготавливает следующее важное

\begin{definition}
\label{6.5.} \hfill {} \\
{\em 1)}\quad Минимальный диссеминатор класса
\[
    SIN_{m}^{<\alpha _{1}}\left[ X\right] \cap
    SIN_{m-1}^{<\alpha_{1}}
\]
будет называться производящим диссеминатором с базой данных \ $X$
\ ниже \ $\alpha_1$ \ и будет обозначаться общим символом \
$\check{\delta}^{X}$, \ или \ $\check{\delta}$ \ или \ $\delta$;

\noindent {\em 2)}\quad без этого условия минимальности
диссеминатор этого класса будет называться плавающим и будет
обозначаться общим символом \ $\widetilde{\delta }^{X}$, \ или,
короче, символом \ $\widetilde{\delta}$ \ или \ $\delta$, \
указывая на \ $\alpha_1$ \ и на другие его атрибуты в контексте.
\hspace*{\fill} $\dashv$
\end{definition}

Как обычно, индексы \ $m$, \ $\alpha _{1}$,$X$ \ будут опускаться,
если они произвольны либо восстанавливаются из контекста.

\noindent Эти термины ``производящий диссеминатор'' и ``плавающий
диссеминатор'' оправдываются следующим обстоятельством:
\textit{производящий диссеминатор} \textit{единственным образом}
определяется через его базу \ $X$ \ ниже \ $\alpha_1$ \ как
минимальный, в то время как \textit{плавающий диссеминатор} может
и не быть единственным и его значение не установлено (он
``плавает''); более того, термин ``производящий диссеминатор''
оправдывает также следующая лемма

\begin{lemma}
\label{6.6.} \hfill {} \\
\hspace*{1em} Пусть
\[
    X\subseteq \alpha _{0}, \quad  \alpha _{0}<\check{\delta}<\alpha _{1},
\]
\ $\check{\delta}$ \ \smallskip это производящий диссеминатор \
$\in SIN_{m}^{<\alpha _{1}}\left[ X\right] $ \ и \  $\varphi (
\alpha ,\overrightarrow{a} ) $ \ это \ $\Sigma _{m-1}$-формула
\smallskip с ординальной переменной \ $\alpha $ \ и кортежем
\  $\overrightarrow{a}$ \  \smallskip констант~~$\in X$.
\\
Предположим, что
\[
    \forall \alpha \in \left] \alpha _{0} ,\check{\delta}
    \right[~\varphi ^{\vartriangleleft \check{\delta}}
    ( \alpha ,\overrightarrow{a} ),
\]
\noindent тогда существует некоторое \ $\alpha _{0}^{\prime }\in
\left[ \alpha _{0}^{{}},\check{\delta}\right[ $ \  такое, что

\[
    \forall \alpha \in \left] \alpha _{0}^{\prime },\alpha _{1}
    \right[ ~~\varphi^{\vartriangleleft \alpha _{1}}
    ( \alpha ,\overrightarrow{a} ).
\]
\end{lemma}
\textit{Доказательство.} \ Здесь следует отметить, что переменная
\ $\alpha$ \ под квантором \ $\forall$ \ пробегает не всё \
$\check{\delta}$, \ как это полагалось в лемме~6.3, но только его
заключительный сегмент\ $\left[ \alpha _{0}^{{}},\check{\delta}
\right[ $, и тем не менее этот \ $\check{\delta}$ \
\emph{производит} рапространение \ $\varphi ( \alpha
,\overrightarrow{a} )$\ вплоть до \ $\alpha_{1}$.\
\\
Предварительно нужно сделать следующее замечание: для всяких
ординалов\ $\alpha \in SIN_{e}$, \ $e>1$, \ $\beta $: \qquad
\qquad $ \beta < \alpha \longleftrightarrow \beta \vartriangleleft
\alpha$.
\\
Теперь достаточно сначала рассмотреть случай, когда \ $\alpha_{0}
\in SIN_{m-1}^{<\alpha_{1}} $ \ (см. лемму 6.7 ниже). Из
минимальности \ $\check{\delta}$ \ следует, что \ $\alpha_{0}
\notin SIN_{m}^{<\alpha_{1}}[X] $ \ и по лемме 6.4 это означает,
что существует \ $\Sigma_{m-1}$-формула \ $\varphi_{1}(\alpha,
\overrightarrow{a_{1}}) $ \  с кортежем \ $\overrightarrow{a_{1}}
$ \  констант $\in X $ \ такая, что

\begin{equation}
\label{e6.1} \forall \alpha < \alpha_{0} \ \
\varphi_{1}^{\triangleleft \alpha_{0}}(\alpha,
\overrightarrow{a_{1}})~,
\end{equation}
\begin{equation}
\label{e6.2} \exists \alpha^{\prime} < \alpha_{1} \neg
\varphi_{1}^{\triangleleft \alpha_{1}}(\alpha^{\prime},
\overrightarrow{a_{1}})~.
\end{equation}
\vspace{0pt}

\noindent Так как  \  $\alpha_{0},\check{\delta} \in
SIN_{m-1}^{<\alpha_{1}} $, \ то утверждение (6.1) эквивалентно

\vspace{6pt}
\begin{equation}
\label{e6.3} \forall \alpha < \alpha_{0} \ \
\varphi_{1}^{\triangleleft \check{\delta}}(\alpha,
\overrightarrow{a_{1}})~.
\end{equation}
\vspace{0pt}

\noindent Кардинал  \  $\check{\delta}$ \  есть диссеминатор и
поэтому (6.2) влечёт

\vspace{6pt}
\[
    \exists \alpha^{\prime} < \check{\delta} \neg
    \varphi_{1}^{\triangleleft \check{\delta}}(\alpha^{\prime},
    \overrightarrow{a_{1}})~.
\]
\vspace{0pt}

\noindent Отсюда и из(6.3) следует, что существует ординал
 \  $\alpha_{0}^{\prime} \in [ \alpha_{0}, \check{\delta}[ $ \  такой, что

\begin{equation}
\label{e6.4} \forall \alpha < \alpha_{0}^{\prime} \ \
\varphi_{1}^{\triangleleft \check{\delta}}(\alpha,
\overrightarrow{a}_{1})\wedge \neg \varphi_{1}^{\triangleleft
\check{\delta}}(\alpha_{0}^{\prime}, \overrightarrow{a}_{1})~.
\end{equation}
\vspace{0pt}

\noindent Теперь из условия этой леммы легко заключить утверждение

\vspace{6pt}
\[
    \forall \alpha^{\prime}, \alpha < \check{\delta}   \left ( \neg
    \varphi_{1}^{\triangleleft \check{\delta}}(\alpha^{\prime},
    \overrightarrow{a}_{1}) \wedge \alpha > \alpha^{\prime}
    \longrightarrow  \varphi^{\triangleleft \check{\delta}}(\alpha,
    \overrightarrow{a}) \right )~.
\]
\vspace{0pt}

\noindent Оно содержит константы $\in X $, \ поэтому по
лемме~6.4~3) диссеминатор \ $\check{\delta} $ \ продолжает его до
\  $ \alpha_{1} $:

\begin{equation}
\label{e6.5} \forall \alpha^{\prime}, \alpha < \alpha_{1}   \left
( \neg \varphi_{1}^{\triangleleft \alpha_{1}}(\alpha^{\prime},
\overrightarrow{a}_{1}) \wedge \alpha > \alpha^{\prime}
\longrightarrow  \varphi^{\triangleleft \alpha_{1}}(\alpha,
\overrightarrow{a}) \right )~.
\end{equation}
\vspace{0pt}

\noindent Так как  \  $\check{\delta} \in SIN_{m-1}^{<\alpha_{1}}
$, \ то (6.4) влечёт  \ $\neg \varphi_{1}^{\triangleleft
\alpha_{1}}(\alpha_{0}^{\prime}, \overrightarrow{a}_{1}) $ \ и
\smallskip поэтому (6.5) влечёт  \  $\forall \alpha \in [
\alpha_{0}^{\prime}, \alpha_{1}[ ~~ \varphi^{\triangleleft
\alpha_{1}}(\alpha, \overrightarrow{a}) $.\hfill $\dashv$
\\

Используя эту лемму будем говорить, как и выше, что \
$\check{\delta} $ \ \textit{расширяет или продолжает} утверждение
\ $\forall \alpha \ \ \varphi $ \ до \ $\alpha_{1} $, \ или, в
обращённой форме, что \ $\check{\delta}$ \ \textit{ограничивает
или релятивизирует} утверждение \ $\exists \alpha \neg \varphi $ \
ниже \ $\alpha_{1} $, \ указывая ординалы \ $\alpha_{0},
\alpha_{0}^{\prime} $ \ в контексте.

\begin{corollary}
\label{6.7.} \hfill {} \\
\hspace*{1em} Пусть \ $\check{\delta}\in SIN_{m}^{<\alpha
_{1}}\left[ X\right] $ \ это производящий диссеминатор ниже  \
$\alpha_{1} $, \ тогда
\[
    \sup  ( \check{\delta}\cap SIN_{m-1}^{<\alpha _{1}} ) =\check{\delta}.\
\]
\end{corollary}
\textit{Доказательство.} \ Пусть, напротив,
\[
    \alpha_{0} = \sup  \left ( \check{\delta}\cap SIN_{m-1}^{<
    \alpha_{1}} \right ) < \check{\delta},
\]
тогда \ $\alpha_{0}\in SIN_{m-1}^{<\alpha_{1}} $ \ и можно
применить аргумент из предыдущего доказательства. Диссеминатор \
$\check{\delta}$ \ ограничивает предложение
\[
    \exists \alpha > \alpha_{0}  \ \ SIN_{m-1}(\alpha),
\]
которое выполняется ниже \ $\alpha_{1} $, \ так как сам \ $
\check{\delta}\in SIN_{m-1}^{<\alpha_{1}} $. \ А тогда существует
\[
    \alpha \in \; ] \alpha_{0}, \check{\delta}[ \; \cap SIN_{m-1}^{< \alpha_{1}}
\]
вопреки предположению. \hfill $\dashv$
\\
Теперь ясно, что в доказательстве леммы 6.6 условие $\alpha_{0}\in
SIN_{m-1}^{<\alpha_{1}} $ \ можно опустить.
\\
\quad \\
В дальнейшем мы будем исследовать взаимное расположение
диссеминаторов одного и того же уровня \ $m$ \ ниже одного и того
же \ $\alpha _{1}$, \ но с разными базами данных.  Как эти базы
влияют на их расположение? Обсуждение подобных вопросов более
прозрачно, когда\ база данных \ $X$ \ это некоторый ординал \
$\rho$ \ (то есть множество меньших ординалов). Чтобы подчеркнуть
подобные случаи мы будем писать
\[
    SIN_{m}^{<\alpha _{1}}\left[ <\rho \right]
    \mbox{\it \quad \emph{вместо} \quad} SIN_{m}^{<\alpha _{1}}
    \left[ X\right].
\]
Далее все базы \ $\rho $ \  это \textit{предельные} ординалы (если
не подразумевается другой случай).
\\
\hfill {} \\
Следующая лемма, хотя и очевидная, предоставляет, однако,
несколько важных технических инструментов диссеминаторного
анализа:

\begin{lemma}
\label{6.8.}
\quad \\
{\em 1)}\quad Пусть
\[
    \delta < \alpha _{1}, \quad \delta \in SIN_{m}^{<\alpha _{1}}\left[ X\right]
    \cap SIN_{m-1}^{<\alpha_{1}}
\]
\vspace{0pt}

\noindent тогда всякий \ $SIN_{m-1}^{<\alpha_{1}}$-кардинал из \
$[ \delta, \alpha_1 [$ \ это снова диссеминатор (плавающий) с той
же базой \ $X$ \  данных ниже \ $\alpha_1$:

\[
    \alpha \in SIN_{m-1}^{<\alpha_{1}} \wedge \alpha > \delta
    \longrightarrow \alpha \in SIN_{m}^{<\alpha_{1}}\left[
    X\right];
\]
\vspace{0pt}

\noindent {\em 2)}\quad возрастание базы данных \ $\rho$ \ влечёт
неубывание соответствующих производящих диссеминаторов ниже \
$\alpha_1$:

\[
    \rho _{1}<\rho _{2}\longrightarrow \check{\delta}^{\rho
    _{1}}\leq \check{\delta}^{\rho _{2}};\quad \rho _{1}^{+}\leq
    \rho_{2}\longrightarrow \check{\delta}^{\rho_{1}}<
    \check{\delta}^{\rho_{2}};
\]
\vspace{0pt}

\noindent {\em 3)}\quad предельный переход в базах данных влечёт
предельный переход в соответствующих производящих диссеминаторах
ниже \ $\alpha_1$:
\[
    \lim_{i} \rho _{i}=\rho \longrightarrow \lim_i \check{\delta}^{\rho
    _i}=\check{\delta}^{\rho }.
\]
\end{lemma}
\textit{Доказательство.} 1). Верхний индекс \ $\alpha_{1} $ \ для
краткости будет опускаться. Обратимся к определению~6.1 (или к
лемме~6.4) и рассмотрим утверждение
\[
    \exists x \ \
    U_{m-1}^{\Pi}( \mathfrak{n},x,\overrightarrow{a})
\]
для произвольных \ $\mathfrak{n} \in \omega_0$ \ и кортежа \
$\overrightarrow{a} $ \ констант~~$\in X $.
\\
Предположим, что
\[
    \exists x \vartriangleleft \alpha_1
    ~U_{m-1}^{\Pi}(\mathfrak{n},x,\overrightarrow{a}).
\]
Так как \ $\delta $ \ это диссеминатор, то это влечёт
\[
    \exists x \vartriangleleft \delta ~
    U_{m-1}^{\Pi}(\mathfrak{n},x,\overrightarrow{a})
\]
и, таким образом,
\[
    \exists x \vartriangleleft \alpha \ \
    U_{m-1}^{\Pi}(\mathfrak{n},x,\overrightarrow{a}) .
\]
Здесь \ $U_{m-1}^{\Pi}$ \ следует \ $\vartriangleleft$-ограничить
кардиналом  \ $\alpha $, \ так как \ $\alpha \in SIN_{m-1} $; \
поэтому \ $\alpha \in SIN_{m}[X] $. \ Утвеждение 2) следует из
определений и вместе с 1) влечёт 3). \hfill $\dashv$\\
\quad \\
Чтобы получить более детальную информацию о расположении
диссеминаторов  естественно рассмотреть их применительно к
матричным носителям.
\\
Далее через  $\widehat{ \rho }$ \  обозначается замыкание \ $\rho
$ \  относительно функции пары.

\begin{definition}
\label{6.9.} \hfill {} \\
\hspace*{1em} Пусть \ $\chi<\alpha<\alpha_1$ \ и\ $S$ \ матрица,
редуцированная к кардиналу \ $\chi$ \ и сингулярная на носителе \
$\alpha $.
\\
{\em 1)}\quad Мы будем называть диссеминатором для \ $S$ \ на \
$\alpha $ \ (или диссеминатором для этого носителя) уровня \ $m$ \
с базой данных \  $X$ \ всякий диссеминатор \ $\delta $ \ ниже
кардинала предскачка \ $\alpha_{\chi}^{\Downarrow }$ \ после \
$\chi$ \ с теми же параметрами, то есть всякий
\[
    \delta \in SIN_{m}^{<\alpha_{\chi}^{\Downarrow }}\left[ X\right]
    \cap SIN_{m-1}^{<\alpha_{\chi}^{\Downarrow }}.
\]
В таком случае мы также будем говорить, что матрица \ $S$ \ на \
$\alpha$ \ (или её носитель \ $\alpha $) \ обладает этим
диссеминатором (с этими параметрами).
\\
\quad \\
{\em 2)}\quad Мы будем называть {\sl собственным диссеминатором}
матрицы \ $S$ \ на \ $\alpha$ \ уровня \ $m$ \ всякий её
диссеминатор этого уровня с базой данных \ $\rho =
\widehat{\rho_1}$, $\rho_1 = Od(S)$, \ и будем обозначать его
общим символом \ $\widetilde{\delta}^S$ \ или \ $\delta^S$, \ а
его базу данных \ $\rho$ --- символом\ $\rho^S$; \ если\
$\widetilde{\delta}^S$ \ минимален с этой базой \ $\rho^S$, \ то
он будет называться {\sl производящим собственным диссеминатором}
для \ $S$ \ на \ $\alpha$ \ и обозначаться через \
$\check{\delta}^S$.
\\
\quad \\
{\em 3)}\quad Матрица \ $S$ \ будет называться диссеминаторной
матрицей, или, короче, \ $\delta $-матрицей уровня \ $m$, \
допустимой на носителе \ $\alpha$ для \ $\gamma =
\gamma_{\tau}^{<\alpha_1}$ \ ниже \ $\alpha_1$ \ если она обладает
некоторым диссеминатором  \ $ \delta < \gamma$ \ уровня \ $m$ \ с
некоторой базой данных \ $\rho \le \chi^{\ast +}$ \ на этом
носителе такой, что \ $S \vartriangleleft \rho, \ \rho =
\widehat{\rho }$ \ (также ниже\ $\alpha_1$). \label{c6}
\endnote{
\ стр. \pageref{c6}. \ Здесь следует обратить внимание на
следующее интересное понятие безэксцессивной базы данных:
\textit{база\ $\rho$ \ матрицы \ $S$, \ допустимой на носителе \
$\alpha$,\ называется безэксцессивной}, \ если если уменьшение \
$\rho$ \ влечёт уменьшение производящего диссеминатора для \ $S$ \
на том же носителе \ $\alpha$:
\[
    \forall \rho^{\prime} \left ( \rho^{\prime} < \rho
    \longrightarrow \check{\delta}^{\rho^{\prime}} <
    \check{\delta}^{\rho} \right ).
\]

\noindent Дело в том, что среди баз диссеминаторов матрицы \ $S$ \
на носителе \ $\alpha$ \ возможны базы \ $\rho$, \ обладающие
некоторой ``избыточностью'' информации в том смысле, что для
некоторого \ $\rho^{\prime} < \rho$ \ всё ещё
\[
    \check{\delta}^{\rho^{\prime}} = \check{\delta}^{\rho}
\]
на \ $\alpha$ \ и, выходит, вместо  \ $\rho$ \  \textit{может быть
использована меньшая база \ $\rho^{\prime} < \rho $ \ без всякой
потери} для расположения диссеминатора (и, следовательно, для его
действия). Естественно использовать базы, свободные от подобной
избыточности, то есть безэксцессивные базы. Следует отметить, что
в подобном случае использование меньшей базы, может быть,
возможно, \textit{но только} когда соответствующий диссеминатор
уменьшен.

\noindent Нетрудно доказать лемму, которая устанавливает
замечательное свойство: не только такая база определяет
расположение производящего диссеминатора, но и наоборот --  такой
диссеминатор определяет его безэксцессивную базу даже для разных
матриц на разных носителях:
\\

\noindent {\bf Лемма}
\\
\hspace *{1em} \em \noindent Пусть \ $S_1$, \ $S_2$ \ это матрицы
на носителях \ $\alpha_1$, \ $\alpha_2$, \ обладающие
производящими диссеминаторами \ $\check{\delta}_1$, \
$\check{\delta}_2$ \ уровня \ $m$ \ с безэксцессивными базами \
$\rho_1$, \ $\rho_2$ \ на этих носителях соответственно, тогда:
\[
    \mbox{\it если \quad} \check{\delta}_1 = \check{\delta}_2 \in
    SIN_{m-1}^{<\alpha_1^{\Downarrow}} \cap
    SIN_{m-1}^{<\alpha_2^{\Downarrow}}
    ,\quad \mbox{\it то} \quad \rho_1 = \rho_2.
\]
\em 
\\
} 
\\

\noindent {\em 4)}\quad В дальнейшем через \ $Lj^{<\alpha } ( \chi
) $ \ обозначается следующая \ $\Delta _{1}$-формула:
\[
    \chi <\alpha \wedge SIN_{n-1}^{<\alpha } ( \chi  ) \wedge
    \sum rng \left ( \widetilde{\mathbf{S}}_{n}^{\sin
    \vartriangleleft \chi } \right ) \in B_{\chi }\wedge
    \sup dom \left (\widetilde{\mathbf{S}}_{n}^{\sin
    \vartriangleleft \chi } \right )=\chi,
\]
и кардинал \ $\chi$ \ с этим свойством называется насыщенным ниже\
$\alpha$.
\\
Теперь примем следующее важное соглашение:
\\
везде далее мы будем рассматривать сингулярные матрицы  \ $S$ \ на
носителях \ \mbox{$\alpha \in SIN_{n-2}$}, \ редуцированные к
некоторым насыщенным кардиналам \ $\chi$,\ обладающих свойством:
\\
\hspace*{9em} $\chi^\ast \le \chi \wedge Lj^{<\alpha } ( \chi )$.
\hspace*{\fill} $\dashv$
\end{definition}

Следует отметить, что кардиналы этого вида существуют, например,
кардинал \ $\chi = \chi ^{\ast }$ \ для каждого \ $\alpha >\chi
^{\ast }$, $\alpha \in SIN_{n-2}$, \ как это можно видеть из
лемм~\ref{5.5.}, \ref{3.8.} (последняя для \ $n-1$ \ в роли \
$n$).
\\
Это соглашение делает возможным применение леммы~\ref{4.6.} о
спектральном типе к таким кардиналам.
\\
 \textit{Зафиксируем такой кардинал} \ $ \chi$ \ вплоть до
специального замечания, так что нижний индекс \ $\chi $ \ будет
как правило опускаться; например, кардинал скачка \
$\alpha_{\chi}^{\Downarrow} $ \ будет обозначаться просто через \
$\alpha^{\Downarrow} $ \ и так далее.
\\

\noindent Следующая лемма анализирует множество матричных
диссеминаторов с базами, превосходящими матричный тип самой
матрицы; её доказательство демонстрирует типичное применение леммы
4.6 о спектральном типе:

\begin{lemma}
\label{6.10.} \hfill {} \\
\hspace*{1em} Пусть матрица \ $S$ \ на носителе \ $\alpha$ \
редуцирована к кардиналу \ $\chi$ \ и обладает некоторым
диссеминатором уровня
 \ $m$ \ с
базой \ $\in \left[ OT ( S ) ,\chi ^{+}\right[ $\;.
\\
Тогда \ $S$ \ обладает диссеминаторами этого уровня со всеми
базами $<\chi ^{+}$ \ на этом носителе; аналогично для
производящих и плавающих диссеминаторов.
\end{lemma}

\noindent \textit{Доказательство.} \ Допустим,  это неверно, тогда
существует ординал
\[
    \rho_{0} \in [ OT(S), \chi^{+} [
\]
такой, что для всякого   \  $\rho \in [\rho_{0}, \chi^{+} [ $ \
матрица  \  $S $ \  на \ $\alpha$ \ не обладает ни одним
диссеминатором с базой \ $\rho $. \ В этом случае стоя на \
$\alpha^{\downarrow} $ \ можно обнаружить ординал \ $\rho_{1} \in
[ OT(S), \chi^{+} [ $, \ который может быть определён через  \
$\chi $ \ ниже \ $\alpha^{\Downarrow} $ \ формулой:
\[
    \exists \delta < \alpha^{\Downarrow} \ SIN_{m}^{<
    \alpha^{\Downarrow} } [<\rho_{1}] (\delta) \wedge \forall \rho
    \in \;] \rho_{1}, \chi^{+} [ \ \forall \delta < \alpha^{\Downarrow}
    \neg SIN_{m}^{< \alpha^{\Downarrow} } [< \rho] (\delta)~.
\]

\noindent Но тогда по лемме 4.6 получается противоречие: \ $OT(S)
\leq \rho_{1} < OT(S) $. \hfill $\dashv$
\\

\noindent Аналогичным рассуждением просто устанавливается

\begin{lemma}
\label{6.11.} \hfill {} \\
\hspace*{1em} Пусть матрица \ $S$ \ на носителе \ $\alpha$ \
редуцирована к кардиналу \ $\chi$ \ и пусть она обладает
производящими диссеминаторами \ $ \check{\delta}^{\rho _{0}}$, \
$\check{\delta}^{\rho _{1}}$ \ уровня \ $m$ \ с базами \
$\rho_{0}$, \ $\rho_{1}$ \ соответственно на этом носителе такими,
что

\[
    \rho_{0} < OT(S) \leq \rho_{1} <\chi ^{+}.
\]

\noindent Тогда \qquad \qquad\qquad
$\check{\delta}^{\rho_{0}}<\check{\delta}^{\rho _{1}}$ \mbox{\it
\quad или \quad} $\check{\delta}^{\rho _{0}}= \check{\delta}^{\rho
_{1}}=\check{\delta}^{\chi ^{+}}$. \hspace*{\fill} $\dashv$

\end{lemma}

\noindent Опишем некоторые способы построения \ $\delta $-матриц,
которые используются в дальнейшем. Следует обратить самое
пристальное внимание на следующее важное условие
\[
    A_n^{\vartriangleleft \alpha_1}(\chi^{\ast}) =
    \left\| u_{n}^{\vartriangleleft \alpha _{1}}
    ( \underline{l}) \right\|~,
\]
где, напомним, формула  \ $u_n(\underline{l})$ \ -- это
спектрально универсальная \ $\Sigma_n$-формула \ $u_{n}^{\Sigma }$
\ без индивидных констант, но с вхождениями  \ $\underline{l}$ \
(см. определение 2.6); это условие играет существенную роль в
дальнейшем и означает, что никакое \ $\Sigma_n$-утверждение \
$\varphi(\underline{l})$ \ не имеет ординалов скачка после \
$\chi^{\ast}$ \ ниже \ $\alpha_1$.

\begin{lemma}
\label{6.12.} {\em (О порождении сингулярных матриц)}
\\
\hspace*{1em} Пусть кардиналы \ $\chi <\alpha _{0}<\alpha _{1}$ \
выполняют условия:
\\

(i)\quad $ \alpha _{0}\in SIN_{n-1}^{<\alpha _{1}}$~; \medskip

(ii)\quad $ A_{n}^{\vartriangleleft \alpha _{1}} ( \chi^{\ast}  )
=\left\| u_{n}^{\vartriangleleft \alpha _{1}} ( \underline{l} )
\right\| $~.
\\
\quad \\
Тогда существует матрица \ $S_{0}$, \ которая редуцирована к
кардиналу \ $\chi$ \ и сингулярна на носителе \ $\alpha
_{0}^{\prime }\in \left] \alpha _{0},\alpha _{1} \right[ $ \
таком, что \ $\alpha _{0}=\alpha _{0\chi }^{\prime \Downarrow }$.
\\
В этом случае матрица \  $S_{0}$ \ называется порождённой
кардиналом \ $\alpha _{0}$ \ ниже \ $\alpha _{1}$.
\end{lemma}
\textit{Доказательство.} \ Применим лемму 5.12. По этой лемме
выполняется \
$\Pi_n$-утверждение \ $\varphi$ 

\[
    \forall \chi^\prime \forall \alpha _{0}^\prime \bigl(  ( \chi^\prime
    \mbox{\it \ предельный кардинал }>\omega _{0} ) \longrightarrow
    \exists \alpha _{1}^\prime > \alpha _{0}^\prime ~
    \sigma  ( \chi^\prime ,\alpha _{1}^\prime )  \bigr) ;
\]
\vspace{0pt}

\noindent  следовательно, оно выполняется ниже каждого \
$\Pi_{n-1}$-кардинала по лемме 3.2 для \ $n-1$ \ в роли \ $n$); \
но более того  --  оно выполняется также ниже любого
$\Pi_{n-2}$-кардинала \ $\alpha_1$, \ выполняющего условия этой
леммы.
\\
Этот факт обеспечивается условием $(ii)$, которое влечёт
сохранение этого утверждения \ $\varphi$ \ ниже \ $\alpha_1$. \
Действительно, нетрудно видеть, что в противном случае \
$\Sigma_n$-утвеждение \ $\neg \varphi$ \ получает ординал скачка
\textit{после \ $\chi^\ast$ \ ниже \ $\alpha_1$} -- именно,
\textit{минимальную} пару \ $(\chi^\prime, \alpha_0^\prime)$ \
такую, что утверждение

\[
    \bigl(  ( \chi^\prime
    \mbox{\it \ предельный кардинал }>\omega _{0} ) \ \wedge \
    \forall \alpha _{1}^\prime > \alpha _{0}^\prime ~
    \neg \sigma  ( \chi^\prime ,\alpha _{1}^\prime )  \bigr)
\]
\vspace{0pt}

\noindent выполняется ниже \ $\alpha_1$; \ поэтому \
$\Sigma_n$-спектрально универсальное утверждение \
$u_n(\underline{l})$ \ также получает этот ординал скачка
благодаря лемме~\ref{2.7.}. Это влечёт следующее нарушение условия
\ $(ii)$:
\[
    A_n^{\vartriangleleft \alpha _{1}} ( \chi^{\ast} ) < \|
    u_n^{\vartriangleleft \alpha _{1} } (\underline{l}) \|~.
\]
Итак, \ $\varphi^{< \alpha_1}$ \ выполняется; отсюда следует
существование кардинала

\[
    \alpha _{0}^{\prime }=\min \{ \alpha \in
    \left] \alpha _{0},\alpha _{1}\right[ :\exists S~\sigma
    ( \chi,\alpha ,S ) \}
\]
\vspace{0pt}

\noindent и матрицы \  $S_{0}$ \  на носителе \ $\alpha
_{0}^{\prime }$; \ эта минимальность обеспечивает \ \mbox{$\alpha
_{0}=\alpha _{0\chi }^{\prime \Downarrow }$}. \label{c7}
\endnote{
\ стр. \pageref{c7}. \ Здесь совершенно необходимо обратить
внимание на тот факт, что условие $(ii)$ этой леммы
\[
A_{n}^{\vartriangleleft \alpha _{1}} ( \chi^{\ast}  ) =\left\|
u_{n}^{\vartriangleleft \alpha _{1}} ( \underline{l} ) \right\|
\]
является достаточным для проведения и этого доказательства, и
также для всех дальнейших рассуждений. Чтобы проиллюстрировать, к
чему может привести непонимание этого условия, автор хотел бы
представить (в качестве примера) отзыв одного рецензента в 2000
году:
\\
\noindent {\em Рецензент}:  ``Лемма 6.12 неверна. Вот контрпример.
Пусть \ $\chi$ \ это ординал \ $\chi^{\ast}$, \ определённый на
стр.25 [здесь определение \ref{5.4.} на стр. \pageref{5.4.} ---
\emph{автор}]. Тогда для каждого \ $\alpha_1$, \ условие $(ii)$
леммы 6.12 выполняется. Пусть \ $\alpha_1$ \ это
\textit{наименьший} ординал \ $\alpha$, \ обладающий следующими
тремя свойствами:
\\
\\
(a) \ $\alpha > \chi$;
\\
(b) \ $\alpha \in SIN_{n-2}$;
\\
(c) \ $\sup \left\{ \beta < \alpha: \beta \in SIN_{n-1}^{< \alpha}
\right\} = \alpha$.
\\
\\
Очевидно, такой \ $\alpha$ \ существует, и всякий такой выполняет
условие \ $\alpha$ \ $(i)$ леммы 6.12. Пусть \ $\alpha_0$ \ это
любой элемент \ $SIN_{n-1}^{<\alpha_1}$, \ превосходящий \ $\chi$.
\ Тогда посылка 6.12 удовлетворена. Но заключение 6.12 нарушено.
Так как если \ $S_0$ \ и \ $\alpha_0^{\prime}$ \ выполняют
заключение, то по определению ``сингулярной матрицы'' и
``носителя'', \ $\alpha_0^{\prime}$ \ выполняет (a)-(c); но\
$\alpha_0^{\prime} < \alpha_1$ \ в противоречии с минимальностью \
$\alpha_1$.
\\
\\
Ошибка автора в доказательстве леммы 6.12 состоит в использовании
леммы 5.12. Действительно, 5.12 выполняется в модели 5.12 \ $L_k$,
\ но не выполняется в модели \ $L_{\alpha_1}$. \ Переформулируя
это в обозначениях автора, если \ $\varphi$ \ это утверждение из
5.12, то \ $\varphi$ выполняется, а \ $\varphi^{<\alpha_1}$ \ не
выполняется. \ $\varphi$ --- как автор указывает в доказательстве
5.12 --- это \ $\Pi_n$\ -утверждение, и \ $\alpha_1$ \ только
уровня \ $n-2$, \ а не уровня \ $n-1$, субнедостижимый; поэтому не
существует обоснования для заключения, что\ $\varphi$ \ может быть
ограничено к \ $\alpha_1$. \ Контрпример показывает, что это
невозможно.''
\\
\\
{\em Комментарий автора в его ответе в 2000 году}: ``Нет, этот
контрпример сам неверен: не для всякого \ $\alpha_1$ \ условие
$(ii)$ выполняется; например, для \ $\alpha_1$, \ использованного
самим рецензентом! Действительно, рецензент прав, что в этом
случае утверждение \ $\varphi$ из 5.12, использованное им,
выполняется в \ $L_k$ \ (и поэтому сохраняется ниже\
$\chi^{\ast}$) \ и не выполняется ниже \ $\alpha_1$, \ но отсюда
он делает ложный вывод. Наоборот, противоположная ситуация имеет
место: утверждение\ $\neg \varphi$ \ получает ординал скачка после
 \ $\chi^{\ast}$ \ ниже \
$\alpha_1$; поэтому \ $\Sigma_n$-спектрально универсальное
утверждение \ $u_n(\underline{l})$ \ также получает такой ординал
скачка (более того, бесконечно много таких ординалов в условиях
этого контрпримера). Это означает нарушение условия (ii), то есть
в действительности
\[
    A_{n}^{\vartriangleleft \alpha_{1}} ( \chi^{\ast}  ) < \left\|
    u_{n}^{\vartriangleleft \alpha_{1}} ( \underline{l} ) \right\|.
\]
Да, действительно, \ ``$\varphi$ \ это \ $\Pi_n$-утверждение, и \
$\alpha_1$ \ только уровня \ $n-2$, \ а не уровня \ $n-1$'', но
тем не менее \ $\varphi$ сохраняется ниже \ $\alpha_1$, \ если
условие $(ii)$ выполняется''.
\\
\quad \\
} 
\hfill $\dashv$ 

\begin{lemma}
\label{6.13.} {\em (О порождении $\delta $-матриц).}
\\
\hspace*{1em} Пусть кардиналы \ $\chi <\alpha _{0}<\alpha _{1}$ \
выполняют условия:

\begin{itemize}
\item[(i)] $\alpha _{0}\in SIN_{n-1}^{<\alpha _{1}}$~; \medskip

\item[(ii)] $A_{n}^{\vartriangleleft \alpha _{1}} ( \chi^{\ast} )
=\left\| u_{n}^{\vartriangleleft \alpha _{1}} ( \underline{l} )
\right\| $~; \medskip

\item[(iii)] $\alpha _{0}= \sup (\alpha_0 \cap \sup
SIN_{m-1}^{<\alpha _{1}})$~; \medskip

\item[(iv)] на заключительном сегменте \ $\alpha _{0}$ \ некоторая формула  \ $\psi
^{\vartriangleleft \alpha _{0}} ( \beta ,\gamma) $ \ определяет
функцию \ $f ( \beta  ) =\gamma <\chi ^{+}$, \ которая неубывает
до \ $\chi ^{+}$.
\end{itemize}
Тогда существует матрица \  $S_{0}$, \ порождённая кардиналом \
$\alpha _{0} $ \ ниже \ $\alpha _{1}$, \ редуцированная к
кардиналу \ $\chi$ \ и  сингулярная на носителе \
$\alpha_0^\prime$, \ которая обладает всеми порождающими
диссеминаторами на этом носителе
\[
    \check{\delta}^{\rho }\in SIN_{m-1}^{<\alpha_{0}}\quad,
    \quad \check{\delta}^{\rho} < \alpha_{0}
\]
уровня \  $m$ \  со всеми базами \ $\rho <\chi ^{+}$.
\\
В этом случае мы будем говорить, что \ $\delta $-матрица \ $S_{0}$
\ порождена кардиналом \ $\alpha _{0}$ \ ниже \ $\alpha _{1}$.
\\
Если кроме того \  $\psi \in \Sigma _{m}$, \ то множество этих
диссеминаторов конфинально кардиналу \ $\alpha _{0}$.
\end{lemma}
\textit{Доказательство.} \ Напомним, что по лемме~\ref{6.12.} \
$\alpha_0$ \ это ординал предскачка \ $\alpha_0^{\prime
\Downarrow}$; \ согласно лемме~\ref{6.10.} достаточно установить
существование диссеминатора \ $\delta ^{\rho _{0}}\in
SIN_{m-1}^{<\alpha _{0}}$ \ с базой $\rho _{0}=OT ( S_{0} ) $. \
Из \textit{(iii)}, \textit{\ (iv)} следует, что существует
кардинал
\[
    \beta ^{S}\in SIN_{m-1}^{< \alpha_{0}},~~ \beta^S > \chi,
    \mbox{\it \quad для которого \quad} \gamma^{S}=f ( \beta ^{S} )
    >\rho _{0};
\]
теперь из лемм ~\ref{6.4.}, \ref{4.6.} следует, что кардинал
\[
    \beta^{S}\in SIN_{m}^{<\alpha _{0}} \left[ <\rho _{0}\right].
\]
является требуемым диссеминатором.
\\
Действительно, допустим, что это неверно, тогда существует
некоторая  \ $\Sigma_{m-1}$-формула \  $\varphi(\alpha,
\overrightarrow{a}) $ \ с кортежем  \ $\overrightarrow{a} $ \
констант \ $< \rho_{0} $ \ и такая, что

\vspace{6pt}
\begin{equation}
\label{e6.6} \forall \alpha < \beta^{S} \varphi^{\triangleleft
\alpha_{0}} (\alpha, \overrightarrow{a}) \wedge \exists \alpha
\in [\beta^{S}, \alpha_{0} [ \ \ \neg \varphi^{\triangleleft
\alpha_{0}} (\alpha, \overrightarrow{a}) ~.
\end{equation}
\vspace{0pt}

\noindent Для некоторой краткости возьмём кортеж \
$\overrightarrow{a} $, \ состоящий только из одного \
$\alpha_{1}^{\prime} < \rho_{0} $. \ Следует воспользоваться
преимуществом кардинала скачка
\[
    \alpha_{2} \in dom \left ( \widetilde{\mathbf{S}}_{n}^{\sin
    \triangleleft \alpha_{0}^{\prime}} \overline{\overline{\lceil}}
    \chi \right ),
\]
 -- где, напомним, \ $\alpha_{0}^{\prime}$ \ это минимальный носитель
$
> \alpha_{0} $ \ of \ $S_{0} $, --  такого, что

\vspace{6pt}
\begin{equation}
\label{e6.7} OT  \left ( \alpha_{2} \cap dom  \left (
\widetilde{\mathbf{S}}_{n}^{\sin \triangleleft
\alpha_{0}^{\prime}} \overline{\overline{\lceil}} \chi \right )
\right ) = \alpha_{1}^{\prime} ~;
\end{equation}
\vspace{0pt}

\noindent существование такого  \ $\alpha_{2} $ \ очевидно из
определения 5.7. После этого, используя (6.6), (6.7), можно
определить некоторый ординал \ $\gamma \in [ \gamma^{S}, \chi^{+}
[ $ \ через \ $\chi, \alpha_{0}, \alpha_{2} $ \  следующей
формулой
\[
    \exists \alpha^{\prime}, \beta < \alpha_{0} \ ( OT
    \left ( dom \left ( \widetilde{\mathbf{S}}_{n}^{\sin
    \triangleleft \alpha_{2}} \overline{\overline{\lceil}}
    \chi \right ) \right ) = \alpha^{\prime} \wedge \forall
    \alpha < \beta \ \ \varphi^{\triangleleft \alpha_{0}}
    (\alpha,\alpha^{\prime}) \wedge
\]
\[
    \qquad\qquad\qquad\qquad\qquad\qquad\qquad\qquad
    \wedge \neg \varphi^{\triangleleft \alpha_{0}} (\beta,
    \alpha^{\prime} ) \wedge  \psi^{\triangleleft \alpha_{0}}
    (\beta,\gamma) ) ~.
\]
\vspace{0pt}

\noindent Но тогда лемма 4.6 влечёт противоречие:

\[
    OT(S_{0}) < \gamma < OT(S_{0}) ~.
\]
\vspace{0pt}

\noindent Итак, \ $\beta^S$ \ это диссеминатор с базой  \ $\rho_0$
\ и поэтому \ $S_0$ \ обладает производящими диссеминаторами со
всеми базами
 \ $\rho <\chi^+$.
\\
Теперь допустим, что \ $\psi \in \Sigma_{m} $, \ но напротив --
множество всех диссеминаторов не конфинально кардиналу  \
$\alpha_{0} $. \ Согласно леммам~6.10, 6.8~3) в этом случае
существует диссеминатор

\[
    \beta_{1} \in SIN_{m}^{< \alpha_{0}} [< \chi^{+}]
\]
\vspace{0pt}

\noindent с базой  \ $\chi^{+} $. \ Для  \
$\gamma_{1}=f(\beta_{1}) $ \ кардинал \ $\beta_{1} $ \ продолжает
до  \  $\alpha_{0} $ \  \ $\Pi_{m}$-утверждение

\[
    \forall \beta, \gamma \quad (\psi(\beta, \gamma)
    \longrightarrow \gamma < \gamma_{1} )
\]
\vspace{0pt}

\noindent и функция  \  $f $ \  становится ограниченной  ординалом
\ $\gamma_{1} < \chi^{+} $ \ вопреки \ $(iv) $.\hfill $\dashv$
\\
\quad \\
Сходными рассуждениями устанавливается следующая важная

\begin{lemma}
\label{6.14.} {\em (О порождении матричных \ $\delta$-функций)}
\\
\hspace*{1em} Пусть кардиналы \ $\chi <\alpha _{1}$ \ выполняют
условия: \medskip
\begin{itemize}
\item[(i)] $A_{n}^{\vartriangleleft \alpha _{1}} ( \chi^{\ast} ) =\left\|
u_{n}^{\vartriangleleft \alpha _{1}} ( \underline{l} )
\right\|$~; \medskip

\item[(ii)] $\alpha _{1}=\sup SIN_{m-1}^{<\alpha _{1}}=\sup
SIN_{n-1}^{<\alpha _{1}}$~; \medskip

\item[(iii)] $cf ( \alpha_{1} ) >\chi ^{+}$~; \medskip

\item[(iv)] кардинал \ $\chi \ge \chi^{\ast}$ \ определяется через \
$\chi^\ast$ \ формулой класса \ $\Sigma_{n-2} \cup \Pi_{n-2}$~.
\end{itemize}
Тогда существует кардинал \ $\gamma _{\tau _{0}}^{<\alpha _{1}}$ \
такой, что на множестве
\[
    T=\{ \tau : \gamma _{\tau_{0}}^{<\alpha _{1}}
    \leq \gamma _{\tau }^{<\alpha _{1}}\}
\]
определена функция  \  $\mathfrak{A}$, \ которая  \quad для
каждого \ $\tau \in T$ \ обладает свойствами:
\\

{\em 1)} \quad\ $\mathfrak{A} ( \tau  ) $ \ это \ $\delta
$-матрица \ $S$ \ уровня \ $m$ \ редуцированная к \ $\chi$ \ и
допустимая на некотором носителе \ $\alpha \in \left] \gamma
_{\tau }^{<\alpha _{1}}, \gamma _{\tau + 1}^{<\alpha _{1}} \right[
$ \ для \ $\gamma _{\tau }^{<\alpha _{1}}$ \ ниже \ $\alpha_1$;
\\

{\em 2)}\quad множество всех порождающих диссеминаторов матрицы \
$S$ \ на \ $\alpha$ \ этого уровня конфинально кардиналу \ $\alpha
_{\chi }^{\Downarrow }$;
\\

{\em 3)}\quad $cf ( \alpha _{\chi }^{\Downarrow } ) =\chi ^{+}$;
\\

{\em 4)}\quad $( \gamma _{\tau }^{<\alpha _{1}}+1 ) \cap
SIN_{m-1}^{<\alpha _{1}}\subseteq SIN_{m-1}^{<\alpha _{\chi
}^{\Downarrow }}$~;
\\
\quad \\
В этом случае мы будем говорить, что функция \ $\mathfrak{A}$ \
порождена кардиналом \ $ \alpha _{1}$ \ и будем называть её
матричной \ $\delta$-функцией.
\end{lemma}
\textit{Доказательство} состоит в повторениии построения функции \
$S_{f}^{< \alpha_{1}} $ \ ниже \ $ \alpha _{1}$ \ из определения
\ref{5.14.}, но немного изменённого. Верхний индекс\ $<\alpha_{1}
$ \ будет опускаться для краткости. Применим лемму~6.12 к всякому
кардиналу \ $\gamma_{\tau_{1}}
> \chi $, \ $\gamma_{\tau_{1}} \in SIN_{m-1} $ \  конфинальности
 \ $\neq \chi^{+} $; \ существование таких кардиналов следует из
условия \ $(iii) $ \ и леммы~3.4~1) (для \ $n=m-1$). Очевидно,
матрица \ $S $ \ на носителе \ $\alpha \in \; ]
\gamma_{\tau_{1}},\gamma_{\tau_{1+1}} [ $, \ порождённая
кардиналом \ $\gamma_{\tau_{1}} $, \ обладает следующим свойством
для \ $\tau = \tau_{1} $:
\begin{equation}
\label{e6.8}cf(\alpha^{\Downarrow}) \neq \chi^{+} \wedge \forall
\gamma \leq \gamma_{\tau}  \left ( SIN_{m-1}(\gamma)
\longrightarrow SIN_{m-1}^{< \alpha^{\Downarrow}} (\gamma) \right
)~.
\end{equation}

По лемме~3.2 об ограничении (для  \ $n-1 $ \ вместо \ $n $) такие
носители \ $\alpha $ \ существуют в каждом  \  $ [
\gamma_{\tau},\gamma_{\tau_{+1}} [ $ \ для
 \  $\chi < \gamma_{\tau_{+1}} \leq  \gamma_{\tau_{1}}
$ \  (см. доказательство леммы~5.17~2)~(ii)~). Действительно,
рассмотрим \ $\Pi_{n-2}$-формулу\ $\varphi(\alpha, \chi,
\gamma^{m}_{\tau}, \gamma_{\tau}, S)$:

\vspace{6pt}
\[
    \gamma_{\tau} < \alpha  \wedge cf(\alpha^{\Downarrow})
    \neq \chi^{+} \wedge \sigma (\chi, \alpha, S) \wedge
    SIN_{m-1}^{< \alpha^{\Downarrow}} (\gamma_{\tau}^{m}),
\]

\noindent где

\[
    \gamma_{\tau}^{m} = \sup  \left ( (\gamma_{\tau}+1)
    \cap SIN_{m-1} \right ).
\]
\vspace{0pt}

\noindent  \ $ \Sigma_{n-1} $-утверждение \ $ \exists \alpha \ \
\varphi $ \  содержит константы
\[
    \chi, \quad \gamma_{\tau}^{m}, \quad \gamma_{\tau} < \gamma_{\tau+1},
    \quad S \vartriangleleft \gamma_{\tau+1},
\]
поэтому  \ $SIN_{n-1} $-кардинал  \  $ \gamma_{\tau+1} $ \
ограничивает это утверждение и матрица \ $S $ \ появляется на
носителе \  $ \alpha \in [ \chi_{\tau}, \chi_{\tau+1} [ $, \ что и
устанавливает (6.8).
\\
После этого повторим определение 5.14 матричной функции \ $ S_{f}
$ \ ниже \ $\alpha_{1} $, \ но накладывая дополнительное
требование (6.8) на носитель \ $\alpha $; \ в результате
получается функция \ $S_{f}^{m} = (S_{\tau}^{m})_{\tau} $,\
принимающая значения:

\vspace{6pt}
\[
    S_{\tau}^{m}=\min_{\underline{\lessdot }}\{ S:\exists
    \alpha < \alpha_{1}(\gamma_{\tau}<\alpha \wedge
    cf(\alpha^{\Downarrow}) \neq \chi^{+} \wedge
    \sigma ( \chi ,\alpha ,S ) \wedge
\]

\vspace{-6pt}

\[
    \qquad \qquad \qquad \qquad \qquad
    \wedge \forall \gamma \leq \gamma_{\tau} \ \
    (SIN_{m-1}(\gamma) \longrightarrow SIN_{m-1}^{<
    \alpha^{\Downarrow}}(\gamma))) \} ~.
\]
\vspace{0pt}

\noindent Эта функция  \  $\underline{\lessdot } $-неубывает как и
функция  \ $S_{f} $ \ в лемме 5.17, поэтому благодаря
\textit{(iii)} существует ординал \ $\tau^{m} \in dom \left (
S_{f}^{m} \right )$ \ такой, что
\[
    \gamma_{\tau^{m}} \in SIN_{m-1} \mbox{\it \quad и для всякого \quad}
    \tau \geq \tau^{m} \quad S_{\tau}^{m} = S_{\tau^{m}}^{m}~.
\]
Рассмотрим произвольный  \  $\tau
> \tau^{m} $ \  и носитель  \  $\alpha_{\tau}^{m} \in [
\gamma_{\tau}, \gamma_{\tau+1}[ $ \  матрицы  \ $S_{\tau}^{m} $. \
Теперь определение функции  \ $S_{f}^{m} $ \ следует повторить
ниже  \ $\alpha_{\tau}^{1} = \alpha_{\tau}^{m \Downarrow} $, \
тогда получающаяся фунция \ $S_{f}^{m<\alpha_{\tau}^{1}} $ \  по
прежнему \ $\underline{\lessdot } $-неубывает и совпадает с \
$S_{f}^{m} $ \ на  \ $\tau $. \ Рассуждение из конца предыдущего
параграфа показывает с помощью условия  $(iv)$, что она \
$\underline{\lessdot } $-неубывает до\ $\chi^{+} $ \ и поэтому
появляется кардинал

\[
    \alpha_{0 \tau} = \sup \left \{ \gamma_{\tau^{\prime}}^{<
    \alpha_{\tau}^{1}} : \tau^{\prime} \in dom \left (
    S_{f}^{m < \alpha_{\tau}^{1}} \right ) \right \}
\]
такой, что
\[
    \alpha_{0 \tau} = \sup  \left ( \alpha_{0 \tau} \cap
    SIN_{m-1}^{< \alpha_{\tau}^{1}} \right ).
\]
\vspace{0pt}

\noindent Очевидно, что \ $ cf(\alpha_{0 \tau}) = \chi^{+} $, \
поэтому условие \ $ cf(\alpha_{\tau}^{1}) \neq \chi^{+} $ \ влечёт
\[
    \gamma_{\tau}< \alpha_{0 \tau} < \alpha_{\tau}^{1}.
\]
Теперь можно применить лемму~6.13 к  \ $\alpha_{0 \tau}$, \
$\alpha_{\tau}^{1} $ \ в роли \ $\alpha_{0}$, \ $\alpha_{1} $ \
соответственно и к функции
\[
    f(\beta)= \sup  \left \{ Od  \left ( S_{\tau^{\prime}}^{m <
    \alpha_{\tau}^{1}} \right ) : \gamma_{\tau^{\prime}}^{<
    \alpha_{\tau}^{1}} < \beta \right \} ~,
\]
которые выполняют все условия этой леммы.

\noindent Матрица \ $S = S_{0 \tau} $, \ порождённая кардиналом \
$\alpha_{0 \tau} $ \ на носителе\ $\alpha_{0 \tau}^{\prime} $ \
ниже  \ $\alpha_{\tau}^{1} $ \ обладает производящими
диссеминаторами уровня \ $m $ \ со всеми базами \  $< \chi^{+} $.
\ Нетрудно видеть, что множество этих диссеминаторов конфинально \
$\alpha_{0 \tau} = \alpha_{0 \tau}^{\prime \Downarrow} $ \ и что
она имеет некоторый производящий диссеминатор \ $\check{\delta} <
\gamma_{\tau} $, \ допустимый для \ $\gamma_\tau$.
\\
Следовательно, получающаяся функция  \  $\mathfrak{A}$ \ со
значениями  \ $\mathfrak{A}(\tau) = S_{0 \tau} $ \ выполняет
заключения этой леммы полностью. \hspace*{\fill} $\dashv$
\\
                               \\                \hfill {} \\

Теперь базовая теория разработана достаточно для дальнейшего
анализа матричных функций, который будет предпринят в последующей
Части II этой работы. \label{c8}
\endnote{
\ стр. \pageref{c8}. \ Этот раздел может быть закончен ещё одним
комментарием на расположение диссеминаторов. В леммах~\ref{6.13.},
~\ref{6.14.} порождённые \ $\delta $-матрицы на носителях \
$\alpha $ \ имеют производящие диссеминаторы, расположенные
конфинально \ $\alpha _{\chi }^{\Downarrow }$,\  и \ $cf ( \alpha
_{\chi }^{\Downarrow } ) =\chi ^{+}$. \ Используя методы
рассуждений из доказательств лемм  6.13,  6.14  можно показать,
что для \ $\underline{\lessdot }$-минимальных \ $\delta $-матриц и
\ $m \geq n+1$ \ это неизбежно:
\\

\noindent {\bf Лемма} \\
\hspace*{1em} \em Пусть кардиналы \ $\chi <\gamma
_{\tau}^{<\alpha_{1}}$ \  и матрица \  $S $ \ выполняют условия:
\\

(i)\quad $ A_{n}^{\vartriangleleft \alpha _{1}} ( \chi^{\ast} )
=\left \| u_{n}^{\vartriangleleft \alpha _{1}} ( \underline{l} )
\right\| $ ~;
\\

(ii)\quad $S$ \ это \ $\delta $-матрица уровня \ $m \geq n+1$, \
допустимая на некотором носителе \ $\alpha \in \left] \gamma
_{\tau }^{<\alpha _{1}},\alpha _{1}\right[ $ \ для \ $\gamma
_{\tau }^{<\alpha _{1}}$ \ ниже \ $\alpha_1$;
\\

(iii)\quad $S$ \ это\ $\underline{\lessdot }$-минимальная матрица
изо всех  \ $\delta $-матриц с этим свойством.
\\
\quad \\
Тогда \ $S$ \ на носителе \  $\alpha $ \ обладает производящими
диссеминаторами уровня \ $m$ \ со всеми базами \ $\rho <\chi
^{+}$, \  расположенными конфинально\ $\alpha _{\chi }^{\Downarrow
}$, \ и \ $ cf ( \alpha _{\chi }^{\Downarrow } ) =\chi ^{+}$. \em
\\

\noindent Эта лемма не используется в дальнейшем и поэтому её
доказательство опускается.
\\
\quad \\
} 

\theendnotes

\newpage

\thispagestyle{empty}

\selectlanguage{russian}

\begin{center}
\quad \\
\quad \\
\quad \\
\quad \\
\quad \\
\quad \\
{Научное издание} \\
\quad \\
{\bf Киселев} Александр Анатольевич\\
\quad \\
{\bf НЕДОСТИЖИМОСТЬ \\
И \\
СУБНЕДОСТИЖИМОСТЬ}\\
\quad \\
{В двух частях}\\
{Часть I}\\
\quad \\
{\small Ответственный за выпуск {\em Т. Е. Янчук}} \\
\quad \\
{\footnotesize Подписано в печать 10.06.2011. Формат
60$\times$84 1/16. Бумага офсетная.} \\
{\footnotesize  Ризография. Усл. печ. л. 6,74.
Уч.-изд. л. 6,5.} \\
{\footnotesize Тираж 100 экз. Зак. 472.} \\
\quad \\
{\footnotesize Республиканское унитарное предприятие} \\
{\footnotesize <<Издательский центр Белорусского государственного университета>>} \\
{\footnotesize ЛИ \No 02330/0494361 от 16.03.2009.
\\
Ул. Красноармейская, 6, 220030, Минск.}\\
\quad \\
{\footnotesize Отпечатано с оригинала-макета заказчика} \\
{\footnotesize в республиканском унитарном предприятии} \\
{\footnotesize <<Издательский центр Белорусского государственного университета>>.} \\
{\footnotesize ЛП \No 02330/0494178 от 03.04.2009.
\\
Ул. Красноармейская, 6, 220030, Минск.}
\end{center}
\label{end}

\selectlanguage{english}

\end{document}